\documentclass[12pt]{amsart}
\usepackage[]{amscd,amsfonts,amsmath,amsxtra,amssymb}
\usepackage{graphics}
\usepackage[all]{xy}
\usepackage{hyperref}
\usepackage[french]{babel}
\usepackage[T1]{fontenc}
\newtheorem{sub}{}[section]
\newtheorem{subsub}{}[sub]
\usepackage[active]{srcltx}
\def\ov#1{\overline{#1}}

\def\coker{\mathop{\rm coker}\nolimits}
\def\Hom{\mathop{\rm Hom}\nolimits}
\def\HHom{\mathop{\mathcal Hom}\nolimits}
\def\Ext{\mathop{\rm Ext}\nolimits}
\def\EExt{\mathop{\mathcal Ext}\nolimits}

\def\Pic{\mathop{\rm Pic}\nolimits}
\def\Aut{\mathop{\rm Aut}\nolimits}
\def\AAut{\mathop{\mathcal Aut}\nolimits}

\def\End{\mathop{\rm End}\nolimits}

\def\GL{\mathop{\rm GL}\nolimits}

\def\imm{\mathop{\rm im}\nolimits}
\def\deg{\mathop{\rm deg}\nolimits}

\def\DIV#1#2{\frac{\partial #1}{\partial #2}}
\def\spec{\mathop{\rm spec}\nolimits}
\def\lra{\longrightarrow}

\def\sigg{\mathop{\hbox{$\displaystyle\sum$}}\limits}

\def\hfl#1#2{\smash{\mathop{\ \hbox to 12mm{\rightarrowfill}}
\limits^{\scriptstyle#1}_{\scriptstyle#2} \ }}
\def\hflb#1#2{\smash{\mathop{\hbox to 12mm{\leftarrowfill}}
\limits^{\scriptstyle#1}_{\scriptstyle#2}}}

\def\m#1{{\hbox{$#1$}}}

\def\ot{\otimes}
\def\og{\leavevmode\raise.3ex\hbox{$\scriptscriptstyle\langle\!\langle$}}
\def\fg{\leavevmode\raise.3ex\hbox{$\scriptscriptstyle\,\rangle\!\rangle$}}

\def\nsp{\lbrace 0\rbrace}

\def\Ssect#1#2{\pagebreak[3]\begin{sub}\label{#2}{\sc\small\small  #1}\rm\medskip}

\def\sepsec{\vskip 2.5cm}
\def\sepsub{\vskip 1.5cm}
\def\sepsubsub{\vskip 1cm}
\def\sepprop{\vskip 0.8cm}

\def\xmat#1{\[\xymatrix{#1}\]}

\def\flinc{\ar@{^{(}->}}
\def\fleq{\ar@{=}}
\def\flon{\ar@{->>}}
\def\fmaps{\ar@{|-{>}}}
\def\Nligne{\hfil\break}

\font\tte=cmbsy10
\def\paragra{{\tte \char120}}
\def\para{\paragra~\hskip -2pt}
\def\DB#1{\lbrack\lbrack#1\rbrack\rbrack}
\newcommand{\N}{{\mathbb N}}

\newcommand{\C}{{\mathbb C}}

\renewcommand{\P}{{\mathbb P}}

\newcommand{\F}{{\mathbb F}}

\renewcommand{\L}{{\mathbb L}}

\def\T{{\mathbb T}}
\newcommand{\ka}{{\mathcal A}}

\newcommand{\kc}{{\mathcal C}}
\newcommand{\kd}{{\mathcal D}}
\newcommand{\ke}{{\mathcal E}}
\newcommand{\kf}{{\mathcal F}}
\newcommand{\kg}{{\mathcal G}}

\newcommand{\ki}{{\mathcal I}}

\newcommand{\ko}{{\mathcal O}}

\newcommand{\ks}{{\mathcal S}}

\newcommand{\ku}{{\mathcal U}}
\newcommand{\kv}{{\mathcal V}}

\begin{document}

\def\refname{R\'ef\'erences}
\def\contentsname{Sommaire}
\def\proofname{D\'emonstration}
\def\abstractname{R\'esum\'e}
\def\emailaddrname{\it Adresse email}

\author{Jean--Marc Dr\'{e}zet}
\address{
Institut de Math\'ematiques de Jussieu,
Case 247,
4 place Jussieu,
F-75252 Paris, France}
\email{drezet@math.jussieu.fr}
\urladdr{http://www.math.jussieu.fr/\~{}drezet}

\begin{abstract}
The primitive curves are the multiple curves that can be locally embedded
in smooth surfaces (we will always suppose that the associated reduced curves
are smooth). These curves have been defined and studied by
C.~B\u anic\u a and O.~Forster in 1984. In 1995, D.~Bayer and D.~Eisenbud
gave a complete description of the double curves. We give here a
parametrization of primitive curves of arbitrary multiplicity. Let ${\bf
Z}_n={\rm spec}(\C[t]/(t^n))$. The curves of multiplicity $n$ are obtained by
taking an open cover $(U_i)$ of a smooth curve $C$ and by glueing
schemes of type $U_i\times{\bf Z}_n$ using automorphisms of $U_{ij}\times {\bf
Z}_n$ that leave $U_{ij}$ invariant. This leads to the study of the sheaf of
nonabelian groups $\kg_n$ of automorphisms of $C\times {\bf Z}_n$ that leave
the reduced curve invariant, and to the study of its first cohomology set. We
prove also that in most cases it is the same to extend a primitive curve $C_n$
of multiplicity $n$ to one of multiplicity $n+1$, and to extend the quasi
locally free sheaf $\kd_n$ of derivations of $C_n$ to a rank 2 vector bundle on
$C_n$.
\end{abstract}

\title[{\tiny Param\'etrisation des courbes multiples primitives}]
{Param\'etrisation des courbes multiples primitives}
\maketitle
\tableofcontents

\section{Introduction}
\label{intro}

Une {\em courbe primitive} est une vari\'et\'e $Y$ de
Cohen-Macaulay telle que la sous-vari\'et\'e r\'eduite associ\'ee \m{C=Y_{red}}
soit une courbe lisse irr\'eductible, et que tout point ferm\'e de
$Y$ poss\`ede un voisinage pouvant \^etre plong\'e dans une surface lisse. Ces
courbes ont \'et\'e d\'efinies et \'etudi\'ees par C.~B\u anic\u a et
O.~Forster dans \cite{ba_fo}. On s'int\'eresse plus particuli\`erement ici au
cas o\`u $C$, et donc $Y$, sont projectives.

Les courbes primitives les plus simples, c'est-\`a-dire les courbes de multiplicit\'e 2, ont \'et\'e pr\'ec\'edemment utilis\'ees par D.~Ferrand dans
\cite{fe}. Les courbes multiples apparaissent aussi fr\'equemment dans l'\'etude des fibr\'es vectoriels (ou des faisceaux r\'eflexifs) stables sur \m{\P_3} (cf. \cite{ha3}, \cite{ha4}, \cite{nu_tr}). Les courbes doubles abstraites (c'est-\`a-dire non plong\'ees) ont \'et\'e classifi\'ees par D.~Bayer et D.~Eisenbud dans \cite{ba_ei}. Le probl\`eme de leur d\'eformation en courbes lisses a \'et\'e abord\'e dans \cite{gon}.

Soient $P$ un point ferm\'e de $Y$, et $U$ un voisinage de $P$ pouvant
\^etre plong\'e dans une surface lisse $S$. Soit $z$ un \'el\'ement de
l'id\'eal maximal de l'anneau local \m{\ko_{S,P}} de $S$ en $P$ engendrant
l'id\'eal de $C$ dans cet anneau. Il existe alors un unique entier $n$,
ind\'ependant de $P$, tel que l'id\'eal de $Y$ dans \m{\ko_{S,P}} soit
engendr\'e par \m{(z^n)}. Cet entier $n$ s'appelle la {\em multiplicit\'e} de
$Y$. Si \m{n=2} on dit que $Y$ est une {\em courbe double}. Soit \m{\ki_C} le
faisceau d'id\'eaux de $C$ dans $Y$. Alors le faisceau conormal de $C$,
\m{L=\ki_C/\ki_C^2} est un fibr\'e en droites sur $C$, dit {\em associ\'e} \`a
$Y$. Il existe une filtration canonique
\[C=C_1\subset\cdots\subset C_n=Y\ ,\]
o\`u au voisinage de chaque point $P$ l'id\'eal de \m{C_i} dans \m{\ko_{S,P}}
est \m{(z^i)}.

Soit \m{L\in\Pic(C)}. On plonge $C$ dans le fibr\'e dual \m{L^*} vu comme un
surface lisse, au moyen de la section nulle. Alors le $n$-i\`eme voisinage
infinit\'esimal de $C$ dans \m{L^*} est une courbe primitive de multiplicit\'e
$n$ et de fibr\'e en droites associ\'e $L$. On l'appelle la {\em courbe
primitive triviale} de fibr\'e associ\'e $L$. Les autres exemples simples de
courbes primitives sont les courbes de Cohen-Macaulay qui sont plong\'ees dans
une surface lisse, et dont la courbe r\'eduite associ\'ee est lisse. Mais il en
existe beaucoup d'autres.

Dans le cas des courbes doubles, D.~Bayer et D.~Eisenbud ont obtenu dans
\cite{ba_ei} la classification suivante : si $Y$ est de multiplicit\'e 2, on a une suite exacte de fibr\'es vectoriels sur $C$
\begin{equation}\label{ecs}
0\lra L\lra\Omega_{Y\mid C}\lra\omega_C\lra 0\end{equation}
qui est scind\'ee si et seulement si $Y$ est la courbe triviale. En
particulier si $C$ n'est pas projective, $Y$ est toujours triviale. Dans le cas
o\`u $C$ est projective et $Y$ non triviale, cette courbe est enti\`erement
d\'etermin\'ee par la droite de \m{\Ext^1_{\ko_C}(\omega_C,L)} induite par la
suite exacte pr\'ec\'edente. Les courbes primitives non triviales de
multiplicit\'e 2 et de fibr\'e en droites associ\'e $L$ sont donc
param\'etr\'ees par l'espace projectif \m{\P(\Ext^1_{\ko_C}(\omega_C,L))}.

Le but du pr\'esent article est de donner une description analogue des courbes
primitives de multiplicit\'e quelconque, et d'\'etendre certains r\'esultats de
\cite{ba_ei}. On cherche \`a d\'ecrire les classes d'\'equivalence de courbes
primitives de multiplicit\'e $n$ de courbe r\'eduite associ\'ee $C$, deux
telles courbes \m{C_n}, \m{C'_n} \'etant dites \'equivalentes s'il existe un
isomorphisme \m{C_n\simeq C'_n} induisant l'identit\'e sur $C$.

\sepsub

\Ssect{Construction des courbes primitives}{intro_1}

\begin{subsub}Courbes multiples primitives abstraites -- \rm
En g\'en\'eral il n'existe pas de r\'etraction \m{Y\to C} (voir cette question
trait\'ee dans un cadre plus g\'en\'eral dans \cite{la}). Mais c'est vrai {\em
localement}. Pour tout ouvert $U$ de $C$ on note \m{Y(U)} l'ouvert
correspondant de $Y$. On montre (th\'eor\`eme \ref{pr0}) que pour tout point
ferm\'e $P$ de $C$ il existe un ouvert $U$ de $C$ contenant $P$ tel qu'il
existe une r\'etraction \m{Y(U)\to U}. On a alors
\[Y(U) \ \simeq \ U\times{\bf Z}_n ,\]
avec \m{{\bf Z}_n=\spec({\C[t]/(t^n)})}. La courbe $Y$ est donc obtenue en
recollant des vari\'et\'es du type \m{U\times{\bf Z}_n}. Cela conduit \`a la
notion de {\em courbe multiple primitive abstraite}. Soit \m{(U_i)} un
recouvrement ouvert \m{(U_i)} de $C$, et pour tous $i$, $j$, \m{\sigma_{ij}} un
automorphisme de \m{U_{ij}\times{\bf Z}_n} laissant \m{U_{ij}} invariant. On
suppose qu'on a la relation de cocycle :
\m{\sigma_{jk}\circ\sigma_{ij}=\sigma_{ik}} sur \m{U_{ijk}}. Le sch\'ema obtenu
en recollant les \m{U_i\times{\bf Z}_n} au moyen des \m{\sigma_{ij}} est une
courbe multiple primitive abstraite.

On montre qu'un tel sch\'ema peut \^etre plong\'e dans \m{\P_3} (th\'eor\`eme
\ref{plong2}). Les courbes multiples primitives abstraites sont donc identiques
aux courbes multiples primitives de \cite{ba_fo}.
\end{subsub}

\sepsubsub

\begin{subsub}Faisceaux de groupes d'automorphismes -- \rm Pour tout ouvert
$U$ de $C$ soit \m{\kg_n(U)} le groupe des automorphismes de \m{U\times
{\bf Z}_n} laissant $U$ invariant. On obtient ainsi un faisceau de groupes (non
ab\'eliens) \m{\kg_n} sur $C$. L'ensemble de cohomologie \m{H^1(C,\kg_n)}
s'identifie \`a celui des classes d'isomorphisme de courbes multiples
primitives de multiplicit\'e $n$ et de courbe r\'eduite associ\'ee $C$. La
construction et les propri\'et\'es de la cohomologie des faisceaux de groupes
non n\'ecessairement ab\'eliens sont rappel\'ees dans le chapitre \ref{n_ab}.

Soient \m{g\in H^1(C,\kg_n)}, et \m{C_n} la courbe de multiplicit\'e $n$
induite. De l'action de \m{\kg_n} sur lui-m\^eme par conjugaison on d\'eduit un
nouveau faisceau de groupes \m{(\kg_n)^g} (faisceau obtenu par {\em
recollement}, cf. \ref{n_ab}, \cite{fr}), qui s'identifie au faisceau de
groupes \m{\AAut_C(C_n)} des automorphismes de \m{C_n} laissant $C$ invariante.
Cette construction s'applique aussi \`a tous les sous-faisceaux de groupes
distingu\'es de \m{\kg_n} (cf. chapitre \ref{n_ab}). On note \m{\Aut_C(C_n)}
le groupe des sections globales de \m{\AAut_C(C_n)}.

La structure de \m{\kg_n} est \'etudi\'ee dans le chapitre \ref{fg_aut}. On
commence par traiter le cas de la $\C$-alg\`ebre \m{\C\DB{x,t}/(t^n)}.
Soit \m{\rho:\C\DB{x,t}/(t^n)\to\C\DB{x}} la projection. On s'int\'eresse aux
automorphismes de $\C$-alg\`ebres $\phi$ de \m{\C\DB{x,t}/(t^n)} tels que
\m{\rho\circ\phi=\rho}. On montre qu'ils sont du type \m{\phi_{\mu\nu}}, avec
\m{\mu,\nu\in\C\DB{x,t}/(t^{n-1})}, $\nu$ inversible,
\[\phi_{\mu\nu}(\alpha) \ = \ \sigg_{k=0}^{n-1}\frac{1}{k!}(\mu
t)^k\frac{d^k\alpha}{dx}\]
pour tout \m{\alpha\in\C\DB{x}}, et \m{\phi_{\mu\nu}(t)=\nu t} (th\'eor\`eme
\ref{prox1}). Les automorphismes de \m{U\times{\bf Z}_n} ont une description
semblable pourvu que \m{\omega_{C\mid U}} soit trivial.

Les \'el\'ements de \m{\kg_n(U)} laissent invariant l'id\'eal \m{(t)}. On en
d\'eduit un morphisme de restriction
\[\rho_n:\kg_n\to\kg_{n-1}\]
qui est surjectif. L'application induite
\[H^1(\rho_n):H^1(C,\kg_n)\lra H^1(C,\kg_{n-1})\]
associe \`a une courbe de multiplicit\'e $n$ la courbe de multiplicit\'e
sous-jacente de multiplicit\'e \m{n-1}.

En associant \m{\nu_{\mid C}} \`a \m{\phi_{\mu\nu}} on d\'efinit un
morphisme surjectif \ \m{\xi_n:\kg_n\lra\ko_C^* .}
Le morphisme induit
\[H^1(\xi_n):H^1(C,\kg_n)\lra H^1(C,\ko_C^*)=\Pic(C)\]
associe \`a la courbe multiple $Y$ le fibr\'e en droites associ\'e $L$.
\end{subsub}

\sepsubsub

\begin{subsub}Courbes doubles -- \rm Soient \m{L\in\Pic(C)}, \m{g_0\in
H^1(C,\kg_2)} l'\'el\'ement correspondant \`a la courbe double triviale
\m{C_2^0} de fibr\'e en droites associ\'e $L$. On a des suites exactes
\xmat{0\ar[r] & T_C\ar[r] & \kg_2\ar[r]^-{\xi_2} & \ko_C^*\ar[r] & 0 ,}
\xmat{0\ar[r] & (T_C)^{g_0}=T_C\ot L\ar[r] & (\kg_2)^{g_0}=\AAut_C(C_2^0)
\ar[r]^-{\xi_2^{g_0}} & \ko_C^*\ar[r] & 0 ,}
Il en d\'ecoule une application surjective
\[\lambda_{g_0}:H^1(C,T_C\ot L)\lra H^1(\xi_2)^{-1}(L)\]
dont les fibres sont les orbites de l'action de \m{\C^*} (par multiplication)
sur \m{H^1(C,T_C\ot L)}. L'ensemble \m{H^1(\xi_2)^{-1}(L)} param\`etre les
courbes doubles de courbe r\'eduite $C$ et de fibr\'e en droites associ\'e $L$.
On peut retrouver ainsi la classification des courbes doubles de \cite{ba_ei}.
Si \m{g\in H^1(\xi_2)^{-1}(L)}, on note \m{E(g)} le fibr\'e vectoriel de rang 2
figurant dans une suite exacte
\[0\lra T_C\lra E(g)\lra L^*\lra 0\]
associ\'ee \`a un \'el\'ement de \m{\Ext^1_{\ko_C}(L^*,T_C)=H^1(C,T_C\ot L)} de
\m{\lambda_{g_0}^{-1}(g)}. Si \m{C_2} est la courbe double correspondant \`a
$g$, on a \ \m{E(g)=(\Omega_{C_2\mid C})^*}.
\end{subsub}

\sepsubsub

\begin{subsub}\label{intro1_0}Courbes de multiplicit\'e \m{n>2} -- \rm On
montre qu'on a 
\[\ker(\rho_n) \ \simeq \ T_C\oplus\ko_C\]
(proposition \ref{prox2}). Soient \m{g\in H^1(C,\kg_n)}, \m{C_n} la courbe
multiple correspondante, de fibr\'e en droites associ\'e $L$. Soient \m{g_k\in
H^1(C,\kg_k)} l'image de $g$, pour \m{2\leq k<n}. La description pr\'ecise de
\m{\kg_n} permet de montrer qu'on a
\[(T_C\oplus\ko_C)^g \ = \ E(g_2)\ot L^{n-1}\]
(proposition \ref{g_aut2_p1}). On a donc une suite exacte de faisceaux de
groupes
\[0\lra E(g_2)\ot L^{n-1}\lra\AAut_C(C_n)\lra\AAut_C(C_{n-1})\lra 0 .\]
Cela permet de d\'ecrire l'ensemble \m{H^1(\xi_n)^{-1}(g_{n-1})}, qui
param\`etre les courbes de multiplicit\'e $n$ qui sont des prolongements de
\m{C_{n-1}}. On a une application surjective
\[\lambda_g:H^1(C,E(g_2)\ot L^{n-1})\lra H^1(\xi_n)^{-1}(g_{n-1})\]
envoyant $0$ sur $g$, et dont les fibres sont les orbites d'une action de
\m{\Aut_C(C_{n-1})}.

L'espace \m{H^1(C,E(g_2)\ot L^{n-1})} param\`etre les classes d'\'equivalence
de prolongements de \m{C_{n-1}} en courbe de multiplicit\'e $n$ si on
consid\`ere que deux tels prolongements \m{C_n}, \m{C'_n} sont \'equivalents
s'il existe un isomorphisme \m{C_n\simeq C'_n} induisant l'identit\'e sur
\m{C_{n-1}}. L'action de \m{\Aut_C(C_{n-1})} sur \m{H^1(C,E(g_2)\ot L^{n-1})}
est la suivante : si \m{\alpha\in\Aut_C(C_{n-1})} et si \m{i:C_{n-1}\to C_n}
est un prolongement de \m{C_{n-1}}, alors \m{\alpha.C_n} est le prolongement
\xmat{C_{n-1}\ar[r]^-{\alpha^{-1}} & C_{n-1}\flinc[r]^-i & C_n .}

Le fait que \m{\ker(\rho_n)^g} est un faisceau coh\'erent sur $C$ entraine
aussi que \m{H^1(\rho_n)} est surjective. On en d\'eduit que toute courbe de
multiplicit\'e \m{n-1} peut \^etre prolong\'ee en courbe de multiplicit\'e $n$
(proposition \ref{cpa_def1}).

En utilisant le fait que toute courbe non projective est affine on en d\'eduit
aussi que si $C$ n'est pas projective les seules courbes multiples primitives
dont la courbe r\'eduite associ\'ee est $C$ sont les courbes triviales (ce
r\'esultat est prouv\'e pour les courbes doubles dans \cite{ba_ei}).

Soit \m{C_{n-1}} une courbe primitive de multiplicit\'e \m{n-1}, de courbe
r\'eduite associ\'ee $C$ projective. Les cas les plus simples o\`u on peut
d\'ecrire compl\`etement les prolongements de \m{C_{n-1}} en courbe de
multiplicit\'e $n$ sont \'enum\'er\'es en \ref{cl_simp}. Le cas o\`u
\m{C_{n-1}} est triviale est trait\'e dans \ref{act_cn1_2} : on obtient par
exemple que si $C$ est de genre positif, et si le degr\'e du fibr\'e en droites
associ\'e \`a \m{C_{n-1}} est n\'egatif, alors les prolongements de \m{C_{n-1}}
en courbes de multiplicit\'e $n$ non triviales sont param\'etr\'es pas un {\em
espace projectif tordu}. Les courbes triples sont compl\`etement classifi\'ees
dans le chapitre \ref{mult_3}.
\end{subsub}

\sepsubsub

\begin{subsub}\label{intro1_1}\'Eclatements -- \rm On d\'ecrit en \ref{ecl} la
proc\'edure d'\'eclatement d'un point $P$ d'une courbe multiple primitive en
terme de faisceaux de groupes. Elle se traduit par un morphisme surjectif
\m{b_{n,1}^P:H^1(C,\kg_n)\to H^1(C,\kg_n)}. On en d\'eduit une
g\'en\'eralisation du cas des courbes doubles trait\'e dans \cite{ba_ei}
(theorem 1.9).
\end{subsub}

\sepsubsub

\begin{subsub}\label{intro1_2} Courbes scind\'ees -- \rm Soient \m{n\geq 2} un
entier et \m{C_n} une courbe multiple primitive de multiplicit\'e $n$ et de
courbe r\'eduite associ\'ee $C$. On dit que \m{C_n} est {\em scind\'ee} s'il
existe une r\'etraction \m{C_n\to C}. Les exemples les plus simples sont les
courbes triviales. D'apr\`es \cite{ba_ei} les seules courbes doubles scind\'ees
sont les courbes triviales. Ceci n'est plus vrai en multiplicit\'e sup\'erieure
\`a 2. On peut classifier les courbes scind\'ees en \'etudiant un autre
faisceau de groupes non ab\'eliens, qui est un sous-faisceau de groupes (non
distingu\'e) de \m{\kg_n}, celui qui correspond aux automorphismes de la forme
\m{\phi_{0\nu}}.

On peut d\'ecrire enti\`erement les courbes scind\'ees de multiplicit\'e 3. On
consid\`ere une extension
\[0\lra\ko_C\lra E\lra L^*\lra 0\]
sur $C$, et soit \m{\sigma\in\Ext^1_{\ko_C}(L^*,\ko_C)} l'\'el\'ement
associ\'e. On consid\`ere le plongement \m{C\subset\P(E)} d\'efini par
l'inclusion \m{\ko_C\subset E}. Soit \m{C_3} la courbe de multiplicit\'e 3
correspondante dans la surface \m{\P(E)}. C'est une courbe scind\'ee de
fibr\'e en droites associ\'e $L$, et on montre qu'elles sont toutes de ce type.
La classe d'isomorphisme de \m{C_3} ne d\'epend que de \m{\C\sigma}. Les
courbes scind\'ees non triviales de multiplicit\'e $3$ et de fibr\'e en droites
associ\'e $L$ sont donc naturellement param\'etr\'ees par \m{\P(H^1(C,L))}
(proposition \ref{mult_3_4_1}).
\end{subsub}

\sepsubsub

\begin{subsub}Plongement des courbes multiples primitives dans des surfaces --
\rm (ce sujet n'est abord\'e que dans l'Introduction).
La surjectivit\'e du morphisme \m{b_{n,1}^P} de \ref{intro1_1} (autrement
dit le fait que le ``blowing-down'' est toujours possible) est \`a mettre en
relation avec le fait qu'une courbe multiple primitive ne peut pas
forc\'ement \^etre plong\'ee dans une surface lisse. D'apr\`es \cite{ba_ei},
theorem 7.1, la seule courbe double non triviale de courbe r\'eduite
associ\'ee \m{\P_1} pouvant \^etre plong\'ee dans une surface lisse est la
courbe double d\'eduite d'une conique plane. La d\'emonstration utilise
\cite{ha2}, theorem 4.1. Le m\^eme r\'esultat permet ais\'ement de prouver que
si \m{C_n} est une courbe multiple primitive de multiplicit\'e \m{n\geq 2} de
courbe r\'eduite associ\'ee $C$ et de fibr\'e en droites associ\'e $L$
plong\'ee dans une surface lisse, alors on a
\[\deg(L) \ \geq -4g-5 ,\]
$g$ d\'esignant le genre de $C$, ou alors \m{C_n} est la courbe de
multiplicit\'e $n$ induite par un plongement de $C$ dans un fibr\'e en espaces
projectifs comme indiqu\'e dans \ref{intro1_2}. Mais si $L$ est de degr\'e
tr\`es n\'egatif la construction des courbes multiples primitives indique qu'il
en existe beaucoup d'autres (les espaces \m{H^1(E(g_2)\ot L^k)} deviennent
tr\`es grands lorsque le degr\'e de $L$ diminue).
\end{subsub}
\end{sub}

\sepsub

\Ssect{Faisceaux des d\'erivations et fibr\'es tangents restreints}{intro_3}

D'apr\`es \cite{ba_ei} une courbe double \m{C_2} de courbe r\'eduite associ\'ee
$C$ est enti\`erement d\'etermin\'ee par le fibr\'e \m{\Omega_{C_2\mid C}} de
rang 2 sur $C$ et par la suite exacte (\ref{ecs}). On va g\'en\'eraliser
partiellement ce r\'esultat aux courbes
de multiplicit\'e sup\'erieure. Soit \m{C_n} une courbe multiple primitive de
multiplicit\'e \m{n>2}, de courbe r\'eduite associ\'ee $C$ projective et de
fibr\'e en droites associ\'e $L$. On pose
\[\kd_n \ = \ (\Omega_{C_n})^* ,\]
qu'on appelle le {\em faisceau des d\'erivations de \m{C_n}}. C'est un faisceau
{\em quasi localement libre} de rang g\'en\'eralis\'e 2 sur \m{C_n} (cf.
\cite{dr2}), de type \m{(0,\ldots,0,1,1)} (c'est-\`a-dire qu'il est localement
isomorphe \`a \m{\ko_{C_n}\oplus\ko_{C_{n-1}}}).

On pose
\[\T_{n-1} \ = \ \kd_{n\mid C_{n-1}} ,\]
qui est un faisceau localement libre de rang 2 sur \m{C_{n-1}}. On a aussi
\m{\T_{n-1}=(\Omega_{C_n\mid C_{n-1}})^*}. On l'appelle le {\em fibr\'e tangent
restreint} de \m{C_n}. Si \m{C_n} est plong\'ee dans une surface lisse, on a
en effet \m{\T_{n-1}=T_{S\mid C_{n-1}}}. Bien que ce soit un fibr\'e sur
\m{C_{n-1}} il d\'epend effectivement de \m{C_n}, et plus pr\'ecis\'ement de
l'inclusion \m{C_{n-1}\subset C_n}. La classe d'isomorphisme de \m{C_n} seule
d\'etermine \m{\T_{n-1}} \`a l'action de \m{\Aut_C(C_{n-1})} pr\`es.

On a des inclusions canoniques
\[\kd_{n-1}\ot\ki_C \ \subset \ \T_{n-1}\ot\ki_C \ \subset \ \kd_n\]
donc \m{\T_{n-1}} est un prolongement de \m{\kd_{n-1}} en fibr\'e vectoriel de
rang 2 sur \m{C_{n-1}}. On montre (corollaire
\ref{prol_dn2}) que si \m{\Omega_{C_2\mid C}} est stable, si \m{\deg(L)\leq
0}, et si en cas d'\'egalit\'e on a \m{L^k\not\simeq\ko_C} pour \m{1\leq k<2n},
alors 
le prolongement \m{\T_{n-1}} de \m{\kd_{n-1}} d\'etermine compl\`etement la
courbe \m{C_n}. De plus tout prolongement de \m{\kd_{n-1}} en fibr\'e de rang 2
sur \m{C_{n-1}} correspond \`a un prolongement de \m{C_{n-1}} en courbe de
multiplicit\'e $n$.

Les hypoth\`eses du corollaire \ref{prol_dn2} entrainent que
\m{\Aut_C(C_{n-1})} est trivial. La d\'emonstration repose sur une analyse dans
\ref{pl_ll} des prolongements de faisceaux quasi localement libres en faisceaux
localement libres. Cela conduit dans le cas de \m{\kd_{n-1}} sur \m{C_{n-1}}
\`a une application lin\'eaire
\[\theta:\Hom(T_C,E\ot L^{n-1})\lra H^1(C,E\ot L^{n-1})\]
(o\`u \m{E=\Omega_{C_2\mid C}^*}) telle que \m{\coker(\theta)} s'identifie
naturellement aux prolongements de \m{\kd_{n-1}} en fibr\'e vectoriel de rang 2
sur \m{C_{n-1}}. Si \m{u\in H^1(E\ot L^{n-1})}, on note \m{\T(u)} le
prolongement de \m{\kd_{n-1}} induit par $u$.

D'autre part \m{H^1(E\ot L^{n-1})} param\`etre aussi les prolongements de
\m{C_{n-1}} en courbe de multiplicit\'e $n$ (cf. \ref{intro1_0}). De $u$
on d\'eduit donc un prolongement \m{C_n} de multiplicit\'e $n$, d'o\`u un
fibr\'e de rang 2 \m{\T_{n-1}} prolongement de \m{\kd_{n-1}}. On note
\m{\T_{n-1}(u)} ce fibr\'e \m{\T_{n-1}}. On montre que
\m{\T_{n-1}(u)=\T(-(n-1)u)} (th\'eor\`eme \ref{prol_dn1}), d'o\`u on d\'eduit
le corollaire \ref{prol_dn2}.
\end{sub}

\sepsub

\Ssect{Automorphismes des courbes multiples primitives}{intro_2}

Soit \m{C_n} une courbe multiple primitive de multiplicit\'e \m{n\geq 2}, de
courbe r\'eduite associ\'ee $C$ projective et de fibr\'e en droites associ\'e
$L$. Il existe un morphisme canonique
\[\alpha:\Aut_C(C_n)\lra\C^*\]
associant \`a un automorphisme de \m{C_n} l'automorphisme induit de $L$, qui
est une homoth\'etie. On pose
\[\Aut^0_C(C_n) \ = \ \ker(\alpha) .\]
On montre (th\'eor\`eme \ref{fa2_3}) que si \m{C_n} n'est pas triviale, alors
\m{\imm(\alpha)} est fini. Mais \m{\Aut^0_C(C_n)} peut \^etre non trivial. Pour
le d\'ecrire il faut introduire une autre description des automorphismes
\m{\phi_{\mu\nu}} de \m{\C\DB{x,t}/(t^n)} tels que \m{\nu(x,0)=1}. Soit $D$ une
d\'erivation de \m{\C\DB{x,t}/(t^n)}. 
On suppose que
\begin{equation}D(\C\DB{x,t})\subset (t) \ \ \ \ {\rm et} \ \ \
\ D((t))\subset(t^2) .\end{equation}
Les d\'erivations qui ont ces propri\'et\'es sont celles qui se mettent sous la
forme
\[D \ = \ at\DIV{}{x}+bt^2\DIV{}{t} ,\]
avec \m{a,b\in\C\DB{x,t}}. On pose
\[\chi_D \ = \ \sigg_{k=0}^{n-1}\frac{1}{k!}D^k \ : \
\C\DB{x,t}/(t^n)\lra\C\DB{x,t}/(t^n) .\]
On d\'efinit ainsi un automorphisme de \m{\C\DB{x,t}/(t^n)} tel que
\m{\rho\circ\phi=\rho}. On montre (th\'eor\`eme \ref{autr_2}) que pour tout
automorphisme \m{\phi_{\mu\nu}} de \m{\C\DB{x,t}/(t^n)} tel que \m{\nu(x,0)=1}
il existe une unique telle d\'erivation $D$ telle que \m{\phi_{\mu\nu}=
\chi_D}.

Cette repr\'esentation des automorphismes a l'avantage de se globaliser. Pour 
\'etudier \m{\Aut^0_C(C_n)} il faut remplacer les d\'erivations pr\'ec\'edentes
par des sections de \m{\kd_{n-1}\ot\ki_C} (cf. \ref{intro_3}), et on d\'efinit
de mani\`ere analogue les \'el\'ements \m{\chi_D} de \m{\Aut^0_C(C_n)}
associ\'es aux sections $D$ de ce fibr\'e. On montre (th\'eor\`eme
\ref{fa2_3_8}) que les automorphismes \m{\chi_D} sont les seuls \'el\'ements de
 \m{\Aut^0_C(C_n)}. Comme ensemble, on a donc
\[\Aut^0_C(C_n) \ = \ H^0(C_{n-1}, \kd_{n-1}\ot\ki_C) ,\]
mais \'evidemment la structure de groupe de \m{\Aut^0_C(C_n)} n'est pas
l'addition (ce groupe n'est en g\'en\'eral pas commutatif).

On donne dans le corollaire \ref{fa2_3_7} de nombreux cas o\`u \m{\Aut_C(C_n)}
est trivial.

Soit \m{C_{n-1}} une courbe primitive de multiplicit\'e \m{n-1\geq 2} et de
courbe r\'eduite associ\'ee $C$ projective. On \'etudie en \ref{act_cn1}
l'action de \m{\Aut_C(C_{n-1})} sur \m{H^1(E\ot L^{n-1})} (o\`u
\m{E=\Omega_{C_2\mid C}^*}). On en d\'eduit la classification des prolongements
de \m{C_{n-1}} en courbes de multiplicit\'e $n$ lorsque \m{C_{n-1}} est
triviale. On obtient aussi un r\'esultat qui est utilis\'e dans la
classification des courbes triples.
\end{sub}

\sepsub

\Ssect{Plan des chapitres suivants}{intro_5}

Le chapitre \ref{prelim} rassemble des rappels ou des r\'esultats techniques
utilis\'es dans les autres chapitres. En particulier on d\'ecrit dans
\ref{pl_ll} les prolongements de faisceaux quasi localement libres en faisceaux
localement libres (cf. \cite{dr2}). Dans \ref{cech} on donne des r\'esultats
concernant la cohomologie de \v Cech qui permettront en particulier
d'identifier certains faisceaux par la suite.

\medskip

Le chapitre \ref{n_ab} est consacr\'e \`a des rappels de r\'esultats
concernant la cohomologie des faisceaux de groupes non n\'ecessairement
ab\'eliens. On s'inspire ici de \cite{fr}.

\medskip

Dans le chapitre \ref{fg_aut} on \'etudie en d\'etail les faisceaux de groupes
\m{\kg_n} mentionn\'es dans l'Introduction.

\medskip

Le chapitre \ref{HMPA} est consacr\'e \`a la construction et la
param\'etrisation des courbes multiples primitives en utilisant les faisceaux
de groupes \m{\kg_n}.

\medskip

Le chapitre \ref{fa_der} traite de la relation entre les prolongements d'une
courbe primitive \m{C_{n-1}} de multiplicit\'e \m{n-1} en courbes de
multiplicit\'e $n$ d'une part, et les prolongements du faisceau des
d\'erivations de \m{C_{n-1}} en fibr\'e de rang 2 d'autre part.

\medskip

Dans le chapitre \ref{fa2} on \'etudie les groupes d'automorphismes des
courbes multiples primitives et les actions de ces groupes intervenant dans la
param\'etrisation des courbes.

\medskip

Le chapitre \ref{mult_3} est consacr\'e \`a la description des courbes
primitives triples.
\end{sub}

\sepsub

\Ssect{Notations}{notations}

Soit \m{n\geq 1} un entier. On pose \ \m{A_n=\C[t]/(t^n)} ,
\m{{\mathbf Z}_n=\spec(\C[t]/(t^n))} .

Si \m{\alpha\in\C[t]}, on note \m{\alpha_i} le coefficient de \m{t^i} dans
\m{\alpha}. On emploie une notation analogue pour d'autres anneaux, par exemple
un anneau du type \m{A[t]/(t^n)}, $A$ \'etant un anneau commutatif unitaire.

Si $X$, $Y$ sont des vari\'et\'es alg\'ebriques, on note \m{p_X}, \m{p_Y} les
projections \m{X\times Y\to X}, \m{X\times Y\to Y}.

Si $Y$ est une sous-vari\'et\'e ferm\'ee de $X$, on notera \m{\ki_{Y,X}} ou
plus simplement \m{\ki_Y} le faisceau d'id\'eaux de $Y$ dans $X$.

\end{sub}

\sepsec

\section{Pr\'eliminaires}\label{prelim}

\Ssect{Courbes multiples primitives}{cmpr}

(cf. \cite{ba_fo}, \cite{va}, \cite{man}, \cite{be_fr}).

Soient $X$ une vari\'et\'e alg\'ebrique lisse connexe de dimension 3, et
\m{C\subset X} une courbe lisse connexe. On appelle {\em courbe multiple de
support $C$} un sous-sch\'ema de Cohen-Macaulay \m{ Y\subset X} tel que
l'ensemble des points ferm\'es de $Y$ soit $C$. Autrement dit, \m{Y_{red}=C}.

Soit $n$ le plus petit entier tel que \m{Y\subset C^{(n-1)}}, \m{C^{(k-1)}}
d\'esignant le $k$-i\`eme voisinage infinit\'esimal de $C$, c'est-\`a-dire \
\m{\ki_{C^{(k-1)}}=\ki_C^{k}} .
On a une filtration \ \m{C=C_1\subset C_2\subset\cdots\subset C_{n}=Y} \
o\`u $C_i$ est le plus grand sous-sch\'ema de Cohen-Macaulay contenu dans
\m{Y\cap C^{(i-1)}}. On appelle $n$ la {\em multiplicit\'e} de $Y$.

On dit que $Y$ est {\em primitive} si  pour tout point ferm\'e $x$ de $C$,
il existe une surface $S$ de $X$ contenant un voisinage de $x$ dans $Y$ et
lisse en $x$. Dans ce cas, \m{L=\ki_C/\ki_{C_2}} est un fibr\'e en droites sur
$C$ et on a \ \m{\ki_{C_{j}}/\ki_{C_{j+1}}=L^j} \ pour \m{1\leq j<n}.
Soit \m{P\in C}. Alors il existe des \'el\'ements $x$, $y$, $t$ de
\m{m_{X,P}} (l'id\'eal maximal de \m{\ko_{X,P}}) dont les images dans
 \m{m_{X,P}/m_{X,P}^2} forment une base, et que pour \m{1\leq i<n} on ait
\ \m{\ki_{C_i,P}=(x,y^{i})} .

Le cas le plus simple est celui o\`u $Y$ est contenue dans une surface lisse
$S$ de $X$. Dans ce cas il est m\^eme inutile de mentionner la vari\'et\'e
ambiente $X$.
Supposons $Y$ de multiplicit\'e $n$. Soient \m{P\in C} et \m{f\in\ko_{S,P}} une
\'equation locale de $C$. Alors on a \ \m{\ki_{C_i,P}=(f^{i})} \ pour \m{0\leq
j<n}, en particulier \m{I_{Y,P}=(f^n)}, et \ \m{L=\ko_C(-C)} .

On notera \m{\ko_n=\ko_{C_n}} et on consid\`erera \m{\ko_i} comme un faisceau
sur \m{C_n} de support \m{C_i} si \m{1\leq i<n}.

\end{sub}

\sepsub

\Ssect{Prolongements en faisceaux localement libres}{pl_ll}

On  utilise ici les notations et les r\'esultats de \cite{dr2}.

Soient $C$ une courbe projective lisse irr\'eductible, \m{Y=C_n} une courbe
multiple primitive de multiplicit\'e \m{n\geq 2} de courbe r\'eduite associ\'ee
$C$, et \m{L=\ki_C/\ki_{C_2}}, qui est un fibr\'e en droites sur $C$.

Soient $r$, $s$ des entiers positifs et $\ke$ un faisceau quasi localement
libre de type \m{(0,\ldots,0,r,s)} sur \m{C_n}, c'est-\`a-dire que pour tout
point ferm\'e $P$ de \m{C_n} il existe un voisinage $U$ de $P$ tel que
\m{\ke_{\mid U}\simeq r\ko_{n-1\mid U}\oplus s\ko_{n\mid U}}. Soient
\m{E=\ke^{(2)}} (c'est-\`a-dire le sous-faisceau de $\ke$ annulateur de
\m{\ki_C^{n-1}}) et \m{F=\ke/E}. Alors $E$ est un faisceau localement libre sur
\m{C_{n-1}} de rang \m{R=r+s} et $F$ est localement libre sur $C$ de rang $s$.
Le morphisme canonique \m{\ke\ot\ki_C\to\ke} induit un morphisme injectif de
fibr\'es vectoriels sur $C$
\[\tau_\ke : F\ot L\lra E_{\mid C} \ .\]
R\'eciproquement soient $E$ un fibr\'e vectoriel de rang \m{r+s} sur
\m{C_{n-1}} et $F$ un fibr\'e vectoriel de rang $s$ sur $C$. On s'int\'eresse
aux extensions
\[0\lra E \lra\ke\lra F\lra 0\]
sur \m{C_n}.

\sepprop

\begin{subsub}\label{pl_lc1}{\bf Proposition : } On a une suite exacte canonique
\begin{equation}\label{pl_equ1}
\xymatrix{0\ar[r] & \Ext^1_{\ko_C}(F,E_{\mid C}\ot L^{n-2})\ar[r] &
\Ext^1_{\ko_n}(F,E)\ar[r]^-\beta & \Hom(F\ot L,E_{\mid C})\ar[r] & 0 .}
\end{equation}
Soient \m{\sigma\in\Ext^1_{\ko_n}(F,E)} et \ \m{0\to E\to\ke\to F\to 0} \
l'extension correspondante. Alors $\ke$ est quasi localement libre de type 
\m{(0,\ldots,0,r,s)} si et seulement si \m{\beta(\sigma)} est injectif. On a
dans ce cas \m{\beta(\sigma)=\tau_\ke} et \m{E=\ke^{(2)}}.
\end{subsub}

\begin{proof}
D'apr\`es la suite spectrale des Ext (cf. \cite{go}, 7.3) on a une suite exacte
\[0\lra H^1(\HHom(F,E))\lra\Ext^1_{\ko_n}(F,E)\lra H^0(\EExt^1_{\ko_n}(F,E))
\lra 0 .\]
Soient $\L$ une extension de \m{\ki_C} en fibr\'e en droites sur \m{C_n} et
$\F$ une extension de $F$ en fibr\'e vectoriel sur \m{C_n} (ces fibr\'es
existent d'apr\`es \cite{dr2}, th\'eor\`eme 3.1.1). On a alors une r\'esolution
localement libre de $F$ sur \m{C_n}
\xmat{\cdots\F\ot\L^n\ar[r]^-{\alpha_2} & \F\ot\L\ar[r]^-{\alpha_1} &
\F\ar[r]^-{\alpha_0} & F\ar[r] & 0}
(par exemple \m{\alpha_0} est la restriction \m{\F\to\F_{\mid C}=F},
\m{\ker(\alpha_0)=\F\ot\ki_C} et \m{\alpha_1} est d\'efini par la restriction
\m{\L\to\ki_C}, etc...) On calcule \m{\HHom(F,E)} et \m{\EExt^1_{\ko_n}(F,E)}
au moyen de la r\'esolution pr\'ec\'edente, et on obtient
\m{\HHom(F,E)=\ker(A_1)}, \m{\EExt^1_{\ko_n}(F,E)=\coker(A_1)}, \m{A_1} \'etant
le morphisme \m{\HHom(\F,E)\to\HHom(\F\ot\L,E)} d\'eduit de \m{\alpha_1}. On en
d\'eduit imm\'ediatement que \m{\HHom(F,E)\simeq\HHom(F,E_{\mid C}\ot
L^{n-2})}, \m{\EExt^1_{\ko_n}(F,E)\simeq\HHom(F\ot L,E_{\mid C})}, ces
morphismes \'etant ind\'ependants du choix de $\L$ et $\F$, d'o\`u la suite
exacte (\ref{pl_equ1}).

Soient \m{\sigma\in\Ext^1_{\ko_n}(F,E)} et \ \m{0\to E\to\ke\to F\to 0} \
l'extension associ\'ee. On a \Nligne\m{\beta(\sigma)\in
H^0(\EExt^1_{\ko_n}(F,E))}, et si \m{P\in C},
\m{\beta(\sigma)(P)\in\EExt^1_{\ko_{nP}}(F_P,E_P)} correspond \`a l'extension
de \m{\ko_{nP}}-modules \ \m{0\to F_P\to\ke_P\to E_P\to 0} . Pour d\'emontrer
la seconde partie de la proposition il suffit de prouver le r\'esultat suivant
: si
\[0\lra (r+s)\ko_{n-1,P}\lra M\lra s\ko_{CP}\lra 0\]
est une suite exacte de \m{\ko_{nP}}-modules correspondant \`a
\[\alpha\in\Ext^1_{\ko_{nP}}(s\ko_{CP},(r+s)\ko_{n-1,P})\simeq\Hom(s\ko_{CP},
(r+s)\ko_{CP}) ,\]
alors $M$ est quasi libre de type \m{(0,\ldots,0,r,s)} si et seulement si
$\alpha$ est injectif. Soit \m{z\in\ko_{nP}} un g\'en\'erateur de \m{\ki_{CP}}.
La th\'eorie classique des extensions montre que $M$ est isomorphe au conoyau
du morphisme
\[\ov{\alpha}\oplus(\times z) : s\ko_{nP}\lra (r+s)\ko_{n-1,P}\oplus s\ko_{nP}
\ ,\]
o\`u \ \m{\ov{\alpha}:s\ko_{nP}\to (r+s)\ko_{n-1,P}} \ induit $\alpha$ et
\m{\times z} est la multiplication par $z$. Le r\'esultat en d\'ecoule
ais\'ement.
\end{proof}

\sepprop

\begin{subsub}\label{pl_lc1b} Prolongement en faisceau localement libre --
\rm On consid\`ere une extension \Nligne \m{0\to E\to\ke\to F\to 0}, avec $\ke$
quasi localement libre de type \m{(0,\ldots,0,r,s)}. On cherche \`a construire
des faisceaux localement libres $\kf$ sur \m{C_n} tels que \m{\ke\subset\kf} et
\m{\kf_2=\ke_2}. Pour un tel $\kf$ on a alors
\[\kf_C \ \simeq \ E_{\mid C}\ot L^* ,\]
et un diagramme commutatif avec lignes exactes
\xmat{ 0\ar[r] & E\ar[r]\fleq[d] &\ke\ar[r]\flinc[d] &
F\ar[r]\flinc[d]^{\tau_\ke\ot I_{L^*}} & 0 \\
0\ar[r] & E\ar[r] & \kf\ar[r] & E_{\mid C}\ot L^*\ar[r] & 0}
Soient \m{\sigma\in\Ext^1_{\ko_n}(F,E)}, \m{\sigma'\in\Ext^1_{\ko_n}(E_{\mid
C}\ot L^*,E)} associ\'es aux suites exactes pr\'ec\'edentes. Alors, si
\[\pi:\Ext^1_{\ko_n}(E_{\mid C}\ot L^*,E)\lra\Ext^1_{\ko_n}(F,E)\]
est le morphisme induit par \m{\tau_\ke\ot I_{L^*}}, on a 
\m{\pi(\sigma')=\sigma}.

R\'eciproquement, si \m{\sigma'\in\Ext^1_{\ko_n}(E_{\mid C}\ot L^*,E)} est tel
que \m{\pi(\sigma')=\sigma} et que dans l'extension correspondante \ \m{0\to
E\to\kf\to E_{\mid C}\ot L^*\to 0}, $\kf$ soit localement libre, alors on a une
inclusion naturelle \m{\ke\subset\kf} telle \m{\ke^{(2)}=\kf^{(2)}}. La
construction de $\kf$ se ram\`ene donc \`a celle de \m{\sigma'}.

On a un diagramme commutatif avec lignes exactes
\xmat{0\ar[r] & \Ext^1_{\ko_C}(E_{\mid C}\ot L^*,E_{\mid C}\ot
L^{n-2})\ar[r]^-{i'}\flon[d]^\psi & \Ext^1_{\ko_n}(E_{\mid C}\ot
L^*,E)\ar[r]^-{p'}\ar[d]^\pi &\End(E_{\mid C})\ar[r]\ar[d] & 0\\
0\ar[r] & \Ext^1_{\ko_C}(F,E_{\mid C}\ot L^{n-2})\ar[r]^-i &
\Ext^1_{\ko_n}(F,E)\ar[r]^-p & \Hom(F\ot L,E_{\mid C})\ar[r] & 0}
les fl\`eches verticales \'etant induites par \m{\tau_\ke}.

Si \m{\sigma'\in\Ext^1_{\ko_n}(E_{\mid C}\ot L^*,E)} et si \
\m{0\to E\to\kf\to E_{\mid C}\ot L^*\to 0} \ est l'extension correspondante,
alors $\kf$ est localement libre si et seulement si \m{p'(\sigma')} est un
isomorphisme (cf. \cite{dr2}), et dans ce cas l'extension pr\'ec\'edente est,
\`a un automorphisme de \m{E_{\mid C}} pr\`es, \'equivalente \`a l'extension
canonique \ \m{0\to\kf^{(2)}\to\kf\to\kf_{\mid C}\to 0}. D'autre part,
\m{p(\sigma)} est l'inclusion \m{F\ot L\subset E_{\mid C}}.
\end{subsub}

\sepprop

\begin{subsub}\label{pl_lc2}{\bf Lemme : } 
Il existe \m{\sigma'\in\Ext^1_{\ko_n}(E_{\mid C}\ot L^*,E)} tel que
\m{\pi(\sigma')=\sigma} et \m{p'(\sigma')=I_{E_{\mid C}}}.
\end{subsub}

\begin{proof}
Soit \m{\sigma_0\in\Ext^1_{\ko_n}(E_{\mid C}\ot L^*,E)} tel que \m{p'(\sigma_0)
=I_{E_{\mid C}}}. Alors on a \Nligne\m{p\circ\pi(\sigma_0)=p(\sigma)}. Donc il
existe \m{\lambda\in\Ext^1_{\ko_C}(F,E_{\mid C}\ot L^{n-2})} tel que \
\m{\sigma=\pi(\sigma_0)+i(\lambda)}. Puisque $\psi$ est surjectif il existe
\m{\mu\in\Ext^1_{\ko_C}(E_{\mid C}\ot L^*,E_{\mid C}\ot L^{n-2})} tel que
\m{\lambda=\psi(\mu)}. On a alors
\[\sigma=\pi(\sigma_0)+i\circ\psi(\mu)=\pi(\sigma_0+i'(\mu)) ,\]
et il suffit de prendre \m{\sigma'=\sigma_0+i'(\mu)}.
\end{proof}

\sepprop

\begin{subsub}\label{pl_cl2b} Classification des prolongements -- \rm
Il existe donc bien des prolongements localement libres $\kf$ de $\ke$ tels que
\m{\kf^{(2)}=\ke^{(2)}}. On note \m{{\bf Pr}(\ke)} l'ensemble des classes
d'isomorphisme de prolongements de $\ke$.

On s'int\'eresse maintenant aux choix possibles pour les $\kf$ pr\'ec\'edents.
Soit \m{\Gamma=(E_{\mid C}\ot L^*)/F}. On a une suite exacte
\xmat{
\Hom(F,E_{\mid C}\ot L^{n-2})\ar[r]^-\theta &
\Ext^1_{\ko_C}(\Gamma,E_{\mid C}\ot L^{n-2})\ar[r]^-\tau &
\Ext^1_{\ko_C}(E_{\mid C}\ot L^*,E_{\mid C}\ot L^{n-2})}
\xmat{\null \ar[r]^-\psi &
\Ext^1_{\ko_C}(F,E_{\mid C}\ot L^{n-2})\ar[r] & 0 \ .}
d\'eduite de la suite exacte \ \m{0\to\Gamma\to E_{\mid C}\ot L^*\to F\to 0}.
Il d\'ecoule de ce qui pr\'ec\`ede que les \m{\sigma''\in\Ext^1_{\ko_n}(E_{\mid
C}\ot L^*,E)} tels que \m{p'(\sigma'')=I_{E_{\mid C}}} et \m{\pi(\sigma'')=
\sigma} sont de la forme \m{\sigma''=\sigma'+i'(\mu)}, avec \
\m{\mu\in\coker(\theta)\subset\Ext^1_{\ko_C}(E_{\mid C}\ot L^*,E_{\mid C}\ot
L^{n-2})}. Les prolongements de possibles de $\ke$ sont donc les fibr\'es
d\'efinis pas ces \m{\sigma''}. On donc obtenu une application surjective
\[\coker(\theta)\lra{\bf Pr}(\ke) .\]

Le r\'esultat suivant permet sous certaines hypoth\`eses de conclure que les
prolongements de $\kf$ en fibr\'e vectoriel $\ke$ de rang 2 tel que
\m{\kf^{(2)}=\ke^{(2)}} sont classifi\'es par \m{\coker(\theta)}.
\end{subsub}

\sepprop

\begin{subsub}\label{pl_cl2b2} {\bf Proposition : } Soient $\kv$ un faisceau
localement libre de rang $R$ sur \m{C_n} et \m{V=\kv_{\mid C}}. On suppose que
\m{\deg(L)\leq 0}, \m{L^{Rk}\not\simeq\ko_C} pour \m{1\leq k<n} et que $V$
est stable. Alors $\kv$ est simple.
\end{subsub}

\begin{proof}
On fait une d\'emonstration par r\'ecurrence sur $n$, le r\'esultat \'etant
trivial si \m{n=1}. Supposons le vrai sur \m{C_{n-1}}. En consid\'erant les
suites exactes
\[0\lra\kv\ot\ki_C^{k+1}\lra\kv\ot\ki_C^k\lra V\ot L^k\lra 0\]
on voit ais\'ement que les hypoth\`eses entrainent que
\m{\Hom(V,\kv\ot\ki_C)=\nsp}. On voit aussi que \m{H^0(\ko_n)=\C}. Soit
\m{\phi\in\End(\kv)}. On en d\'eduit un diagramme commutatif avec lignes exactes
\xmat{0\ar[r] & \kv\ot\ki_C\ar[r]^-i\ar[d]^{\phi_0} & \kv\ar[r]^-p\ar[d]^\phi &
V\ar[r]\ar[d]^f & 0\\
0\ar[r] & \kv\ot\ki_C\ar[r]^-i & \kv\ar[r]^-p & V\ar[r] & 0}
D'apr\`es l'hypoth\`ese de r\'ecurrence et la stabilit\'e de $V$, il existe
\m{\lambda_0,\lambda\in\C} tels que\Nligne \m{\phi_0=\lambda_0I_{\kv\ot\ki_C}},
\m{f=\lambda I_V}. En utilisant l'inclusion \m{V\ot L^{n-1}\subset\kv\ot\ki_C}
on voit que \m{\lambda=\lambda_0}. Il en d\'ecoule qu'il existe un morphisme \
\m{\theta:V\to\kv\ot\ki_C} \ tel que \ \m{\phi-\lambda I_\kv=i\circ\theta\circ
p}. Puisque \m{\Hom(V,\kv\ot\ki_C)=\nsp} on a \m{\phi=\lambda I_\kv}.
\end{proof}

\sepprop

\begin{subsub}\label{pl_cl2b3} {\bf Corollaire : } On suppose que
\m{\deg(L)\leq 0}, \m{L^{Rk}\not\simeq\ko_C} pour \m{1\leq k<n} et que
\m{E_{\mid C}} est stable. Alors l'application canonique
\[\coker(\theta)\lra{\bf Pr}(\ke)\]
est bijective. Donc les prolongements de $\kf$ en fibr\'e
vectoriel $\ke$ de rang $R$ tel que \m{\kf^{(2)}=\ke^{(2)}} sont classifi\'es
par \m{\coker(\theta)}.
\end{subsub}

\begin{proof}
Soient \m{\sigma,\sigma'\in\Ext^1_{\ko_n}(E_{\mid C}\ot L^*,E)} tels que les
extensions correspondantes soient des prolongements localement libres $\ke$,
\m{\ke'} isomorphes de $\kf$. Un isomorphisme \m{\phi:\ke\simeq\ke'} induit un
diagramme commutatif
\xmat{0\ar[r] & E\ar[r]\ar[d]^{\phi'} & \ke\ar[r]\ar[d]^\phi &
E_{\mid C}\ot L^*\ar[r]\ar[d]^{\phi''} & 0\\
0\ar[r] & E\ar[r] & \ke'\ar[r] & E_{\mid C}\ot L^*\ar[r] & 0}
o\`u la suite exacte du haut correspond \`a $\sigma$ et celle du bas
\`a $\sigma'$. D'apr\`es la proposition \ref{pl_cl2b2}, $\phi$, \m{\phi'} et
\m{\phi''} sont des homoth\'eties de m\^eme rapport, d'o\`u \m{\sigma=\sigma'}.
\end{proof}

\sepsubsub

\begin{subsub}\label{pl_cl2c} Interpr\'etation en termes de 1-cocycles -- \rm
Soit \m{\ku=(U_i)} un recouvrement ouvert de \m{C_n} tel que $\kf$ soit
repr\'esent\'e par un 1-cocycle \m{(\phi_{ij})} de $\ku$, \m{\phi_{ij}:
R.\ko_{n\mid U_{ij}}\simeq R.\ko_{n\mid U_{ij}}}. On a donc des isomorphismes
\m{\phi_i:\kf_{\mid U_i}\simeq R.\ko_{n\mid U_{ij}}} tels que
\m{\phi_{ij}=\phi_i\phi_j^{-1}}. On suppose aussi que \m{L_{\mid U_i}} est
trivial~: soit \m{\xi_i\in H^0(U_i,L)} une section ne s'annulant en aucun point
et \m{\nu_{ij}=\xi_i^{-1}\xi_j\in\ko_n(U_{ij})^*}.

On suppose que $E$ et $F$ sont triviaux sur les \m{U_i} (en tant que fibr\'es
vectoriels sur \m{C_{n-1}} et $C$ respectivement). On peut dans ce cas mettre
\m{\phi_{ij}} sous la forme
\[\phi_{ij} \ = \ \begin{pmatrix} A_{ij} & B_{ij} \\ \xi_iC_{ij} & D_{ij}
\end{pmatrix} \ ,\]
avec \m{A_{ij}\in\ko_n(U_{ij})\ot\End(\C^r)},
\m{B_{ij}\in\ko_n(U_{ij})\ot L(\C^s,\C^r)},
\m{C_{ij}\in\ko_{n-1}(U_{ij})\ot L(\C^r,\C^s)},\Nligne
\m{D_{ij}\in\ko_n(U_{ij})\ot\End(\C^s)}, $\ke$ \'etant repr\'esent\'e par
\m{(\psi_{ij})},
\[\psi_{ij}:r.\ko_{n\mid U_{ij}}\oplus s.\ko_{n-1\mid U_{ij}} \ \simeq \
 r.\ko_{n\mid U_{ij}}\oplus s.\ko_{n-1\mid U_{ij}} ,\]
de la forme
\[\psi_{ij} \ = \ \begin{pmatrix} A_{ij} & B_{ij} \\ C_{ij} & \ov{D}_{ij}
\end{pmatrix} \ ,\]
o\`u \m{\ov{D}_{ij}\in\ko_{n-1}(U_{ij})\ot\End(\C^s)} est induit par
\m{D_{ij}}, c'est-\`a-dire que $\ke$ est obtenu en recollant au moyen des
\m{\phi_{ij}} les sous-faisceaux \ \m{r.\ko_{n\mid U_{ij}}\oplus
s.\ko_{n-1\mid U_{ij}}} \ de \m{R.\ko_{n\mid U_{ij}}}. De m\^eme les
\'el\'ements de \m{\Ext^1_{\ko_C}(E_{\mid C}\ot L^*,E_{\mid C}\ot L^{n-2})}
sont repr\'esent\'es par des cocycles de la forme
\[\mu_{ij} \ = \ \begin{pmatrix} a_{ij} & b_{ij} \\ c_{ij} & d_{ij}
\end{pmatrix} \ ,\]
et ceux qui proviennent de \m{\Ext^1_{\ko_C}(\Gamma,E_{\mid C}\ot L^{n-2})}
sont repr\'esent\'es par des cocycles tels que \m{a_{ij}=0}, \m{c_{ij}=0}.
Les autres fibr\'es vectoriels $\ke$ sont obtenus en prenant des 1-cocycles
\m{(\phi'_{ij})} de la forme
\[\phi'_{ij} \ = \ \begin{pmatrix} A_{ij} & B_{ij}+\xi_i^{n-1}b_{ij} \\
\xi_iC_{ij} & D_{ij}+\xi_i^{n-1}d_{ij}
\end{pmatrix} \ ,\]
avec \m{b_{ij}\in\ko_C(U_{ij})\ot L(\C^s,\C^r)},
\m{b_{ij}\in\ko_C(U_{ij})\ot\End(\C^s)}. Le fibr\'e vectoriel correspondant est
associ\'e \`a l'image dans \m{\coker(\theta)} de l'\'el\'ement de
\m{\Ext^1_{\ko_C}(\Gamma,E_{\mid C}\ot L^{n-2})} d\'efini par le cocycle
\m{(\xi_i^{n-1}b_{ij},\xi_i^{n-1}d_{ij})}.
\end{subsub}
\end{sub}

\sepsub

\Ssect{Cohomologie de \v Cech et extensions}{cech}

Soient $X$ une vari\'et\'e projective lisse irr\'eductible, et
\m{\ku=(U_i)_{i\in I}} un recouvrement de $X$ par des ouverts affines.
On emploie les notations habituelles :
\m{U_{i_0\ldots i_p}=U_{i_0}\cap\cdots\cap U_{i_p}}, etc.

On sait qu'on peut utiliser $\ku$ pour calculer la cohomologie des faisceaux coh\'erents sur $X$ (cf. \cite{ha}, chapter III, \para 4). On peut d\'efinir des \'el\'ements de leurs groupes de cohomologie en utilisant des {\em cocycles} au sens habituel. Mais parfois apparaissent des familles qu'on apellera aussi cocycles bien que cela n'en soit pas toujours au sens strict. Par exemple \'etant donn\'es un fibr\'e vectoriel $E$ et un fibr\'e en droites $L$ sur $X$, on peut repr\'esenter des \'el\'ements de \m{H^1(E\ot L)} de plusieurs fa\c cons : la premi\`ere consiste \`a consid\'erer des familles \m{(e_{ij}\ot\lambda_{ij})}, o\`u \m{e_{ij}\in H^0(U_{ij},E)}, \m{\lambda_{ij}\in H^0(U_{ij},L)}. Les relations de cocycle sont dans ce cas
\[e_{ij}\ot\lambda_{ij}+e_{jk}\ot\lambda_{jk} = e_{ik}\ot\lambda_{ik}\]
sur \m{U_{ijk}}. Mais on peut aussi partir d'un cocycle \m{(\gamma_{ij})}, \m{\gamma_{ij}\in\ko^*_X(U_{ij})} d\'efinissant $L$. Il existe donc des isomorphismes \m{\gamma_i:L_{\mid U_i}\simeq \ko_{U_i}} tels que \m{\gamma_{ij}=\gamma_i\gamma_j^{-1}}. Posons \m{e'_{ij}=\gamma_i(\lambda_{ij})e_{ij}\in H^0(U_{ij},E)}. Alors la famille \m{(e'_{ij})} repr\'esente le m\^eme \'el\'ement de \m{H^1(E\ot L)} que
\m{(e_{ij}\ot\lambda_{ij})}, mais avec des relations de cocycles diff\'erentes :
\[e'_{ij}+\gamma_{ij}e'_{jk} = e_{ik}\]
sur \m{U_{ijk}}.

\sepsubsub

\begin{subsub}\label{tors_lin} Torsion par un fibr\'e en droites -- \rm
Soient $E$ un fibr\'e vectoriel de rang $r$ et $L$ un fibr\'e en droites sur
$X$, triviaux sur tout ouvert \m{U_i}. Donc on peut d\'efinir $E$ (resp. $L$)
par un cocycle \m{(\phi_{ij})} (resp. \m{(\nu_{ij})}), avec
\m{\phi_{ij}:\ko_{U_i}\ot\C^r\simeq\ko_{U_i}\ot\C^r}, \m{\nu_{ij}:\ko_{U_i}
\simeq\ko_{U_i}}. Alors il est ais\'e de voir que \m{(\nu_{ij}\phi_{ij})} est
un cocycle d\'efinissant le fibr\'e vectoriel \m{E\ot L}.
\end{subsub}

\sepsubsub

\begin{subsub}{Construction des extensions - }\label{const_ext}\rm
Soient $E$, $F$ des fibr\'es vectoriels sur $X$, triviaux sur tous les ouverts
\m{U_i}, de rangs $r$ et $s$ respectivement. Alors $E$, $F$ sont d\'efinis par
des familles de cocycles \m{(\epsilon_{ij})}, \m{(\phi_{ij})}, o\`u
\m{\epsilon_{ij}:U_{ij}\to \GL(r,\C)}, \m{\phi_{ij}:U_{ij}\to \GL(s,\C)} , qui
sont d\'efinis eux-m\^emes \`a partir de trivialisations
\m{\epsilon_i:E_{\mid U_i}\to\ko_{U_i}\ot\C^r},
\m{\phi_i:F_{\mid U_i}\to\ko_{U_i}\ot\C^s} par
\m{\epsilon_{ij}=\epsilon_i\epsilon_j^{-1}}, \m{\phi_{ij}=\phi_i\phi_j^{-1}}~.

Soient \m{\sigma\in\Ext^1(F,E)}, et
\[0\lra E\lra\Gamma\lra F\lra 0\]
l'extension correspondante. Soit \m{(\delta_{ij})} (o\`u \m{\delta_{ij}:F_{\mid
U_{ij}}\to E_{\mid U_{ij}}}) une famille repr\'esentant $\sigma$. On en d\'eduit
une famille de cocycles \m{(\gamma_{ij})} repr\'esentant $\Gamma$, avec \
\m{\gamma_{ij} : U_{ij}\to\GL(r+s,\C)}, repr\'esent\'e par la matrice \
$\begin{pmatrix}\epsilon_{ij} & \epsilon_i\delta_{ij}\phi_j^{-1}\\
0 & \phi_{ij}\end{pmatrix}$ . Il existe alors une famille \m{(\gamma_i)}
d'isomorphismes, \m{\gamma_i:\Gamma_{\mid U_i}\to\ko_{U_i}\ot\C^{r+s}}, telle
que \m{\gamma_{ij}=\gamma_i\gamma_j^{-1}} pour tous $i$, $j$.
\end{subsub}

\sepsubsub

\begin{subsub}\label{desc_ext} Description des sections d'une extension -- \rm
Soient \m{E'}, \m{E''} des fibr\'es vectoriels sur $X$, \m{L'}, $D$ des fibr\'es
en droites sur $X$. On suppose que leurs restrictions \`a tous les ouverts
\m{U_i} sont triviales. On peut repr\'esenter \m{L'}, $D$ par des cocycles de
$\ku$, \m{(\lambda'_{ij})}, \m{(\delta_{ij})} respectivement. Soient
\[\lambda'_i:L'_{\mid U_i}\simeq\ko_{U_i} , \ \ \
\delta_i:D_{\mid U_i}\simeq\ko_{U_i}\]
des trivialisations telles que \m{\lambda'_{ij}=\lambda'_i{\lambda'_j}^{-1}},
\m{\delta_{ij}=\delta_i\delta_j^{-1}} sur \m{U_{ij}}.
Soit
\[\sigma\in\Ext^1(E''\ot D, E'\ot L'\ot D) ,\]
repr\'esent\'e par un cocycle \m{(\sigma_{ij})}, avec
\[\sigma_{ij}:(E''\ot D)_{\mid U_{ij}}\lra(E'\ot L'\ot D)_{\mid U_{ij}} .\] 
Soit 
\[0\lra E'\ot L'\ot D\lra\Gamma\lra E''\ot D\lra 0 \]
l'extension correspondant \`a $\sigma$. On peut construire $\Gamma$ en
recollant les fibr\'es\Nligne
\m{\lbrack(E''\ot D)\oplus(E'\ot L'\ot D)\rbrack_{\mid U_i}} au moyen des
automorphismes \m{\begin{pmatrix} I & \sigma_{ij}\\ 0 & I\end{pmatrix}}. Une
section de $\Gamma$ peut donc \^etre repr\'esent\'ee par des sections
\m{\epsilon''_i}, \m{\epsilon'_i} de \m{(E''\ot D)_{\mid U_i}},
\m{(E'\ot L'\ot D)_{\mid U_i}} respectivement telles que sur \m{U_{ij}} on ait
\[\epsilon'_i=\epsilon'_j+\sigma_{ij}\epsilon''_j , \quad
\epsilon''_i=\epsilon''_j .\]
Posons \m{\gamma_{ij}=\lambda'_i\sigma_{ij}:E''_{U_{ij}}\to E'_{U_{ij}}},
\m{e'_i=\lambda'_i\delta_i\epsilon'_i\in H^0(U_i,E'_i)},
\m{e''_i=\delta_i\epsilon''_i\in H^0(U_i,E''_i)}. Les relations pr\'ec\'edentes
s'\'ecrivent
\[e'_i=\delta_{ij}(\lambda'_{ij}e'_j+\gamma_{ij}e''_j) ,\quad
e''_i=\delta_{ij}e''_j .\]
R\'eciproquement, des familles \m{(e'_i)}, \m{(e''_i)} v\'erifiant ces relations d\'efinissent une section de $\Gamma$. Cela peut \^etre \'etendu aux sections locales de $\Gamma$. Ce r\'esultat est utilis\'e dans la proposition \ref{g_aut2_p1} pour identifier un fibr\'e vectoriel.
\end{subsub}
\end{sub}

\newpage

\section{Faisceaux de groupes non ab\'eliens}\label{n_ab}

On rappelle ici quelques r\'esultats de \cite{fr} qu'on adaptera \`a nos
besoins.

Soient $X$ un espace topologique et $G$ un faisceau de groupes sur $X$.

\sepsub

\Ssect{Cohomologie de dimension 1}{cohom1}

\begin{subsub}\label{cochai}Cochaines et cocycles -- \rm
Soit \m{\ku=(U_i)_{i\in I}} un recouvrement ouvert de $X$. On appelle {\em
0-cochaine de $\ku$ \`a valeurs dans $G$} une famille \m{(g_i)_{i\in I}}, avec
\m{g_i\in G(U_i)}. On appelle {\em 1-cochaine de $\ku$ \`a valeurs dans $G$}
une famille \m{(g_{ij})_{i,j\in I,i\not=j}}, o\`u \m{g_{ij}\in G(U_{ij})} (avec
\m{U_{ij}=U_i\cap U_j}). On dit que \m{(g_{ij})} est un {\em 1-cocycle} si pour
tous $i$, $j$, $k$ on a \m{g_{ij}g_{jk}=g_{jk}} sur \m{U_{ijk}=U_i\cap U_j\cap
U_k} (en convenant que pour tout $i$, \m{g_{ii}} est la section \'el\'ement
neutre). Deux 1-cochaines \m{(g_{ij})}, \m{(g'_{ij})} sont dites {\em
cohomologues} s'il existe une 0-cochaine \m{(h_i)} de $\ku$ \`a valeurs dans
$G$ telle que pour tous $i$, $j$ ont ait \m{g'_{ij}=h_ig_{ij}h_j^{-1}}. Si
c'est le cas, \m{(g_{ij})} est un 1-cocycle si et seulement si \m{(g'_{ij})} en
est un.
\end{subsub}

\sepsubsub

\begin{subsub}\label{cohom}Cohomologie de dimension 1 -- \rm
La cohomologie est une relation d'\'equivalence dans l'ensemble de 1-cocycles
de $\ku$ \`a valeurs dans $G$. L'ensemble des classes d'\'equivalence est not\'e
\m{H^1(\ku,G)}. En g\'en\'eral, si $G$ n'est pas un faisceau de groupes
ab\'eliens, il n'existe pas de structure naturelle de groupe sur cet ensemble,
mais il contient un \'el\'ement particulier, appel\'e {\em \'el\'ement neutre}
: la classe de la famille des \'el\'ements neutres des \m{G(U_{ij})}.

Si $\kv$ est un recouvrement ouvert de $X$ plus fin que $\ku$ il existe une
application naturelle \m{H^1(\ku,G)\to H^1(\kv,G)} qui est injective et envoie
l'\'el\'ement neutre de \m{H^1(\ku,G)} sur celui de \m{H^1(\kv,G)}. Ces
injections permettent de d\'efinir {\em l'ensemble de cohomologie \m{H^1(X,G)}
de $X$ \`a valeurs dans $F$} comme limite inductive des \m{H^1(\ku,G)}. Cet
ensemble est lui aussi muni d'un {\em \'el\'ement neutre} not\'e $e$. Pour tout
recouvrement ouvert $\ku$ de $X$ on a donc une application injective canonique
\m{H^1(\ku,G)\to H^1(X,G)} qui envoie l'\'el\'ement neutre de \m{H^1(\ku,G)} sur
celui de \m{H^1(X,G)}.

Soit \m{f:G\to H} un morphisme de faisceaux de groupes sur $X$. On en d\'eduit
un morphisme de groupes \m{H^0(f):H^0(X,G)\to H^0(X,H)} et une application
\m{H^1(f):H^1(X,G)\to H^1(X,H)} (d\'efinie de mani\`ere \'evidente au niveau des
1-cocycles) qui envoie l'\'el\'ement neutre de \m{H^1(X,G)} sur celui de
\m{H^1(X,H)}.
\end{subsub}

\sepsubsub

\begin{subsub}\label{acyc}Recouvrements acycliques -- \rm
Un recouvrement ouvert $\ku$ de $X$ est dit {\em acyclique pour $G$} si pour
tout ouvert $U$ de $\ku$ on a \m{H^1(U,G)=\lbrace e\rbrace}. Si $\ku$ est
acyclique pour $G$, alors l'application canonique \m{H^1(\ku,G)\to H^1(X,G)}
est bijective.
\end{subsub}

\sepsubsub

\begin{subsub}\label{act_fa}Actions de faisceaux de groupes -- \rm
Soient $G$, $H$ des faisceaux de groupes sur $X$, et supposons que $H$ agisse
sur $G$ (l'action est suppos\'ee compatible avec les structures de groupe).
Soit \m{z\in H^1(X,H)}, repr\'esent\'e par un 1-cocycle \m{(z_{ij})}
relativement \`a un recouvrement ouvert \m{\ku=(U_i)} de $X$. On d\'efinit un
nouveau faisceau de groupes \m{G^z} sur $X$ de la fa\c con suivante : on
recolle les faisceaux \m{G_{\mid U_i}} au moyen des isomorphismes
\m{z_{ij}:G_{\mid U_{ij}}\to G_{\mid U_{ij}}}. On a donc des isomorphismes
\m{\theta_i:G^z_{\mid U_i}\to G_{\mid U_i}} et des triangles commutatifs
\xmat{ & G^z_{\mid U_{ij}}\ar[ld]_{\theta_j}\ar[rd]^{\theta_i}\\
G_{\mid U_{ij}}\ar[rr]^-{z_{ij}} & & G_{\mid U_{ij}}}
Comme le sugg\`ere la notation employ\'ee, \m{G^z} ne d\'epend que du choix de
$z$.

Les sections de \m{G^z} sont repr\'esent\'ees par des familles \m{(\gamma_i)},
o\`u \m{\gamma_i\in G(U_i)}, telles que \m{\gamma_i=z_{ij}\gamma_j} sur
\m{U_{ij}}.

Les \'el\'ements de \m{H^1(X,G^z)} sont repr\'esent\'es par des familles
\m{(\rho_{ij})}, o\`u \m{\rho_{ij}\in G(U_{ij})}, telles que
\m{\rho_{ij}^{-1}\rho_{ik}=z_{ij}\rho_{jk}} sur \m{U_{ijk}} (il peut \^etre
n\'ecessaire de remplacer $\ku$ par des recouvrements plus fins pour
repr\'esenter ainsi tous les \'el\'ements de \m{H^1(X,G^z)}). Le 1-cocycle
\`a valeurs dans \m{G^z} correspondant \`a \m{(\rho_{ij})} est
\m{\theta_i^{-1}(\rho_{ij})}.

Si $G$ est un faisceau de groupes ab\'eliens, il en est de m\^eme de \m{G^z} et
les groupes de cohomologie \m{H^p(X,G^z)} sont d\'efinis pour tout entier
\m{p\geq 0}. On peut alors repr\'esenter dans ce cas les \'el\'ements de
\m{H^p(X,G^z)} de mani\`ere analogue. Par exemple les \'el\'ements de
\m{H^2(X,G^z)} sont repr\'esent\'es par des familles \m{(\rho_{ijk})}, o\`u
\m{\rho_{ijk}\in G(U_{ijk})}, telles que \m{\rho_{ijk}\rho_{ijl}^{-1}\rho_{ikl}=
z_{ij}\rho_{jkl}} sur \m{U_{ijkl}}.
\end{subsub}
\end{sub}

\sepsub

\Ssect{Suite exacte de cohomologie}{suite_c}

Dans toute la suite on suppose que $X$ est {\em paracompact}. C'est le cas par
exemple si c'est une vari\'et\'e alg\'ebrique munie de la topologie de Zariski.
Soient $G$ un faisceau de groupes sur $X$, \m{\Gamma\subset G} un sous-faisceaux
de groupes {\em distingu\'e} (c'est-\`a-dire que pour tout ouvert $U$ de $X$,
\m{\Gamma(U)} est un sous-groupe distingu\'e de \m{G(U)}. On peut donc d\'efinir
le {\em faisceau de groupes quotient} \m{G/\Gamma}, et on a une suite exacte
\xmat{0\ar[r] & \Gamma\ar[r]^i & G\ar[r]^p & G/\Gamma\ar[r] & 0 .}

Soient \m{A_1,\ldots,A_n} des ensembles munis d'un \'el\'ement privil\'egi\'e
$e$, et
\xmat{A_1\ar[r]^-{f_1} & A_2\ar[r] & \cdots\ar[r] & A_{n-1}
\ar[r]^-{f_{n-1}} & A_n}
une suite d'applications. On dit que cette suite est {\em exacte} si pour
\m{1<i<n} on a\Nligne \m{f_i^{-1}(e)=f_{i-1}(A_{i-1})} .

\sepprop

\begin{subsub}{\bf Proposition :}\label{coh_1} On a une suite exacte canonique
d'ensembles
\xmat{0\ar[r] & H^0(X,\Gamma)\ar[r]^-{H^0(i)} & H^0(X,G)\ar[r]^-{H^0(p)}
& H^0(X,G/\Gamma)\ar[r]^-\delta & \null\quad\quad\quad\quad}
\xmat{\quad\quad\quad\quad\quad\quad H^1(X,\Gamma)\ar[r]^-{H^1(i)} &
H^1(X,G)\ar[r]^-{H^1(p)} & H^1(X,G/\Gamma)\quad .}
La sous-suite
\xmat{0\ar[r] & H^0(X,\Gamma)\ar[r]^-{H^0(i)} & H^0(X,G)\ar[r]^-{H^0(p)}
& H^0(X,G/\Gamma)}
de la pr\'ec\'edente est une suite exacte de groupes. De plus, deux sections
$c$, $c'$ de \m{G/\Gamma} ont m\^eme image par $\delta$ si et seulement si
\m{cc'^{-1}} est l'image par \m{H^0(p)} d'une section de $G$.
\end{subsub}
 
\sepprop

Le groupe \m{H^0(X,G/\Gamma)} agit sur \m{H^1(X,\Gamma)} de la fa\c con
suivante : soient \m{a\in H^1(X,\Gamma)}, repr\'esent\'e par un 1-cocycle
\m{(\gamma_{ij})} d'un recouvrement \m{\ku=(U_i)} de $X$, et \m{c\in
H^0(X,G/\Gamma)}. On peut prendre $\ku$ suffisamment fin pour que pour tout
$i$, \m{c_{\mid U_i}} s'\'etende en \m{c_i\in G(U_i)}. On voit ais\'ement que
\m{(c_i\gamma_{ij}c_j^{-1})} est un 1-cocycle de $\ku$ \`a valeurs dans
$\Gamma$ et que l'\'el\'ement correspondant de \m{H^1(X,\Gamma)} ne d\'epend
que de $a$ et $c$. On notera \m{\sigma(c)a} cet \'el\'ement. On obtient ainsi
une action de \m{H^0(X,G/\Gamma)} sur \m{H^1(X,\Gamma)}, qui permet de
d\'ecrire l'application $\delta$ : on a, pour tout \m{c\in H^0(X,G/\Gamma)}
\[\delta(c)=\sigma(c^{-1})e .\]

\sepprop

\begin{subsub}{\bf Proposition :}\label{coh_2} Soient \m{\gamma,\gamma'\in
H^1(X,\Gamma)}. Alors on a \m{H^1(i)(\gamma)=H^1(i)(\gamma')} si et seulement si
$\gamma$ et $\gamma'$ sont dans la m\^eme \m{H^0(X,G/\Gamma)}-orbite.
\end{subsub}

\sepprop

\begin{subsub}\label{fib_h1p}
Les fibres de \m{H^1(p)} -- \rm On va utiliser l'action du faisceau de groupes
$G$ par conjugaison sur $\Gamma$, $G$ et \m{G/\Gamma}. Soient \m{\omega\in
H^1(X,G/\Gamma)} et \m{g\in H^1(X,G)} tel que \m{H^1(p)(g)=\omega}. Supposons
que $g$ soit repr\'esent\'e par un 1-cocycle \m{(g_{ij})} d'un recouvrement
ouvert \m{\ku=(U_i)} de $X$. En rempla\c cant au besoin $\ku$ par un
recouvrement ouvert plus fin, on voit qu'un autre \'el\'ement de
\m{H^1(p)^{-1}(\omega)} est repr\'esent\'e par un 1-cocycle de la forme
\m{(\gamma_{ij}g_{ij})}, avec \m{\gamma_{ij}\in\Gamma(U_{ij})}.
R\'eciproquement, une 1-cochaine \m{(\gamma_{ij}g_{ij})}, avec
\m{\gamma_{ij}\in\Gamma(U_{ij})} est un 1-cocycle si et seulement si pour tous
\m{i,j,k} on a
\[\gamma_{ij}^{-1}\gamma_{ik} = g_{ij}\gamma_{jk}g_{ij}^{-1} ,\]
c'est-\`a-dire si et seulement si \m{(\gamma_{ij})} induit un 1-cocycle de
$\ku$ \`a valeurs dans \m{\Gamma^g} (cf. \ref{act_fa}). On d\'efinit ainsi une
application surjective
\[\lambda_g:H^1(X,\Gamma^g)\lra H^1(p)^{-1}(\omega)\]
qui envoie l'\'el\'ement neutre de \m{H^1(X,\Gamma^g)} sur $g$. Si on part d'un
autre \'el\'ement \m{g'} de \m{H^1(p)^{-1}(\omega)}, il existe une bijection
\m{\alpha_{g,g'}:H^1(X,\Gamma^g)\to H^1(X,\Gamma^{g'})} telle que
le diagramme suivant soit commutatif :
\xmat{H^1(X,\Gamma^g)\ar[rrd]^{\lambda_g}\ar[dd]_{\alpha_{g,g'}} \\ & &
H^1(p)^{-1}(\omega)\\ H^1(X,\Gamma^{g'})\ar[rru]^{\lambda_{g'}}}
\end{subsub}

\sepsubsub

\begin{subsub}\label{cas_com} Le cas o\`u $\Gamma$ est commutatif -- \rm Dans
ce cas on a \m{\Gamma^g=\Gamma^{g'}}. En effet on peut supposer que \m{g'} est
repr\'esent\'e par un 1-cocycle du type \m{(u_{ij}g_{ij})} de $\ku$, avec
\m{u_{ij}\in\Gamma(U_{ij})} (on a \m{g'=\lambda_g(u)}, $u$ \'etant
repr\'esent\'e par \m{(u_{ij})}) et l'action de \m{g_{ij} :\Gamma(U_{ij})\to
\Gamma(U_{ij})} est la m\^eme que celle de \m{u_{ij}g_{ij}}. Soit \m{u\in
H^1(X,\Gamma^g)} l'\'el\'ement induit par \m{(u_{ij})}. Alors
\m{\alpha_{g,g'}:H^1(X,\Gamma^g)\to H^1(X,\Gamma^g)} est la translation
\m{w\to w-u}.
\end{subsub}

\sepsubsub

On suppose dans toute la suite que $\Gamma$ est {\em commutatif}. 

\sepsubsub

\begin{subsub}\label{imm_h1p}
L'image de \m{H^1(p)} -- \rm On consid\`ere l'action par conjugaison
 de \m{G/\Gamma} sur $\Gamma$. Soit \m{\omega\in H^1(X,G/\Gamma)}. On en
d\'eduit un \'el\'ement \m{\Delta(\omega)} de \m{H^2(X,\Gamma^\omega)},
d\'efini de la fa\c con suivante~: il existe un recouvrement ouvert
\m{\ku=(U_i)} de $X$ tel que $\omega$ soit d\'efini par un 1-cocycle
\m{(\omega_{ij})} de $\ku$ tel que pour tous $i$, $j$ il existe un \m{g_{ij}\in
G(U_{ij})} au dessus de \m{\omega_{ij}}. On suppose que \m{g_{ji}=g_{ij}^{-1}}.
Pour tous indices $i$, $j$, $k$ on a
\m{\gamma_{ijk}=g_{ij}g_{jk}g_{ki}\in\Gamma(U_{ijk})}, \m{(\gamma_{ijk})} est
un 2-cocycle de $\ku$ \`a valeurs dans \m{\Gamma^\omega} et \m{\Delta(\omega)}
est l'\'el\'ement induit de \m{H^2(X,\Gamma^\omega)}.
\end{subsub}

\sepprop

\begin{subsub}{\bf Proposition :}\label{coh_3} On a \m{\Delta(\omega)=0} si et
seulement si $\omega$ est dans l'image de \m{H^1(p)}.
\end{subsub}

\sepprop

\begin{subsub}{\bf Corollaire :}\label{coh_4} Si $X$ est une courbe
alg\'ebrique et si $\Gamma$ est un faisceau coh\'erent sur $X$, alors
\m{H^1(p)} est surjective.
\end{subsub}

\sepprop

\begin{subsub}\label{s_ex0}Suites exactes induites -- \rm Le faisceau de
groupes $G$ agit par conjugaison sur $\Gamma$, $G$ et \m{G/\Gamma}. Pour tout
\m{g\in H^1(X,G)}, \m{\Gamma^g} s'identifie naturellement \`a un sous-faisceau
de groupes de \m{G^g} et on a un isomorphisme canonique \m{(G/\Gamma)^g\simeq
G^g/\Gamma^g}. On a donc une suite exacte
\[0\lra\Gamma^g\lra G^g\lra (G/\Gamma)^g\lra 0 .\]
\end{subsub}

\sepprop

\begin{subsub}{\bf Proposition :}\label{coh_5} Soit \m{\omega=H^1(p)(g)}. Alors
les fibres de \m{\lambda_g:H^1(X,\Gamma^g)\lra H^1(p)^{-1}(\omega)} sont les
orbites de l'action de \m{H^0(X,(G/\Gamma)^g)} sur \m{H^1(X,\Gamma^g)}.
\end{subsub}

\sepprop

Soit \m{g'\in H^1(p)^{-1}(\omega)}. Alors on a \m{(G/\Gamma)^g=(G/\Gamma)^{g'}}
et \m{\Gamma^g=\Gamma^{g'}}, mais les actions de \Nligne
\m{H^0(X,(G/\Gamma)^g)} sur \m{H^1(X,\Gamma^g)} et \m{H^1(X,\Gamma^{g'})} ne
sont pas les m\^emes : si \m{(\gamma,v)\to \gamma v} d\'enote l'action de
\m{H^0(X,(G/\Gamma)^g)} sur \m{H^1(X,\Gamma^g)} et \m{(\gamma,v)\to \gamma *v}
celle sur \m{H^1(X,\Gamma^{g'})}, on a
\[\gamma *v \ = \ \gamma (v+u)-u\]

pour tous \m{g\in H^0(X,(G/\Gamma)^g)}, \m{v\in H^1(X,\Gamma^g)}, $u$
d\'esignant un \'el\'ement de \m{H^1(X,\Gamma^g)} tel que \m{\lambda_g(u)=g'}.

\end{sub}

\sepsec

\section{Faisceaux de groupes d'automorphismes}\label{fg_aut}

Soit $C$ une courbe irr\'eductible lisse. On \'etudie ici le faisceau
\m{\kg_n} des automorphismes de \m{C\times{\mathbf Z}_n} laissant invariants
$C$, et sa cohomologie de dimension 1. On commence par \'etudier les
automorphismes de l'anneau \m{\C\DB{x,t}} laissant invariante la projection
\m{\C\DB{x,t}\to\C\DB{x}}.

\sepsub

\Ssect{Automorphismes de \m{\C\DB{x,t}}}{aut_xt}

Soit \m{\rho:\C\DB{x,t}\to\C\DB{x}} le morphisme de $\C$-alg\`ebres d\'efini
par \ \m{\rho(\sigg_{i\geq 0,j\geq 0}\alpha_{ij}x^it^j)=\sigg_{i\geq
0}\alpha_{i0}x^i} .

Soit $\phi$ un endomorphisme de la $\C$-alg\`ebres \m{\C\DB{x,t}} tel que
\m{\rho\circ\phi=\rho}. Alors on peut \'ecrire
\begin{equation}
\phi(\alpha)=\alpha+\eta(\alpha)t \quad\text{ pour tout\ }\alpha\in\C\DB{x} ,
\ \ \ \ \phi(t)=\nu t ,\label{aut_xt1}\end{equation}
avec \m{\nu\in\C\DB{x,t}}, \m{\eta} \'etant une application
lin\'eaire \m{\C\DB{x}\to\C\DB{x,t}} telle que pour tous \m{\alpha,\beta\in
\C\DB{x}} on ait
\begin{equation}
\eta(\alpha\beta)=\alpha\eta(\beta)+\eta(\alpha)\beta+t\eta(\alpha)\eta(\beta)
.\label{aut_xt2}\end{equation}
R\'eciproquement, si $\eta$ est une telle application et \m{\nu\in\C\DB{x,t}},
il est ais\'e de voir qu'il existe un unique endomorphisme $\phi$ de
\m{\C\DB{x,t}} tel que  \m{\rho\circ\phi=\rho} et que (\ref{aut_xt1}) soit
v\'erifi\'e.

\sepprop

\begin{subsub}\label{lemx0}{\bf Lemme :} L'endomorphisme $\phi$ est un
automorphisme si et seulement si $\nu$ est inversible.
\end{subsub}

\begin{proof} Supposons que $\phi$ soit un automorphisme. Il existe alors
\m{u\in\C\DB{x,t}} tel que \m{\phi(u)=t}. On peut \'ecrire \m{u=u_0+tu_1}, avec
\m{u_0\in\C\DB{x}}. On a alors \m{t=u_0+t(\eta(u_0)+\phi(u_1)\nu)}, donc
\m{u_0=0} et \m{t=\phi(u_1)\nu t}, d'o\`u \m{\phi(u_1)\nu=1} et $\nu$ est
inversible.

R\'eciproquement, supposons que $\nu$ soit inversible. Soit \m{u\in\C\DB{x,t}}.
Il suffit de montrer qu'il existe une suite unique \m{(v_n)_{n\geq 0}} telle que
pour tout \m{n\geq 0} il existe \m{w_n\in\C\DB{x,t}} tel que
\[\phi(v_0+v_1t+\cdots+v_nt^n)=u+t^{n+1}w_n .\]
On montre l'existence et l'unicit\'e de \m{v_n} par r\'ecurrence sur $n$. On a
n\'ecessairement \m{v_0=\rho(u)}. Supposons \m{v_0,\ldots,v_n} construits. On
doit donc avoir
\[\phi(v_0+v_1t+\cdots+v_nt^n+v_{n+1}t^{n+1})=u+t^{n+1}(w_n+\nu^{n+1}
\phi(v_{n+1}))=u+t^{n+2}w_{n+1} .\]
La seule solution est \m{v_{n+1}=-\rho(\nu^{-n-1}w_n)}.
\end{proof}

\sepprop

Soient \m{\mu,\in\C\DB{x,t}}. On d\'efinit \m{\eta_\mu:\C\DB{x}\to\C\DB{x,t}}
par
\begin{equation}
\eta_\mu(\alpha) \ = \ \sigg_{i\geq 1}\frac{1}{i!}\frac{d^i\alpha}{dx^i}\mu^it^i .
\label{aut_xt3}\end{equation}
On v\'erifie ais\'ement que (\ref{aut_xt2}) est vraie pour \m{\eta=\eta_\mu}.
Il est facile de voir qu'on a
\begin{equation}
\alpha+\eta_\mu(\alpha)t \ = \ \alpha(x+\mu t) .
\label{aut_xt4}\end{equation}

Si \m{\nu\in\C\DB{x,t}}, on note \m{\phi_{\mu\nu}} l'endomorphisme de
\m{\C\DB{x,t}} tel que \Nligne\m{\phi_{\mu\nu}(\alpha)=\alpha+\eta_\mu(\alpha)t}
pour tout \m{\alpha\in\C\DB{x}} et \m{\phi_{\mu\nu}(t)=\nu t}.

\sepprop

\begin{subsub}\label{prox0}{\bf Proposition :} Si \m{\mu,\mu'\in\C\DB{x,t}} et
\m{\nu,\nu'\in\C\DB{x,t}}, on a
\[\phi_{\mu'\nu'}\circ\phi_{\mu\nu} \ = \ \phi_{\mu"\nu"} ,\]
avec
\[\mu"=\mu'+\nu'\phi_{\mu',\nu'}(\mu) , \ \ \ \ 
\nu"=\nu'\phi_{\mu',\nu'}(\nu) .\]
\end{subsub}

\begin{proof} D\'ecoule ais\'ement de (\ref{aut_xt3}).
\end{proof}

\sepprop

\'Etant donn\'e que \m{\phi_{0,1}=I_{\C\DB{x,t}}}, on en d\'eduit que si
\m{\nu\in\C\DB{x,t}} est inversible, alors on a \m{\phi_{\mu\nu}^{-1}=
\phi_{\mu',\nu'}}, avec
\begin{equation}\mu'=-\phi_{\mu\nu}^{-1}(\frac{\mu}{\nu}) , \ \ \ \
\nu'=\phi_{\mu\nu}^{-1}(\frac{1}{\nu}) .\label{aut_xt5}\end{equation}

\sepprop

\begin{subsub}\label{prox1}{\bf Th\'eor\`eme : } Pour tout endomorphisme
\m{\phi} de \m{\C\DB{x,t}} tel que \m{\rho\circ\phi=\rho} il existe un unique
\m{\mu\in\C\DB{x,t}} et un unique \m{\nu\in\C\DB{x,t}} tels que
\m{\phi=\phi_{\mu\nu}}.
\end{subsub}

\begin{proof}
La d\'emonstration la plus simple consiste \`a dire que $\phi$ est
enti\`erement d\'etermin\'e par ses valeurs sur $x$ et $t$, qui doivent \^etre
de la forme \m{\phi(x)=x+\mu t}, \m{\phi(t)=\nu t}. On donne ci dessous une
autre d\'emonstration qui peut s'\'etendre, contrairement \`a la pr\'ec\'edente,
\`a d'autres anneaux (cf. \ref{ext_v}).

On sait d\'ej\`a qu'il existe \m{\nu\in\C\DB{x,t}} tel que
\m{\phi(t)=\nu t}. On montre d'abord qu'on peut supposer que $\nu$ est
inversible. Dans le cas contraire on remplace $\phi$ par l'endomorphisme
\m{\phi'} \'egal \`a $\phi$ sur \m{\C\DB{x}} et tel que \m{\phi'(t)=(1+\nu)t}.
S'il existe \m{\mu\in\C\DB{x,t}} tel que \m{\phi'=\phi_{\mu,\nu+1}} alors on a
\m{\phi=\phi_{\mu\nu}}. On cherche \m{\mu=\sigg_{i\geq 0}\mu_it^i} (avec
\m{\mu_i\in\C\DB{x}}), et on va construire \m{\mu_n} par r\'ecurrence sur $n$.
Posons
\[\mu^{(n)}=\sigg_{i=0}^n\mu_it^i , \ \ \ \ \nu^{(n)}=\sigg_{i=0}^n\nu_it^i .\]
On va choisir les \m{\mu_i} de telle sorte que pour tout \m{n\geq 0} on ait
\begin{equation}
\imm(\phi\circ\phi_{\mu^{(n)}\nu^{(n)}}^{-1}-I_{\C\DB{x,t}}) \ \subset \
(t^{n+2}) .\label{aut_xt6}\end{equation}
Remarquons que cette inclusion reste vraie si on ajoute \`a \m{\mu^{(n)}} o\`u
\`a \m{\nu^{(n)}} un multiple de \m{t^{n+1}}. En particulier on aura finalement
\m{\phi\circ\phi_{\mu\nu}^{-1}=I_{\C\DB{x,t}}}, c'est-\`a-dire
\m{\phi=\phi_{\mu\nu}}.

Consid\'erons la d\'ecomposition (\ref{aut_xt1}) de $\phi$. On peut \'ecrire
\m{\eta=\eta_0+t\eta_1}, avec \ \m{\eta_0:\C\DB{x}\to\C\DB{x}}. La relation
(\ref{aut_xt2}) entra\^\i ne que \m{\eta_0} est une d\'erivation, donc est de la
forme \m{\eta_0=\mu_0\frac{d}{dx}}, avec \m{\mu_0\in\C\DB{x}}. On v\'erifie
alors ais\'ement l'inclusion (\ref{aut_xt6}) pour \m{n=0}.

Supposons \m{\mu_0,\ldots,\mu_n} trouv\'es, tels que (\ref{aut_xt6}) soit vraie.
On peut donc \'ecrire
\[\phi\circ\phi_{\mu^{(n)}\nu^{(n)}}^{-1}=I_{\C\DB{x,t}}+t^{n+2}\theta ,\]
avec \m{\theta:\C\DB{x,t}\to\C\DB{x,t}}, qu'on peut mettre sous la forme
\m{\theta=\theta_0+\theta_1t}, avec\Nligne
\m{\theta_0:\C\DB{x,t}\to\C\DB{x}}. Le fait
que \m{\phi\circ\phi_{\mu^{(n)}\nu^{(n)}}^{-1}} est un morphisme d'anneaux
implique que \m{\theta_0} est une d\'erivation. On peut donc \'ecrire
\m{\theta_0=\mu_{n+1}\frac{d}{dx}}, avec \m{\mu_{n+1}\in\C\DB{x}}. Il reste \`a
v\'erifier (\ref{aut_xt6}) pour \m{n+1}. On a
\[\phi\circ\phi_{\mu^{(n)}\nu^{(n)}}^{-1}(t)=\phi(\phi_{\mu^{(n)}\nu^{(n)}}^{-1}
(\frac{1}{\nu^{(n)}}).t)=\phi\circ\phi_{\mu^{(n)}\nu^{(n)}}^{-1}(\frac{1}
{\nu^{(n)}}).\nu t=(\frac{1}{\nu^{(n)}}+t^{n+2}\theta (\frac{1}{\nu^{(n)}}))
.\nu t \]
\[=(\frac{1}{\nu^{(n)}}+t^{n+2}\theta (\frac{1}{\nu^{(n)}}))(\nu^{(n)}t+
\nu_{n+1}t^{n+2}+Kt^{n+3})=t+\frac{\nu_{n+1}}{\nu_0}t^{n+2}+\epsilon t^{n+3}\]
pour des $K$ et $\epsilon$ convenables. Posons
\[\phi'=\phi_{\mu_{n+1}t^{n+1},1+\frac{\nu_{n+1}}{\nu_0}t^{n+1}} .\]
On d\'eduit de ce qui pr\'ec\`ede que
\[\imm(\phi\circ\phi_{\mu^{(n)}\nu^{(n)}}^{-1}-\phi')\subset(t^{n+3}) ,\]
Donc
\[\imm(\phi\circ\phi_{\mu^{(n)}\nu^{(n)}}^{-1}\circ\phi'^{-1}-
I_{\C\DB{x,t}})\subset(t^{n+3}) .\]
On a \m{\phi'\circ\phi_{\mu^{(n)}\nu^{(n)}}=\phi_{\mu'\nu'}}, avec
\[\mu'=\mu_{n+1}t^{n+1}+(1+\frac{\nu_{n+1}}{\nu_0}t^{n+1})\phi'(\mu^{(n)}) , \ \
\ \ \nu'=(1+\frac{\nu_{n+1}}{\nu_0}t^{n+1})\phi'(\nu^{(n)})\]
d'apr\`es la proposition \ref{prox0}. Il est facile de voir que \ \m{\mu'-
\mu^{(n+1)},\nu'-\nu^{(n+1)}\in(t^{n+3})}, donc on a bien
\[\imm(\phi\circ\phi_{\mu^{(n+1)}\nu^{(n+1)}}^{-1}-I_{\C\DB{x,t}}) \ \subset \
(t^{n+3}) .\]
\end{proof}

\sepprop

\begin{subsub}Notation : \rm Soit \m{\kg(\C\DB{x,t})} l'ensemble des
automorphismes $\phi$ de la $\C$-alg\`ebre \m{\C\DB{x,t}} tels que
\m{\rho\circ\phi=\rho}. C'est un sous-groupe du groupe de tous les
automorphismes de \m{\C\DB{x,t}}.
\end{subsub}

\sepprop

\begin{subsub}\label{prox1b}{\bf Proposition : } Soient \m{\mu,\nu\in
\C\DB{x,t}} et \m{\theta\in\C\DB{x,t}}. Alors on a
\[\frac{\partial}{\partial x}\phi_{\mu\nu}(\theta) \ = \ (1+t\DIV{\mu}{x})
\phi_{\mu\nu}(\DIV{\theta}{x}) \ + \
t\DIV{\nu}{x}\phi_{\mu\nu}(\DIV{\theta}{t}) , \]
\[\frac{\partial}{\partial t}\phi_{\mu\nu}(\theta) \ = \ (\mu+t\DIV{\mu}{t})
\phi_{\mu\nu}(\DIV{\theta}{x}) \ + \
(\nu+t\DIV{\nu}{t})\phi_{\mu\nu}(\DIV{\theta}{t})  .\]
\end{subsub}

\begin{proof}
V\'erification imm\'ediate.
\end{proof}

\sepprop
\end{sub}

\sepsub

\Ssect{Autres descriptions des isomorphismes}{chgt}

Soient \m{D:\C\DB{x}\to \C\DB{x}} une d\'erivation et \m{\mu,\nu\in\C\DB{x,t}}.
On note \m{\psi_{\mu\nu}^D} l'endomorphisme de la $\C$-alg\`ebre \m{\C\DB{x,t}}
d\'efini par
\[\psi_{\mu\nu}^D(\alpha)=\sigg_{n\geq 0}\frac{1}{n!}D^n(\alpha)(\mu t)^n\]
pour tout \m{\alpha\in\C\DB{x}}, et \m{\psi_{\mu\nu}^D(t)=\nu t}. En
particulier on a \m{\psi_{\mu\nu}^{d/dx}=\phi_{\mu\nu}}. D'apr\`es la
proposition \ref{prox1} il existe un unique \m{\gamma\in\C\DB{x,t}} tel que
\m{\psi_{\mu\nu}^D=\phi_{\gamma\nu}}. Il existe \m{a\in\C\DB{x}} tel que
\m{D=a\frac{d}{dx}}. On va donner une expression de $\gamma$ en fonction de $a$
et $\mu$.

\sepprop

\begin{subsub}{\bf Proposition : }\label{lemx1}
Soit \m{\gamma=\mu\sigg_{k\geq 0}\frac{1}{(k+1)!}D^k(a)(\mu t)^k} .
Alors on a \m{\psi_{\mu\nu}^D=\phi_{\gamma,\nu}} .
\end{subsub}

\begin{proof}
Cela d\'ecoule imm\'ediatement de la formule \
\m{\phi_{\gamma,\nu}(x)=x+\gamma t} .
\end{proof}

\sepprop

\begin{subsub}\label{autr_0}Utilisation de d\'erivations par rapport \`a $x$ et
$t$ -- \rm Soit \m{D:\C\DB{x,t}\to \C\DB{x,t}} une d\'erivation. On suppose que
\begin{equation}\label{autr_0a}D(\C\DB{x,t})\subset (t) \ \ \ \ {\rm et} \ \ \
\ D((t))\subset(t^2) .\end{equation}
Les d\'erivations qui ont ces propri\'et\'es sont celles qui se mettent sous la
forme
\[D \ = \ at\DIV{}{x}+bt^2\DIV{}{t} ,\]
avec \m{a,b\in\C\DB{x,t}}. On pose
\[\chi_D \ = \ \sigg_{k\geq 0}\frac{1}{k!}D^k \ : \ \C\DB{x,t}\lra\C\DB{x,t} .\]
Cette formule a un sens car on v\'erifie ais\'ement que pour tout \m{k\geq 0}
on a \m{\imm(D^k)\subset(t^{k+1})}.
On note \m{{\rm Der}_0(\C\DB{x,t})} l'espace vectoriel des d\'erivations
poss\`edant les propri\'et\'es (\ref{autr_0a}).
\end{subsub}

\sepprop 

\begin{subsub}\label{autr_1}{\bf Proposition : } 1 - Soient \m{D,D'\in
{\rm Der}_0(\C\DB{x,t})}. Si $D$ et \m{D'} commutent, alors on a \
\m{\chi_{D+D'}=\chi_D\circ\chi_{D'}}.

2 - Pour tout \m{D\in {\rm Der}_0(\C\DB{x,t})}, on a \ \m{\chi_D\in
\kg(\C\DB{x,t})}.
\end{subsub}

\begin{proof} L'assertion 1- est imm\'ediate, et 2- s'en d\'eduit car
\m{\chi_{-D}} est l'inverse de \m{\chi_D} d'apr\`es 1-.
\end{proof}

\sepprop

\begin{subsub}\label{autr_1b}{\bf Proposition : } Soit \m{D\in {\rm
Der}_0(\C\DB{x,t})}. Alors on a \ \m{\chi_D=\phi_{\mu\nu}}, avec
\begin{equation}\label{autr_1c}
\mu \ = \ \frac{1}{t}\sigg_{k\geq 1}\frac{1}{k!}D^k(x) , \ \ \ \
\nu \ = \ 1+\frac{1}{t}\sigg_{k\geq 1}\frac{1}{k!}D^k(t) .\end{equation}
\end{subsub}

Ces formules ont un sens puisque \m{\imm(D)\subset(t)}.

\begin{proof} 
D\'ecoule imm\'ediatement du fait que \m{\chi_D(x)=x+\mu t} et \m{\chi_D(t)=
\nu t}.
\end{proof}

\sepprop

Soit \m{\C\DB{x}^*} le groupe des \'el\'ements inversibles de \m{\C\DB{x}}.
On d\'efinit un morphisme surjectif de groupes
\[\xi_{\C\DB{x,t}}:\kg(\C\DB{x,t})\lra\C\DB{x}^*\]
en associant \`a \m{\phi_{\mu\nu}} le scalaire \m{\nu(x,0)}. On note
\m{\kg_0(\C\DB{x,t})} le sous-groupe de \m{\kg(\C\DB{x,t})} noyau de ce
morphisme. Il contient tous les automorphismes \m{\chi_D}.

\sepprop

\begin{subsub}\label{autr_2}{\bf Th\'eor\`eme : } Pour tout \m{\phi\in
\kg_0(\C\DB{x,t})} il existe un \m{D\in{\rm Der}_0(\C\DB{x,t})}, unique, tel
que \m{\phi=\chi_D}.
\end{subsub}

\begin{proof} 
On cherche $D$ sous la forme \ \m{D=at\DIV{}{x}+bt^2\DIV{}{t}}. Posons
\[a=\sigg_{i\geq 0}a_it^i , \ \ \ b=\sigg_{i\geq 0}b_it^i ,\]
avec \m{a_i,b_i\in\C\DB{x}}, et
\[a^{(n)}=\sigg_{i=0}^na_it^i , \ \ \ b^{(n)}=\sigg_{i=0}^nb_it^i , \ \ \ \
D_n=a^{(n)}t\DIV{}{x}+b^{(n-1)}t^2\DIV{}{t} .\]
Posons \m{\phi=\phi_{\mu\nu}}, avec \m{\nu(x,0)=1}. On pose
\[\mu_n \ = \ \frac{1}{t}\sigg_{k\geq 1}\frac{1}{k!}D_n^k(x) , \ \ \ \
\nu_n \ = \ 1+\frac{1}{t}\sigg_{k\geq 1}\frac{1}{k!}D_n^k(t) .\]
Alors il existe un choix unique de \m{a_0} de telle sorte que
\[\mu=\mu_0 \ ({\rm mod} \ t) \ \ \ \ {\rm et} \ \ \ \ \nu=\nu_0 \ ({\rm mod}
\ t) .\]
C'est \m{a_0=\mu(x,0)}. Il reste \`a prouver que si
\[\mu=\mu_{n-1} \ ({\rm mod} \ t^n) \ \ \ \ {\rm et} \ \ \ \ \nu=\nu_{n-1}
\ ({\rm mod} \ t^n) ,\]
alors il existe \m{a_n,b_{n-1}\in\C\DB{x}}, uniques, tels que
\begin{equation}\label{autr_3}
\mu=\mu_n \ ({\rm mod} \ t^{n+1}) \ \ \ \ {\rm et} \ \ \ \ \nu=\nu_n
\ ({\rm mod} \ t^{n+1}) .\end{equation}
On a
\[D_n(x)=D_{n-1}(x)+a_nt^{n+1}, \ \ \ \ D_n^2(x)=D_{n-1}^2(x)+a_0b_{n-1}t^{n+1}
 \ ({\rm mod} \ t^{n+2}) ,\]
\m{D_n^k(x)=D_{n-1}^k(x) \ ({\rm mod} \ t^{n+2})} \ si \m{k\geq 3},
\[D_n(t)=D_{n-1}(t)+b_{n-1}t^{n+1} \ ({\rm mod} \ t^{n+2}), \ \ \ \
D_n^k(t)=D_{n-1}^k(t) \ \ {\rm si} \ \ \m{k\geq 2} .\]
Donc
\[\mu_n=\mu_{n-1}+(a_n+\frac{a_0b_{n-1}}{2})t^n \ ({\rm mod} \ t^{n+1}) \ \ \ \
{\rm et} \ \ \ \ \nu_n=\nu_{n-1}+b_{n-1}t^n \ ({\rm mod} \ t^{n+1}).\]
Il est donc clair qu'il existe un choix unique de \m{a_n,b_{n-1}} tel qu'on
obtienne (\ref{autr_3}).
\end{proof}

\sepprop

On peut donc identifier comme ensembles \m{\kg_0(\C\DB{x,t})} et \m{{\rm
Der}_0(\C\DB{x,t})}. Mais la structure de groupes sur ce dernier induite par
celle de \m{\kg_0(\C\DB{x,t})} n'est \'evidemment pas l'addition (cf.
\ref{K_neq3}) et n'est m\^eme pas commutative. On notera '*' la loi de groupe
sur \m{{\rm Der}_0(\C\DB{x,t})}, c'est-\`a-dire que \
\m{\chi_{D*D'}=\chi_D\circ\chi_{D'}}.
\end{sub}

\sepsub

\Ssect{Extension des r\'esultats \`a d'autres anneaux}{ext_v}

Soient $C$ une courbe projective irr\'eductible lisse, \m{P\in C} un point
ferm\'e et \m{U\subset C} un ouvert non vide sur lequel le fibr\'e
canonique \m{\omega_C} est trivial. On obtient des r\'esultats analogues aux
pr\'ec\'edents, et on utilisera des notations analogues, si on remplace
\m{\C\DB{x}} par une des $\C$-alg\`ebres suivantes : \m{\C[x]}, \m{\C(x)},
\m{\C((x))}, \m{\ko_{CP}}, \m{\ko_C(U)}. Traitons par exemple le dernier cas.

Soit \m{\frac{d}{dx}} une section de \m{T_{C\mid U}} engendrant
\m{T_C} sur $U$. Soit \m{\rho:\ko_C(U)[t]\to\ko_C(U)} le morphisme
canonique. Soient \m{\mu,\nu\in\ko_C(U)[t]}. On d\'efinit un endomorphisme
\m{\phi_{\mu\nu}} de \m{\ko_C(U)[t]} par
\[\phi_{\mu\nu}(\alpha) \ = \ \sigg_{i\geq 0}\frac{1}{i!}
\frac{d^i\alpha}{dx^i}\mu^it^i \]
pour tout \m{\alpha\in\ko_C(U)}, et \m{\phi_{\mu\nu}(t)=\nu t}.
Il est facile de voir qu'on a
\[\phi_{\mu,\nu}(\alpha) \ = \ \alpha(x+\mu t)\]
et que \m{\phi_{\mu\nu}} est inversible si et seulement si $\nu$ l'est.
De plus, tout endomorphisme $\phi$ de \m{\ko_C(U)[t]} tel que
\m{\rho\circ\phi=\rho} est de la forme \m{\phi=\phi_{\mu,\nu}}, pour des $\mu$,
$\nu$ uniques.

\sepsubsub

\begin{subsub}\label{K_n} Quotients par \m{(t^n)} -- \rm
Soit \m{n\geq 2} un entier. Une autre extension possible des r\'esultats
pr\'ec\'edents consiste \`a quotienter par l'id\'eal \m{(t^n)}, c'est-\`a-dire
que si $A$ est une des $\C$-alg\`ebres \m{\C\DB{x}}, \m{\ko_{CP}},
\m{\ko_C(U)}, on s'int\'eresse aux endomorphismes de \m{A\DB{t}/(t^n)} qui
laissent invariante la projection \m{A\DB{t}/(t^n)\to A}. La seule diff\'erence
est qu'il faut prendre \m{\mu,\nu} dans \m{A\DB{t}/(t^{n-1})}. Le morphisme
\m{\xi_{\C\DB{x}}} de \ref{chgt} se g\'en\'eralise en un morphisme surjectif
\[\xi_A^n:\kg(A\DB{t}/(t^n))\to A^*\]
de noyau \m{\kg_0(A\DB{t}/(t^n))} et le th\'eor\`eme \ref{autr_2} s'\'etend.
\end{subsub}

\sepsubsub

\begin{subsub}\label{K_neq2} Le cas \m{n=2} -- \rm Les automorphismes
\m{\phi_{\mu\nu}} d\'ependent de param\`etres \m{\mu,\nu\in A}. On a
\[\phi_{\mu\nu}^{-1} \ = \ \phi_{-\frac{\mu}{\nu},\frac{1}{\nu}} .\]
\end{subsub}

\sepsubsub

\begin{subsub}\label{K_neq3} Le cas \m{n=3} -- \rm Les automorphismes
\m{\phi_{\mu\nu}} d\'ependent de param\`etres \m{\mu,\nu\in A\DB{t}/(t^2)},
n'ayant donc que deux coordonn\'ees dans \m{A\DB{x}}.

{\em Produit : } on a \ \m{\phi_{\mu',\nu'}\circ\phi_{\mu,\nu}=
\phi_{\mu''\nu''}}, avec
\[\mu''_0=\mu'_0+\mu_0\nu'_0, \ \ \ \mu''_1=\mu'_1+\mu_1{\nu'}_0^2+
\frac{d\mu_0}{dx}\mu'_0\nu'_0+\mu_0\nu'_1, \ \ \ \nu''_0=\nu_0\nu'_0, \ \ \
\nu''_1=\nu_0\nu'_1+\nu_1{\nu'}_0^2+\frac{d\nu_0}{dx}\mu'_0\nu'_0 .\]

{\em Inversion : } on a \ \m{\phi_{\mu\nu}^{-1}
=\phi_{\mu'\nu'}} , avec
\[\nu'_0=\frac{1}{\nu_0} , \ \ \ \mu'_0=-\frac{\mu_0}{\nu_0} , \ \ \
\nu'_1=\frac{1}{\nu_0^3}(\mu_0\frac{d\nu_0}{dx}-\nu_1) , \ \ \
\mu'_1=\frac{1}{\nu_0^3}\big(\mu_0(\nu_1-\mu_0\frac{d\nu_0}{dx}
+\nu_0\frac{d\mu_0}{dx})-\mu_1\big
) .\]

{\em D\'erivations et automorphismes associ\'es : } Les \'el\'ements de
\m{{\rm Der}_0(A\DB{x,t}/(t^3))} sont de la forme
\[D \ = \ at\DIV{}{x}+bt^2\DIV{}{t} ,\]
avec \m{a\in A\DB{t}/(t^2)}, \m{b\in A}. On a \ \m{\chi_D=\phi_{\mu,\nu}}, avec
\[\mu=a_0+\big(a_1+\frac{a_0}{2}(\frac{da_0}{dx}+b)\big)t , \ \ \ \
\nu=1+bt .\]
La loi de groupe sur \m{{\rm Der}_0(A\DB{x,t}/(t^3))} induite par celle sur
\m{\kg_0(A\DB{x,t}/(t^3))} est donn\'ee par : si \m{D=at\DIV{}{x}+
bt^2\DIV{}{t}}, \m{D'=a't\DIV{}{x}+b't^2\DIV{}{t}}, alors on a
\m{D'*D=a''t\DIV{}{x}+b''t^2\DIV{}{t}}, avec
\[a''_0=a_0+a'_0 , \ \ \ \ a''_1=a_1+a'_1+\frac{1}{2}(\frac{da_0}{dx}a'_0-
a_0\frac{da'_0}{dx}+a_0b'-ba'_0) , \ \ \ \ b''=b+b' .\]
On a \ \m{D*D'+D'*D=2(D+D')} (cette relation n'est plus vraie si \m{n>3}).
\end{subsub}

\end{sub}

\sepsub

\Ssect{Un r\'esultat de r\'egularit\'e}{regul}

On travaille ici dans l'anneau \m{\C((x))[[t]]}. Le r\'esultat suivant est
utilis\'e dans \ref{ecl} :

\sepprop

\begin{subsub}\label{reg_0}{\bf Proposition : } Soient \m{p>0} un entier,
\m{\mu,\nu,\mu_0,\nu_0\in\C\DB{x,t}} avec \m{\nu,\nu_0} inversibles. Alors
\[\chi \ = \
\phi_{\mu,x^p\nu}\circ\phi_{\mu_0\nu_0}\circ\phi_{\mu,x^p\nu}^{-1}\]
est un automorphisme de \m{\C\DB{x,t}}, c'est-\`a-dire qu'on peut \'ecrire
\m{\chi=\phi_{\mu_1\nu_1}}, avec \m{\mu_1,\nu_1\in\C\DB{x,t}} et \m{\nu_1}
inversible.
\end{subsub}

\begin{proof}
Soit \m{\nu'\in\C\DB{x,t}}. Alors on a
\[\phi_{\mu,\nu'}(x^p) \ = \ \lambda x^p \ ,\]
avec \m{\lambda\in\C\DB{x,t}} inversible et ind\'ependant de \m{\nu'}. On a
\[\phi_{\mu,\nu'}(x^pt) \ = \ x^p\lambda\nu't \ ,\]
donc si on prend \m{\nu'=\nu/\lambda}, on a
\begin{equation}\label{ecl_equ}
\phi_{\mu,\nu'}\circ\phi_{0,x^p} \ = \ \phi_{\mu,x^p\nu} \ ,
\end{equation}
et il suffit de montrer que
\[\chi_0 \ = \
\phi_{0,x^p}\circ\phi_{\mu_0\nu_0}\circ\phi_{0,x^p}^{-1}\]
est un automorphisme de \m{\C\DB{x,t}}.

On a \m{\phi_{0,x^p}^{-1}=\phi_{0,x^{-p}}}, et d'apr\`es la proposition
\ref{prox0}
\[\phi_{0,x^p}\circ\phi_{\mu_0\nu_0} \ = \ \phi_{\sigma,\tau} \ ,\]
avec \ \m{\sigma=x^p\phi_{0,x^p}(\mu_0)}, \m{\tau=x^p\phi_{0,x^p}(\nu_0)}. Si
\m{\alpha\in\C\DB{x,t}}, on a \ \m{\chi_0(\alpha)=\phi_{\sigma,\tau}(\alpha)\in
\C\DB{x,t}}. D'autre part, on a
\begin{eqnarray*}
\chi_0(t) & = & \phi_{\sigma,\tau}(x^{-p}t) \\
 & = & \tau\sigg_{i\geq 0}\frac{1}{i!}{\sigma}^i\frac{d^i(x^{-p})}{dx^i}t^i\\
 & = & \nu_1t \ ,
\end{eqnarray*}
avec
\[\nu_1 \ = \ \phi_{0,x^p}(\nu_0).\sigg_{i\geq 0}(-1)^i\begin{pmatrix}p-1\\ i
\end{pmatrix}x^{(p-1)i}\phi_{0,x^p}(\mu_0)^it^i \ ,\]
qui est inversible. Donc \m{\chi_0} est bien un automorphisme de \m{\C\DB{x,t}}.
\end{proof}

\end{sub}

\sepsub

\Ssect{Faisceaux de groupes d'automorphismes - La suite exacte de base}{g_aut}

Soient $C$ une courbe alg\'ebrique projective irr\'eductible lisse et \m{n>0}
un entier. On note \m{\ka_n=\ko_{C\times{\bf Z}_n}}. Autrement dit le faisceau
d'anneaux $\ka_n$ sur $C$ est d\'efini par :
\[\ka_n(U) \ = \ \ko_C(U)\ot_\C\C[t]/(t^n)\]
pour tout ouvert $U$ de $C$. On d\'efini aussi le faisceau \m{\ka_\infty} par
\[\ka_\infty(U) \ = \ko_C(U)\ot_\C\C[[t]] .\]
On a des morphismes canoniques de faisceaux d'anneaux si \m{n\geq 1},
\[r_n:\ka_n\to\ka_{n-1}, \quad {\text et} \quad
r_n^\infty:\ka_\infty\to\ka_{n-1} .\]
En particulier on a un morphisme canonique \m{\ka_n\to\ko_C} (pour \m{n\in\N} ou
\m{n=\infty}).

On note $\kg_n$ le faisceau de groupes des automorphismes de \m{C\times{\bf
Z}_n} laissant $C$ invariante. Plus explicitement, pour tout $n$ (fini ou non),
le faisceau de groupes \m{\kg_n} sur $C$ est d\'efini par~: pour tout ouvert
propre \m{U\subset C}, \m{\kg_n(U)} est le groupe des automorphismes de
$\C$-alg\`ebres \m{\theta:\ka_n(U)\to\ka_n(U)} tels que le triangle ci-dessous
soit commutatif
\xmat{
\ka_n(U)\ar[rr]^-\theta\ar[rd] & & \ka_n(U)\ar[ld]\\
& \ko_C(U)
}
Si $n$ est fini et \m{n\geq 2}, tout \'el\'ement de \m{\kg_n(U)} laisse
invariant \m{\ker(\rho_n(U))}, et induit donc un \'el\'ement de
\m{\kg_{n-1}(U)}. On obtient donc des morphismes canoniques
\[\rho_n:\kg_n\to\kg_{n-1}, \quad {\text et} \quad
\rho_n^\infty:\kg_\infty\to\kg_{n-1} .\]

Il d\'ecoule de \ref{ext_v} que si \m{\omega_{C\mid U}} est
trivial, alors tous les \'el\'ements de \m{\kg_n(U)} sont de la forme
\m{\phi_{\mu\nu}}, avec \m{\mu,\nu\in\ka_{n-1}(U)} et $\nu$ inversible (la
d\'efinition de \m{\phi_{\mu\nu}} d\'epend du choix d'une trivialisation de
\m{\omega_{C\mid U}}).

\sepsubsub

\begin{subsub}\label{ker_rho}Le noyau de \m{\rho_n} -- \rm On suppose \m{n\geq
3} fini. Soient \m{\mu,\nu\in\ka_{n-1}(U)} avec $\nu$ inversible. Alors on a
\[\rho_n(\phi_{\mu\nu}) \ = \ \phi_{r_{n-1}(\mu)r_{n-1}(\nu)} .\]
De m\^eme, on a \m{\rho_n^\infty(\phi_{\mu\nu})=
\phi_{r_{n-1}^\infty(\mu)r_{n-1}^\infty(\nu)}} si \m{\mu,\nu\in\ka_\infty(U)}
avec $\nu$ inversible. On en d\'eduit que les morphismes \m{\rho_n} et
\m{\rho_n^\infty} sont surjectifs.

Si \m{\omega_{C\mid U}} est trivial, alors, relativement \`a 
une section \m{\sigma=\frac{d}{dx}} engendrant  \m{T_{C\mid U}},
\m{\ker(\rho_n(U))} est constitu\'e des \m{\phi_{\mu\nu}}, avec \m{\mu,\nu} de
la forme
\[\mu=\theta t^{n-2}, \quad \nu=1+\beta t^{n-2} ,\]
avec \m{\mu,\nu\in\ko_C(U)}. On pose
\[\lambda_{\theta\beta}^\sigma = \phi_{\theta t^{n-2},1+\beta t^{n-2}}
.\]
On utilisera aussi la notation \m{\lambda_{\theta\beta}} s'il n'y a pas
d'ambigu\"it\'e. Si \m{P\in C }, on emploie la m\^eme notation pour les
\'el\'ements de \m{\ker(\rho_{nP})} (dans ce cas \m{\theta,\beta\in\ko_{CP}}).
\end{subsub}

\sepprop

\begin{subsub}\label{prox2} {\bf Proposition : } On a
\m{\ker(\rho_n)\simeq T_C\oplus\ko_C}.
\end{subsub}

\begin{proof} Soit $U$ un ouvert non vide de $C$ tel que \m{\omega_{C\mid U}}
soit trivial. Soit \m{\sigma=\frac{d}{dx}} une section engendrant
\m{T_{C\mid U}}. On voit ais\'ement en utilisant la proposition \ref{prox0}
ou directement que
\[\lambda_{\theta'\beta'}^\sigma\circ\lambda_{\theta\beta}^\sigma=
\lambda_{\theta'+\theta,\beta'+\beta}^\sigma .\]
D'autre part, si on remplace $\sigma$ par \m{a\sigma}, o\`u $a$ est une
fonction r\'eguli\`ere inversible sur $U$, on a d'apr\`es le lemme \ref{lemx1}
\[\lambda_{\theta\beta}^{a\sigma}=\lambda_{a\theta,\beta}^\sigma .\]
Il en d\'ecoule que \m{\ker(\rho_n)} est isomorphe au produit du sous-groupe
correspondant aux \m{\lambda_{\theta 0}^\sigma} et du sous-groupe correspondant
aux \m{\lambda_{0 \beta}^\sigma}, et que le premier (resp. le second)
sous-groupe est isomorphe \`a \m{T_C} (resp. \m{\ko_C}).
\end{proof}

\sepprop

On a donc pour \m{n\geq 3} une suite exacte
\begin{equation}
\xymatrix{
0\ar[r] & T_C\oplus\ko_C\ar[r] & \kg_n\ar[r]^-{\rho_n} & \kg_{n-1}\ar[r] &
0 .}
\end{equation}

\sepsubsub

\begin{subsub}\label{mor_1}Le morphisme canonique \m{\kg_n\to\ko_C^*} -- \rm 
C'est la version globale du morphisme analogue local \m{\xi_{\C\DB{x}}}
de \ref{chgt} (cf. aussi \ref{ext_v}). On suppose que \m{n\geq 2}. Pour tout
ouvert non vide $U$ de $C$ et tout \m{\phi\in\kg_n(U)} il existe un unique
\m{\nu\in\ka_n(U)} inversible tel que \m{\phi(t)=\nu t}. On pose
\m{\xi_n(\phi)=\rho(\nu)}. On d\'efinit ainsi un morphisme de faisceaux de
groupes surjectif \m{\xi_n:\kg_n\to\ko_C^*}. Si \m{n\geq 3} on a
\[\xi_{n-1}\circ\rho_n=\xi_n .\]

L'application induite
\[H^1(\xi_n) : H^1(X,\kg_n)\lra H^1(X,\ko_C^*)=\Pic(C)\]
est surjective. Soit \m{L\in\Pic(C)}, repr\'esent\'e par un 1-cocycle
\m{(s_{ij})} d'un recouvrement ouvert \m{(U_i)} de $C$. On a donc
\m{s_{ij}\in\C^*}. On en d\'eduit le 1-cocycle \m{(\phi_{0,s_{ij}})} \`a
valeurs dans \m{g_n}, repr\'esentant un \'el\'ement de \m{H^1(X,\kg_n)} qui ne
d\'epend que de $L$, et qu'on notera donc \m{e_n(L)}. Alors on a
\[H^1(\xi_n)(e_n(L))=L .\]
\end{subsub}

\sepsubsub

\begin{subsub}\label{g2}Description de \m{\kg_2} -- \rm Le noyau de \m{\xi_2}
s'identifie \`a \m{T_C}, compte tenu de la description des fibres de \m{\kg_2}
: soient \m{P\in C} et \m{x\in\ko_{CP}} un g\'en\'erateur de
l'id\'eal maximal. Alors \m{\kg_{2P}} est le groupe des matrices
\[\phi_{\mu\nu} = \begin{pmatrix}1 & 0\\ \mu\frac{d}{dx}& \nu\end{pmatrix}\quad,
\]
avec \m{\mu,\nu\in\ko_{CP}}, $\nu$ inversible. On a donc une suite exacte
\xmat{
0\ar[r] & T_C\ar[r] & \kg_2\ar[r]^-{\xi_2} & \ko_C^*\ar[r] &
0 .}
Notons qu'il existe une section canonique de \m{\xi_2} : le morphisme
\[{\xymatrix@R=2pt{
\tau:\ko_C^*\ar[r] & \kg_2 \\
\ \ \ \ \nu\fmaps[r] & {\begin{pmatrix}1 & 0\\ 0 & \nu\end{pmatrix}}\quad ,}}\]
mais le sous-groupe \m{\ko_C^*} de \m{\kg_2} n'est pas distingu\'e.
\end{subsub}

\end{sub}

\sepsub

\Ssect{Suites exactes d\'eriv\'ees}{g_aut2}

Soit \m{g\in H^1(C,\kg_n)}. Soit \m{L=H^1(\xi_n)(g)\in\Pic(C)} .

\sepsubsub

\begin{subsub}\label{neq2} Le cas \m{n=2} -- \rm D'apr\`es \ref{s_ex0} on a une
suite exacte
\xmat{0\ar[r] & T_C^g\ar[r] & \kg_2^g\ar[r] & (\ko_C^*)^g\ar[r] & 0 .}
\end{subsub}

\sepprop

\begin{subsub}\label{g_aut2_l0}{\bf Lemme : } 1 -- On a \m{T_C^g=T_C\ot L} et
\m{(\ko_C^*)^g=\ko_C^*}, donc on a une suite exacte
\xmat{0\ar[r] & T_C\ot L\ar[r] & \kg_2^g\ar[r] & \ko_C^*\ar[r] & 0 .}

2 -- L'action de \m{\C^*=H^0(\ko_C^*)} sur \m{H^1(T_C\ot L)} est la
multiplication.
\end{subsub}

\begin{proof} Montrons d'abord que \m{T_C^g=T_C\ot L}. Soient \m{P\in C} et
\m{x\in\ko_{CP}} un g\'en\'erateur de l'id\'eal maximal. Le r\'esultat va
d\'ecouler de la formule
\[\begin{pmatrix}1 & 0\\ \mu\frac{d}{dx}& \nu\end{pmatrix}
\begin{pmatrix}1 & 0\\ \mu_0\frac{d}{dx}& \nu_0\end{pmatrix}
{\begin{pmatrix}1 & 0\\ \mu\frac{d}{dx}& \nu\end{pmatrix}}^{-1} \ = \
\begin{pmatrix}1 & 0\\ (\mu+\nu\mu_0-\mu\nu_0)\frac{d}{dx}& \nu_0\end{pmatrix}
\quad ,\]
pour tous \m{\mu,\mu_0,\nu,\nu_0\in\ko_{CP}} avec $\nu$, \m{\nu_0} inversibles.
Supposons $g$ d\'efini par un 1-cocycle \m{\biggl(\begin{pmatrix}1 & 0\\
\mu_{ij} & \nu_{ij}\end{pmatrix}\biggr)} d'un recouvrement \m{\ku=(U_i)} de
$C$ tel que la restriction de \m{T_C} \`a chaque \m{U_i} soit triviale, avec
\m{\mu_{ij}\in H^0(U_{ij},T_C)}, \m{\nu_{ij}\in\ko_C^*(U_{ij})}. Alors d'apr\`es
la formule pr\'ec\'edente, \m{(T_C)^g} est obtenu en recollant les \m{T_{C\mid
U_i}} au moyen des multiplications par \m{\nu_{ij}} sur \m{U_{ij}}. D'apr\`es
\ref{tors_lin} le r\'esultat obtenu est \m{T_C\ot L}. 

La formule pr\'ec\'edente montre aussi que \m{\kg_2} agit trivialement sur
\m{\ko_C^*}, d'o\`u \m{(\ko_C^*)^g=\ko_C^*}. La partie 2- d\'ecoule aussi
imm\'ediatement de la formule pr\'ec\'edente.
\end{proof}

\sepprop

\begin{subsub}\label{fib_x2} Les fibres de \m{H^1(\xi_2)} -- \rm
Soit \m{g_0=H^1(\tau)(L)\in H^1(C,\kg_2)}. On a d'apr\`es le lemme
\ref{g_aut2_l0} et \ref{fib_h1p} une application surjective
\[\lambda_{g_0} : H^1(T_C\ot L)\lra H^1(\xi_2)^{-1}(L)\]
dont les fibres sont les orbites de l'action de \m{\C^*} par multiplication.
Avec les notations du lemme \ref{g_aut2_l0}, $L$ est d\'efini par
\m{(\nu_{ij})} et \m{g_0} l'est par le 1-cocycle \m{\biggl(\begin{pmatrix}1 &
0\\0 & \nu_{ij}\end{pmatrix}\biggr)}.
Compte tenu de l'\'egalit\'e
\[\begin{pmatrix}1 & 0\\ \mu_{ij} & 1\end{pmatrix}
\begin{pmatrix}1 & 0\\ 0 & \nu_{ij}\end{pmatrix} \ = \
\begin{pmatrix}1 & 0\\ \mu_{ij} & \nu_{ij}\end{pmatrix}
\quad ,\]
on voit que \m{\lambda_{g_0}^{-1}(g)=\C^*u}, o\`u \m{u\in H^1(T_C\ot L)} est
d\'efini par \m{(\mu_{ij})} (les relations de cocycle \'etant ici bien s\^ur
celles d\'ecrites dans \ref{act_fa}). On notera
\begin{equation}\label{fib_x2b} 0\lra T_C\lra E(g)\lra L^*\lra 0\end{equation}
l'extension associ\'ee \`a $u$ (qui ne d\'epend que de $g$).
\end{subsub}

\sepsubsub

\begin{subsub}\label{cas_n} Le cas \m{n>2} -- \rm
On suppose maintenant que \m{n\geq 3}. D'apr\`es \ref{s_ex0} on a une suite
exacte
\xmat{
0\ar[r] & (T_C\oplus\ko_C)^g\ar[r] & \kg_n^g\ar[r] & \kg_{n-1}^g\ar[r] &
0 .}
Soit \m{g_2} l'image de $g$ dans \m{H^1(C,\kg_2)}.
\end{subsub}

\sepprop

\begin{subsub}\label{g_aut2_l1}{\bf Lemme : } Soient \m{P\in C},
\m{\theta,\beta\in\ko_{CP}}, \m{\mu,\nu\in\ka_{n-1,P}} avec $\nu$ inversible,
et \m{\mu_0=\rho(\mu)}, \m{\nu_0=\rho(\nu)}. Alors on a
\[\phi_{\mu\nu}\circ\lambda_{\theta\beta}\circ\phi_{\mu\nu}^{-1} \ =  \
\lambda_{\theta'\beta'} ,\]
avec
\[\beta'=\nu_0^{n-2}\beta ,\quad \theta'=\nu_0^{n-2}(\nu_0\theta-\mu_0\beta) .\]
\end{subsub}

\begin{proof} 
Cela d\'ecoule ais\'ement de la proposition \ref{prox0}.
\end{proof}

\sepprop

\begin{subsub}\label{g_aut2_p1}{\bf Proposition : } On a \
\m{(T_C\oplus\ko_C)^g\simeq E(g_2)\ot L^{n-1}} .
\end{subsub}

\begin{proof} Soit \m{\biggl(\begin{pmatrix}1 & 0\\ \mu_{ij} &\nu_{ij}
\end{pmatrix}\biggr)} un 1-cocycle de $\ku$ repr\'esentant \m{g_2}. D'apr\`es
le lemme pr\'ec\'edent, \m{(T_C\oplus\ko_C)^g} s'obtient en recollant les
faisceaux \m{(T_C\oplus\ko_C)_{\mid U_i}} au moyen des automorphismes sur
\m{U_{ij}} correspondant aux matrices \m{\biggl(\begin{pmatrix}\nu_{ij}^{n-1} &
-\nu_{ij}^{n-2}\mu_{ij}\\ 0 &\nu_{ij}^{n-2}\end{pmatrix}\biggr)} . Le
r\'esultat d\'ecoule donc de \ref{desc_ext}.
\end{proof}

\sepprop

\begin{subsub}\label{fib_gn} Les fibres de \m{H^1(\rho_n)} -- \rm
Soit \m{g_{n-1}=H^1(\rho_n)(g)}. 0n a d'apr\`es la proposition \ref{g_aut2_p1}
et \ref{fib_h1p} une application surjective
\[\lambda_{g} : H^1(E(g_2)\ot L^{n-1})\lra H^1(\rho_n)^{-1}(g_{n-1})\]
envoyant 0 sur $g$ et dont les fibres sont les orbites de l'action de
\m{H^0(C,\kg_{n-1}^g)}.
\end{subsub}

\sepprop

\begin{subsub}\label{surj_gn} Surjectivit\'e de \m{H^1(\rho_n)} -- \rm
Il d\'ecoule de \ref{coh_4} que {\em \m{H^1(\rho_n)} est surjective}.
\end{subsub}

\end{sub}

\sepsec

\section{Courbes multiples primitives abstraites}\label{HMPA}

On utilisera les notations de \ref{fg_aut}. Rappelons que pour \m{n\geq
1}, \m{{\mathbf Z}_n} d\'esigne le sch\'ema \Nligne\m{\spec{(\C[t]/(t^n))}}.

\sepsub

\Ssect{D\'efinition et propri\'et\'es \'el\'ementaires}{cpa_def}

Soit $C$ une courbe irr\'eductible lisse. Soient \m{n\geq 2} un
entier et \m{g\in H^1(C,\kg_n)}, repr\'esent\'e par un 1-cocycle \m{(g_{ij})}
d'un recouvrement ouvert \m{\ku=(U_i)} de $C$. On peut voir \m{g_{ij}} comme un
automorphisme de \m{U_{ij}\times{\mathbf Z}_n} laissant \m{U_{ij}} invariant.
On note \m{C(g)} le sch\'ema obtenu en recollant les \m{U_i\times{\mathbf Z}_n}
au moyen des \m{g_{ij}}. Il est facile de voir que \m{C(g)} ne d\'epend
effectivement que de $g$. Si $C$ est projective c'est un sch\'ema propre et de
dimension 1, c'est donc une vari\'et\'e projective dont la sous-vari\'et\'e
r\'eduite associ\'ee est $C$. 

Pour \m{1\leq k\leq n}, la sous-vari\'et\'e ferm\'ee
\m{U_{ij}\times{\mathbf Z}_k} de \m{U_{ij}\times{\mathbf Z}_n} est invariante
par \m{g_{ij}}. Ces sous-vari\'et\'es se recollent donc et d\'efinissent une
sous-vari\'et\'e \m{C_k(g)} de \m{C(g)}. On a
\[C=C_1(g)\subset C_2(g)\subset\cdots\subset C_{n-1}(g)\subset C_n(g)=C(g) .\]
Si \m{\ki_C} d\'esigne le faisceau d'id\'eaux de $C$ dans \m{C(g)}, celui de
\m{C_k(g)} est \m{\ki_C^k}.

On appelle \m{C(g)} une {\em courbe multiple primitive abstraite} de courbe
r\'eduite associ\'ee $C$, de multiplicit\'e $n$ et de fibr\'e en droites
associ\'e $L$.

\sepprop

\begin{subsub}\label{cpa_def1}{\bf Proposition : }
1 -- Si \m{n'\geq 2} est un entier et si \m{g'\in
H^1(C,\kg_{n'})}, alors il existe un isomorphisme \m{C(g)\simeq C(g')}
induisant l'identit\'e sur $C$ si et seulement si \m{n=n'} et \m{g=g'}.

2 -- On a \m{C(g_k)=C_k(g)}.

3 -- Soit \m{L=H^1(\xi_n)(g)\in\Pic(C)} . Alors on a \ \m{\ki_C/\ki_C^2
\simeq L} .

4 -- Toute courbe primitive abstraite de multiplicit\'e $n$ peut \^etre
\'etendue en une courbe primitive abstraite de multiplicit\'e \m{n+1} de m\^eme
courbe r\'eduite associ\'ee.
\end{subsub}

\sepprop

\begin{proof}
Les trois premi\`eres assertions sont imm\'ediates et la quatri\`eme d\'ecoule
de \ref{surj_gn}.
\end{proof}

\sepprop

Le faisceau de groupes \m{(\kg_n)^g} est canoniquement isomorphe au faisceau
\m{\AAut_C(C(g))} des automorphismes de \m{C(g)} induisant l'identit\'e sur
$C$. En particulier \m{H^0(C,(\kg_n)^g)} s'identifie au groupe \m{\Aut_C(C(g))}
des automorphismes de \m{C(g)} induisant l'identit\'e sur $C$.

Si \m{L\in\Pic(C)}, la courbe \m{C(e_n(L))} s'appelle la {\em courbe triviale}
de multiplicit\'e $n$ et de fibr\'e en droites associ\'e $L$ (cf. \ref{mor_1}).
On en donnera une description plus pr\'ecise en \ref{pr0_c}.

\sepprop

\begin{subsub}\label{cpa_def2}{\bf Proposition : } Si $C$ n'est pas projective,
les seules courbes multiples primitives abstraites de courbe r\'eduite
associ\'ee $C$ sont les courbes triviales.
\end{subsub}

\sepprop

\begin{proof}
Cela d\'ecoule de \ref{fib_x2}, \ref{fib_gn} et du fait que $C$ est affine (cf.
\cite{ha}, chapter IV, ex. 1.4).
\end{proof}

\sepprop

\begin{subsub}\label{cpa_def3} D\'efinition de faisceaux coh\'erents par
recollements -- \rm On note \m{{\bf U}_i} l'ouvert de \m{C(g)} correspondant
\`a \m{U_i\times{\bf Z}_n} (et qui lui est donc isomorphe), et \m{\pi_i:{\bf
U}_i\to U_i\times{\bf Z}_n} l'isomorphisme canonique. On a donc
\m{g_{ij}=\pi_i\circ\pi_j^{-1}} sur \m{U_{ij}\times{\bf Z}_n}.
On se donne des faisceaux coh\'erents \m{\ke_i} sur
\m{U_i\times{\bf Z}_n} et des isomorphismes de faisceaux sur
\m{U_{ij}\times{\bf Z}_n} :
\[\Theta_{ij} : \ke_{j} \ \simeq \ g_{ij}^*(\ke_{i}) \]
tels que pour tous indices distincts deux \`a deux $i$, $j$, $k$ le carr\'e
suivant soit commutatif sur \m{U_{ijk}\times{\bf Z}_n} :
\xmat{\ke_k\ar[r]^-{\Theta_{jk}}\ar[d]^-{\Theta_{ik}} &
g_{jk}^*(\ke_j)\ar[d]^-{g_{jk}^*(\Theta_{ij})} \\
g_{ik}^*(\ke_i)\fleq[r] & g_{jk}^*\big(g_{ij}^*(\ke_i)\big)}
Alors il existe un unique faisceau coh\'erent $\ke$ sur \m{C(g)} tel qu'on ait
des isomorphismes \m{\gamma_i:\ke_{\mid{\bf U}_i}\to\pi_i^*(\ke_i)} tels que
\m{\Theta_{ij}=(\pi_j^{-1})^*(\gamma_i\gamma_j^{-1})} .

Le faisceau \m{\kf=g_{ij}^*(\ke_{i})} peut \^etre d\'ecrit de la fa\c con
suivante : pour tout ouvert $U$ de \m{U_i\times{\bf Z}_n} on a
\m{\kf(U)=\ke_i(U)}, et la structure de module est la suivante : pour tout
\m{\alpha\in\ko_{U_i\times{\bf Z}_n}(U)} et tout \m{v\in\ke_i(U)} le produit
de $v$ par $\alpha$ est en fait \m{(\alpha\circ g_{ij}).v}.
\end{subsub}

\end{sub}

\sepsub

\Ssect{Courbes multiples primitives ordinaires}{courb_o}

Soient $X$ une vari\'et\'e irr\'eductible lisse et \m{Y\subset X}
une hypersurface lisse. On note \m{Y_n} l'hypersurface de multiplicit\'e $n$
associ\'ee et \m{\ko_n} son faisceau structural.

\sepprop

\begin{subsub}{\bf Th\'eor\`eme : }\label{pr0} Il existe un morphisme de
\m{\C}-alg\`ebres \m{\ko_{Y,P}\to\ko_{n,P}} section de la restriction
\m{\ko_{n,P}\to\ko_{Y,P}}. On a donc un isomorphisme de \m{\C}-alg\`ebres \
\m{\ko_{n,P}\simeq\ko_{Y,P}\ot_\C\C[t]/(t^n)}.
\end{subsub}

\begin{proof}
Supposons $X$ plong\'ee dans \m{\P_m}, \m{P=(1,0,\ldots,0)}, $X$ d\'efinie par
les \'equations 
\[f_1(x_1,\ldots,x_m)=\cdots=f_{m-d}(x_1,\ldots,x_m)=0\]
(o\`u \m{d=\dim(X)}) et $Y$ par l'\'equation suppl\'ementaire
\m{f(x_1,\ldots,x_m)=0} au voisinage de $P$. On a alors
\[\ko_{X,P}=\ko_P/(f_1,\ldots,f_{m-d}),\quad\ko_{n,P}=\ko_P/(f_1,\ldots,
f_{m-d},f^n),\quad\ko_{Y,P}=\ko_P/(f_1,\ldots,f_{m-d},f) ,\]
avec \m{\ko_P=\ko_{\P_m,P}}, qu'on peut voir comme la sous-alg\`ebre de
\m{\C(x_1,\ldots,x_m)} constitu\'ee des \'el\'ements \m{A/B}, avec \m{A,B\in
\C[x_1,\ldots,x_m]} et \m{B(0)\not=0}. 

La proposition \ref{pr0} se d\'emontre par r\'ecurrence sur $n$. Le cas \m{n=1}
est \'evident. Supposons qu'elle soit vraie pour \m{n-1\geq 1}.
Un morphisme de \m{\C}-alg\`ebres \m{\phi:\ko_{Y,P}\to\ko_{n,P}} section de la
restriction \m{\ko_{n,P}\to\ko_{Y,P}} est d\'efini par \m{\phi(x_i)}, \m{1\leq
i\leq m}, qui doit \^etre de la forme
\[\phi(x_i) \ = \ x_i+A_if ,\]
avec \m{A_i\in \ko_P}. On doit avoir de plus \ \m{f_j(\phi(x_1),\ldots,
\phi(x_m))=0} \ pour \m{1\leq j\leq m-d} \ et \
\m{f(\phi(x_1),\ldots,\phi(x_m))=0} \ dans \m{\ko_{n,P}}. 
La formule de Taylor donne
\[f_j(x_1+A_1f,\ldots,x_m+A_mf) \ = \ \sigg_{i=1}^{n-1}\biggl(f^i\sigg
C_{\alpha_1\ldots\alpha_m}\frac{\partial^i f_j}{\partial
x_{\alpha_1}\cdots\partial x_{\alpha_i}}
(x_1,\ldots,x_m)
A_{\alpha_1}\cdots A_{\alpha_i}\biggr) \ , \]
\[f(x_1+A_1f,\ldots,x_m+A_mf) = f+\sigg_{i=1}^{n-1}\biggl(f^i\sigg
C_{\alpha_1\ldots\alpha_m}\frac{\partial^j f}{\partial
x_{\alpha_1}\cdots\partial x_{\alpha_i}}(x_1,\ldots,x_m)
A_{\alpha_1}\cdots A_{\alpha_i}\biggr) \ . \]
dans \m{\ko_{nP}} (les \m{C_{\alpha_1\ldots\alpha_m}} \'etant des coefficients
constants) . D'apr\`es l'hypoth\`ese de r\'ecurrence on peut supposer qu'on
peut \'ecrire
\[f_j(x_1+A_1f,\ldots,x_m+A_mf)=D_jf^{n-1} \ , \quad
f(x_1+A_1f,\ldots,x_m+A_mf) = Df^{n-1} \ ,\]
avec \m{D_j,D\in\ko_P}. Rempla\c cons maintenant \m{A_1,\ldots,A_m} par \
\m{A'_1=A_1+f^{n-2}B_1,\ldots},\Nligne\m{A'_m=A_m+f^{n-2}B_m} \ respectivement.
On a alors dans \m{\ko_{nP}} (c'est-\`a-dire modulo \m{(f^n)})
\[f_j(x_1+A'_1f,\ldots,x_m+A'_mf)=f^{n-1}.\biggl(D_j+
\sigg_{i=1}^m\frac{\partial f_j}{\partial x_i}(x_1,\ldots,x_m)B_i\biggr) ,\]
\[f(x_1+A'_1f,\ldots,x_m+A'_mf)=f^{n-1}.\biggl(D+
\sigg_{i=1}^m\frac{\partial f}{\partial x_i}(x_1,\ldots,x_m)B_i\biggr) .\]
Puisque $Y$ est lisse la matrice
\[\begin{pmatrix}
\frac{\partial f_1}{\partial x_1}(0) & \cdots &
\frac{\partial f_1}{\partial x_m}(0)\\
. & . & .\\
. & . & .\\
. & . & .\\
\frac{\partial f_{m-d}}{\partial x_1}(0) & \cdots &
\frac{\partial f_{m-d}}{\partial x_m}(0)\\
\\
\frac{\partial f}{\partial x_1}(0) & \cdots &
\frac{\partial f}{\partial x_m}(0)
\end{pmatrix}\]
est de rang maximal \m{m-d+1}. On peut donc choisir \m{B_1,\ldots,B_m} de telle
sorte que
\[D_j+\sigg_{i=1}^m\frac{\partial f_j}{\partial x_i}(x_1,\ldots,x_m)B_i
\ = \ 0 \ , \ \ \ \ D+\sigg_{i=1}^m\frac{\partial f}{\partial x_i}
(x_1,\ldots,x_m)B_i \ = 0 \ .\]
On a alors
\[f_j(x_1+A'_1f,\ldots,x_m+A'_mf)=0 \ , \ \ \ \ f(x_1+A'_1f,\ldots,x_m+A'_mf)
=0\]
dans \m{\ko_{nP}}.
\end{proof} 

\sepprop

En modifiant l\'eg\`erement la d\'emonstration pr\'ec\'edente on obtient le

\sepprop

\begin{subsub}{\bf Corollaire : }\label{cor1_pr0} Il existe un voisinage $U$ de
$P$ dans $Y$ tel que si \m{U_n} d\'esigne l'ouvert correspondant de \m{Y_n} on
ait un diagramme commutatif
\xmat{ & U\flinc[ld]\flinc[rd] \\
U_n\ar[rr]^-\simeq & & U\times {\mathbf Z}_n}
Autrement dit, \m{Y_n} est localement triviale.
\end{subsub}

\sepprop

\begin{subsub}{\bf Corollaire : }\label{cor2_pr0} Toute courbe multiple
primitive au sens de \ref{cmpr} est une courbe multiple primitive abstraite.
\end{subsub}

\sepprop

\begin{subsub}\label{pr0_b}{\bf Proposition : }
Soient \m{C_n} une courbe primitive de multiplicit\'e \m{n\geq 1}, $C$ la
courbe lisse r\'eduite, \m{C_2} la courbe double associ\'ees, et
\m{L=\ki_C/\ki_{C_2}\in\Pic(C)}. Soit \m{g\in H^1(C,\kg_n)} tel qu'il existe un
isomorphisme \m{C(g)\simeq C_n} induisant l'identit\'e sur $C$. Alors on a
\m{L=H^1(\xi_n)(g)}.
\end{subsub}

\begin{proof}
Soit \m{(g_{ij})} un 1-cocycle d'un recouvrement ouvert \m{(U_i)} de $C$
repr\'esentant $g$. Donc \m{g_{ij}} est un automorphisme de
\m{\ko_C(U_{ij})[t]/(t^n)} tel que \m{\rho g_{ij}=\rho}, $\rho$ d\'esignant le
morphisme canonique \m{\ko_C(U_{ij})[t]/(t^n)\to\ko_C(U_{ij})}. On a alors
\m{g_{ij}(t)=\tau_{ij}t}, avec \m{\tau_{ij}} inversible, donc
\m{\nu_{ij}=\rho(\tau_{ij})\in\ko_C^*(U_{ij})}. Le fibr\'e \m{H^1(\xi_n)(g)}
est d\'efini par le 1-cocycle \m{(\tau_{ij})}. D'apr\`es la construction de
\m{C(g)} il est imm\'ediat que \m{\ki_C/\ki_{C_2}} est d\'efini par le m\^eme
cocycle, d'ou le r\'esultat.
\end{proof}

\sepprop

On a en fait \m{\nu_{ij}\in\ko_C^*(U_{ij})[t]/(t^{n-1})}, et \m{(\tau_{ij})}
est un 1-cocycle d\'efinissant un fibr\'e en droites sur \m{C_{n-1}} isomorphe
\`a \m{\ki_C}.

\sepsubsub

\begin{subsub}\label{pr0_c}Description des courbes triviales -- \rm Soit
\m{L\in\Pic(C)}, et $S$ le fibr\'e en droites \m{L^*} sur $C$, vu comme une
surface. La section nulle de \m{L^*} d\'efinit un plongement \m{C\subset S}.
Soit \m{C_n} la courbe primitive de multiplicit\'e $n$ associ\'ee. Alors
\m{C_n} est la courbe triviale de multiplicit\'e $n$, c'est-\`a-dire que
\m{C_n\simeq C(e_n(L))} (cf. \ref{mor_1}).
\end{subsub}

\sepsubsub

\begin{subsub}\label{pr0_d}Param\'etrisation des courbes doubles -- \rm
Soit \m{C_2=C(g)} une courbe double de courbe lisse sous-jacente $C$ et de
fibr\'e en droites associ\'e $L$. D'apr\`es \ref{fib_x2}, \m{C_2} est triviale
si et seulement si \m{E(g)\simeq T_C\oplus L^*}, c'est-\`a-dire si et seulement
si la suite exacte (\ref{fib_x2b}) est scind\'ee. Les courbes doubles non
triviales de courbe lisse sous-jacente $C$ et de fibr\'e en droites associ\'e
$L$ sont naturellement param\'etr\'ees par \m{\P(H^1(T_C\ot L))}, le point
correspondant \`a \m{C(g)} \'etant associ\'e \`a l'extension (\ref{fib_x2b}).
On retrouve les r\'esultats de \cite{ba_ei}.
\end{subsub}

\sepsubsub

\begin{subsub}\label{pr0_f}Param\'etrisation des prolongements de courbes
multiples -- \rm Soient $n$ un entier tel que \m{n\geq 3}
,\m{C_{n-1}=C(\gamma)} une courbe multiple prmimitive de multiplicit\'e
\m{n-1} et \m{\gamma_2} l'image de $\gamma$ dans \m{H^1(C,\kg_2)}. On note
\m{\kc_n(C_{n-1})} l'ensemble des prolongements de \m{C_{n-1}} en courbe
primitive de multiplicit\'e $n$. Soient \m{C_n} un tel prolongement et $g$
l'unique \'el\'ement de \m{H^1(C,\kg_n)} tel que \m{C_n=C(g)}. D'apr�s
\ref{fib_gn} il existe une surjection canonique
\[\lambda_g : H^1(E(\gamma_2)\ot L^{n-1})\lra\kc_n(C_{n-1})\]
donc les fibres sont les orbites d'une action de \m{\Aut_C(C_{n-1})} sur
\m{H^1(C,\kg_n)} (cf. \ref{cas_com}, \ref{coh_5}).
\end{subsub}

\sepsubsub

\begin{subsub}\label{pr0_e}Cas d'unicit\'e des extensions de courbes -- \rm
Soit \m{C_n} une courbe primitive de multiplicit\'e $n$ de courbe lisse
projective sous-jacente $C$ et de fibr\'e en droites associ\'e $L$. Soient $g$
le genre de $C$ et \m{d=\deg(L)}. Il r\'esulte de \ref{fib_gn} qu'il existe un
{\em unique} prolongement de \m{C_n} en courbe de multiplicit\'e \m{n+1} dans
les cas suivants :
\begin{itemize}
 \item[--] Si $g\geq 1$ : $d>\frac{4g-4}{n}$, ou $d=\frac{4g-4}{n}$ et
$L^n\not=\omega_C^2$. En particulier, si $n>4g-4$ et $d>0$.
\item[--] Si $g=0$ : $n=2,3$ et $d\geq -1$, $n\geq 4$ et $d\geq 0$.
\end{itemize}
\end{subsub}

\end{sub}

\sepsub

\Ssect{Plongement des courbes multiples primitives abstraites}{plong}

Soit $C$ une courbe lisse plong\'ee dans une surface lisse $S$, elle-m\^eme
plong\'ee dans \m{\P_m}, \m{m\geq 4}. Soit \m{n\geq 2} un entier. On note
\m{C_n} la courbe primitive de $S$ de multiplicit\'e $n$ et de courbe r\'eduite
associ\'ee $C$. Si \m{P\in C} on note \m{T_{SP}} le plan de \m{\P_m} tangent
\`a $S$ en $P$. Rappelons qu'on appelle {\em s\'ecante} \`a $C$ une droite de
\m{\P_m} contenant au moins deux points distincts de $C$.

\sepprop

\begin{subsub}\label{plong1}{\bf Proposition : }
Soit \m{O\in\P_m} un point n'appartenant \`a aucun plan tangent \m{T_{SP}},
\m{P\in C}, et \`a aucune s\'ecante \`a $C$. Alors la projection de centre $O$,
\m{\pi_O :\P_m\to\P_{m-1}} induit un plongement \m{C_n\subset\P_{m-1}}.
\end{subsub}

\begin{proof}
On sait d\'ej\`a que \m{\pi_O} induit un plongement \m{C\subset\P_{m-1}} (cf.
\cite{ha}, Chap. IV, Prop. 3.4). Soit \m{P\in C}. Puisque \m{O\not\in T_{SP}}
il existe un voisinage $U$ (au sens de la topologie usuelle) de $P$ tel que
\m{\pi_O} induise un plongement \m{U\cap S\subset\P_{m-1}}. Donc
\m{\pi_O^*:\ko_{\P_{m-1},\pi_O(P)}\to\ko_{C_n,P}} est surjectif, et \m{\pi_O}
induit bien un plongement \m{C_n\subset\P_{m-1}}.
\end{proof}

\sepprop

\begin{subsub}\label{plong2}{\bf Th\'eor\`eme : }
Toute courbe multiple primitive abstraite peut \^etre plong\'ee dans \m{\P_3}.
\end{subsub}

\begin{proof}
Soit $D$ une courbe multiple primitive abstraite, de courbe r\'eduite
associ\'ee $C$. On part d'un plongement \m{D\subset\P_m}, \m{m\geq 4}. Il faut
prouver qu'il existe un point \m{O\in\P_m} tel que la projection de centre $O$,
\m{\pi_O :\P_m\to\P_{m-1}} induise un plongement \m{C_n\subset\P_{m-1}}.

D'apr\`es la structure locale de $D$ il existe des surfaces lisses
\m{S_1,\ldots,S_k\subset\P_m} (non n\'ecessairement projectives) telles que
\m{D\subset S_1\cup\cdots\cup S_k}. La vari\'et\'e des points situ\'es sur
les plans tangents \`a une des \m{S_i} en les points de $C$ est de dimension au
plus 3, de m\^eme que celle des points situ\'es sur les s\'ecantes. Il suffit
donc d'apr\'es la proposition \ref{plong1} de prendre $O$ en dehors de
l'adh\'erence de ces vari\'et\'es.
\end{proof}

\sepprop

On en d\'eduit en utilisant aussi le corollaire \ref{cor2_pr0} le

\sepprop

\begin{subsub}\label{plong3}{\bf Corollaire : }
Les courbes multiples primitives abstraites sont exactement les courbes
multiples primitives au sens de \ref{cmpr}.
\end{subsub}

\end{sub}

\sepsub

\Ssect{\'Eclatements}{ecl}

Soit \m{n\geq 2} un entier. Soit \m{C_n} une courbe multiple primitive de
multiplicit\'e $n$ et de courbe r\'eduite associ\'ee $C$. On peut supposer que
\m{C_n=C(g)}, avec \m{g\in H^1(C,\kg_n)}. Si \m{2\leq i\leq n} on notera
\m{g_i} l'image de \m{g_n} dans \m{H^1(C,\kg_i)}. On a donc \m{C(g_i)=C_i}.

Soit \m{\ku=(U_i)_{i\in I}} un recouvrement ouvert de $C$ On suppose que pour
tout indice $i$, \m{U_i\not=C} et que \m{\omega_{C\mid U_i}} est trivial. Pour
tout ouvert $U$ de $C$, on note $\rho$ le morphisme canonique
\m{\ko_C(U)[t]/(t^n)\to\ko_C(U)}.

Soit \m{(g_{ij})} un 1-cocycle de $\ku$ de $C$ repr\'esentant $g$. Donc
\m{g_{ij}} est un automorphisme de \m{\ko_C(U_{ij})[t]/(t^n)} tel que \m{\rho
g_{ij}=\rho}.

Soient \m{i_0\in I}, \m{P\in U_{i_0}} et \m{x\in\ko_C(U_{i_0})} engendrant
l'id\'eal maximal de \m{\ko_{CP}} et ne s'annulant qu'en $P$. On peut en
modifiant au besoin $\ku$ supposer que \m{i_0} est le seul indice $i$ tel que
\m{P\in U_i}. Soit $A$ un endomorphisme de \m{\ko_C(U_{i_0})[t]/(t^n)} tel que
\m{\rho A=\rho}. On suppose que $P$ est le seul point singulier de $A$,
c'est-\`a-dire que $A$ induit un automorphisme de
\m{\ko_C(U_{i_0}\backslash\lbrace P\rbrace)[t]/(t^n)}, et que \m{\ko_C(-P)} est
trivial sur \m{U_{i_0}}. On d\'efinit un nouveau 1-cocycle
\m{\gamma(A,(g_{ij}))=(g'_{ij})} de $\ku$ par
\[g'_{ij}=g_{ij} \quad {\rm si } \quad i\not=i_0 \ {\rm et } \ j\not=i_0 ,
\ \ \ \ g'_{i_0j}=A.g_{i_0j} ,\ \ \ \ g'_{ji_0}=g_{ji_0}.A^{-1} \ .\]
On note \m{g'} l'\'el\'ement de \m{H^1(C,\kg_n)} repr\'esent\'e par
\m{(g'_{ij})} (on verra plus loin que \m{g'} ne d\'epend que de $g$ et $q$),
et \m{C'_m} les courbes multiples correspondantes.

\sepprop

\begin{subsub}\label{ecl_1}{\bf Proposition : } On suppose que $A$ est de la
forme \m{A=\phi_{\mu,x^q\nu}}, avec \m{q>0} et $\nu$ inversible. Alors

1 - La courbe \m{C(g')} est l'\'eclatement du diviseur de Cartier \m{qP} de
\m{C(g)}.

2 - L'association \m{(g_{ij})\to \gamma(A,(g_{ij}))} induit un morphisme
surjectif
\[b_{n,q}^P : H^1(C,\kg_n)\to H^1(C,\kg_n)\]
ne d\'ependant que de $q$ et $P$.

3 - On a un carr\'e commutatif
\xmat{H^1(C,\kg_n)\ar[rr]^-{H^1(\xi_n)}\ar[d]^{b_{n,q}^P} & &
\Pic(C)\ar[d]^-{\ot\ko_C(qP)}\\
H^1(C,\kg_n)\ar[rr]^-{H^1(\xi_n)} & & \Pic(C)}

4 - On a un carr\'e commutatif
\xmat{H^1(C,\kg_n)\ar[rr]^-{H^1(\rho_n)}\ar[d]^{b_{n,q}^P} & &
H^1(C,\kg_{n-1})\ar[d]^{b_{n-1,q}^P}\\
H^1(C,\kg_n)\ar[rr]^-{H^1(\rho_n)} & & H^1(C,\kg_{n-1})}

5 - Pour tous entiers \m{q,q'>0} on a \ \m{b_{n,q}^P\circ
b_{n,q'}^P=b_{n,q+q'}^P} .
\end{subsub}

\begin{proof}
La d\'emonstration de 1- est analogue \`a celle de \cite{ba_ei}, Theorem 1.9.

L'asserion 2- d\'ecoule de 1-, mais on  peut la d\'emontrer directement.
Rempla\c cons \m{(g_{ij})} par un cocycle cohomologue \m{(h_ig_{ij}h_j^{-1})},
on a
\[A.h_{i_0}g_{i_0j}h_j^{-1}=(Ah_{i_0}A^{-1}).Ag_{i_0j} ,\]
et d'apr\`es la proposition \ref{reg_0}, \m{Ah_{i_0}A^{-1}} est un
automorphisme de \m{\ko_C(U_{i_0})[t]/(t^n)}, donc
\m{\gamma(A,(h_ig_{ij}h_j^{-1}))} est cohomologue \`a \m{\gamma(A,(g_{ij}))}.
Ceci permet de d\'efinir \m{b_{n,q}^P} sans ambiguit\'e. D'autre part la
relation (\ref{ecl_equ}) montre que \m{b_{n,q}^P} ne d\'epend que de $q$. Ceci
d\'emontre 2-.

Soit \m{\sigma\in\ko_C(U)} une \'equation de $P$. Posons
\m{g_{ij}=\phi_{\mu_{ij},\nu_{ij}}}. Alors \m{H^1(\rho_n)(g)} est
repr\'esent\'e par le cocycle \m{(\delta_{ij})=(\rho(\nu_{ij}))}, tandis que
\m{H^1(\rho_n)(b_{n,q}^P(g))} l'est par le cocycle \m{(\epsilon_{ij})} d\'efini
par
\[\epsilon_{ij}=\delta_{ij} \quad {\rm si } \quad i\not=i_0 \ {\rm et } \
j\not=i_0 , \ \ \ \ \epsilon_{i_0j}=\sigma^q\delta_{i_0j} ,\ \ \ \
\delta_{ji_0}=\sigma^{-q}\delta_{ji_0} \ .\]
On en d\'eduit ais\'ement 3-. Les autres assertions sont imm\'ediates.
\end{proof}

\sepprop

On s'int\'eresse maintenant \`a la correspondance obtenue entre les fibres de
\m{H^1(\rho_n)}.\Nligne Soit \m{L=\ki_C/\ki_{C_2}\in\Pic(C)}. D'apr\`es la 
partie 2- du th\'eor\`eme pr\'ec\'edent, on a \Nligne 
\m{\ki_C/\ki_{C'_2}=L(qP)}. Rappelons 
que le fibr\'e \m{E(g_2)} est une extension de \m{L^*} par \m{T_C}. Soit
\m{u\in\Ext^1_{\ko_C}(L^*,T_C)=H^1(T_C\ot L)} correspondant \`a cette extension.

\sepprop

\begin{subsub}\label{ecl_2}{\bf Proposition : } Le fibr\'e \m{E(g'_2)} est
associ\'e \`a l'image de $u$ dans \m{H^1(T_C\ot L(qP))} par le morphisme \
\m{H^1(T_C\ot L)\to H^1(T_C\ot L(qP))} induit par l'inclusion \m{L\subset
L(qP)}.
\end{subsub}

\begin{proof}
Ce r\'esultat est contenu dans \cite{ba_ei}, Theorem 1.9. On en donne une
autre d\'emonstration utilisant des cocycles.
On utilise les notations de \ref{ecl_1}. Alors \m{E(g_2)} est obtenu en
recollant les fibr\'es \m{(T_C\oplus\ko_C)_{\mid U_i}} au moyen des
automorphismes \m{\begin{pmatrix}0 & 1\\ \mu_{ij} & \nu_{ij}\end{pmatrix}}.
Le fibr\'e \m{E(g'_2)} lui est obtenu en recollant les fibr\'es
\m{(T_C\oplus\ko_C)_{\mid U_i}} au moyen des automorphismes \m{\begin{pmatrix}0
& 1\\ \mu_{ij} & \nu_{ij}\end{pmatrix}} si \m{i\not=i_0}, \m{j\not=i_0},
 \m{\begin{pmatrix}0 & 1\\ \sigma^q\mu_{ij} & \sigma^q\nu_{ij}\end{pmatrix}}
si \m{i=i_0}, et \m{\begin{pmatrix}0 & 1\\ \sigma^{-q}\mu_{ij} & \sigma^{-q}
\nu_{ij}\end{pmatrix}} si \m{j=i_0}. Le r\'esultat en d\'ecoule imm\'ediatement.
\end{proof}

\sepprop

On en d\'eduit le diagramme commutatif avec lignes exactes
\xmat{0\ar[r] & T_C\ar[r]\fleq[d] & E(g'_2)\ar[r]\ar[d]^{\tau_q} &
L^*(-qP)\ar[r]\flinc[d] & 0\\
0\ar[r] & T_C\ar[r] & E(g_2)\ar[r] & L^*\ar[r] & 0}
Rappelons qu'on a une application surjective canonique
\[\lambda_g : H^1(E(g_2)\ot L^{n-1})=H^1(E(g_2)^*\ot T_C\ot L^{n-2})\lra
H^1(\rho_n)^{-1}(g)\]
(cf. \ref{fib_gn}). Soit
\[\beta_{n,q}^P:H^1(E(g_2)^*\ot T_C\ot L^{n-2})\lra
H^1(E(g'_2)^*\ot T_C\ot L(qP)^{n-2})\]
le morphisme surjectif d\'eduit de \m{\tau_q} et de \m{L\subset L(qP)}.

\sepprop

\begin{subsub}\label{ecl_3}{\bf Proposition : } Le diagramme suivant est
commutatif :
\xmat{H^1(E(g_2)\ot L^{n-1})\ar[r]^-{\lambda_g}\ar[d]^{\beta_{n,q}^P} &
H^1(\rho_n)^{-1}(g)\ar[d]^{b_{n,q}^P}\\
H^1(E(g'_2)\ot L(qP)^{n-1})\ar[r]^-{\lambda_{g'}} & H^1(\rho_n)^{-1}(g')}
\end{subsub}

\begin{proof}
D\'emonstration analogue \`a la pr\'ec\'edente, utilisant celle de la
proposition \ref{g_aut2_p1}.
\end{proof}

\end{sub}

\sepsub

\Ssect{Courbes multiples scind\'ees}{c_scind}

Soient \m{n\geq 2} un entier et \m{C_n} une courbe multiple primitive de
multiplicit\'e $n$ et de courbe r\'eduite associ\'ee $C$, et de fibr\'e en
droites associ\'e $L$. 

On dit que \m{C_n} est {\em scind\'ee} s'il existe un morphisme \m{C_n\to C}
induisant l'identit\'e sur $C$. Les exemples les plus simples sont les courbes
triviales. Ce sont d'ailleurs les seuls exemples si \m{n=2}, d'apr\`es
\cite{ba_ei} ou \ref{pr0_d}.

On s'int\'eresse aux paires \m{(C_n,\pi)}, o\`u \m{C_n} est une courbe
scind\'ee de multiplicit\'e $n$ et \m{\pi:C_n\to C} est un morphisme induisant
l'identit\'e sur $C$. On dit que deux telles paires, \m{(C_n,\pi)} et
\m{(C'_n,\pi')} sont {\em isomorphes} s'il existe un isomorphisme
\m{\epsilon:C_n\to C'_n} induisant l'identit\'e sur $C$ et tel que
\m{\pi=\pi'\circ\epsilon}. On peut classifier ces paires de la m\^eme mani\`ere
que les courbes primitives, en utilisant des faisceaux de groupes diff\'erents.
On consid\`ere les automorphismes de $\C$-alg\`ebres de \m{\C\DB{x,t}/(t^n)}
laissant invariant \m{\C\DB{x}}. Ils sont du type suivant : soit
\m{\nu\in\C\DB{x,t}/(t^{n-1})} inversible. Alors il existe un unique
automorphisme \m{\psi_\nu} de \m{\C\DB{x,t}/(t^n)} tel que
\m{\psi_\nu(\alpha)=\alpha} pour tout \m{\alpha\in\C\DB{x}} et que
\m{\psi_\nu(t)=\nu t}. Les automorphismes de \m{\C\DB{x,t}/(t^n)} laissant
invariant \m{\C\DB{x}} sont exactement les \m{\psi_\nu}. Le groupe de ces
automorphismes s'identifie donc \`a \m{\big(\C\DB{x,t}/(t^{n-1})\big)^*}, muni
de la loi de groupe
\[\nu'*\nu \ = \ \nu'.\nu(x,\nu't) .\]
On d\'efinit de m\^eme les faisceaux de groupes \m{\ks_n} sur $C$ par : pour
tout ouvert propre $U$ de $C$, \m{\ks_n(U)} est le groupe des automorphismes de
\m{U\times{\bf Z}_n} la projection sur $U$ invariante. On a une suite exacte
canonique de faisceaux de groupes
\xmat{0\ar[r] & \ko_C\ar[r] & \ks_n\ar[r]^-{\rho'_n} & \ks_{n-1}\ar[r] & 0 .}
L'ensemble \m{H^1(C,\ks_n)} s'identifie \`a celui des classes d'isomorphismes
des paires \m{(C_n,\pi)}. On a un morphisme de faisceaux de groupes naturel
\[s_n:\ks_n\lra\kg_n ,\]
Mais \m{\ks_n} n'est pas distingu\'e. Ces inclusions commutent avec les
restrictions \m{\ks_n\to\ks_{n-1}} et \m{\kg_n\to\kg_{n-1}}.

Soient \m{g\in H^1(C,\ks_n)} et \m{\ov{g}} son image dans \m{H^1(C,\kg_n)}.
Soit \m{(C_n,\pi)} la courbe scind\'ee associ\'ee \`a $g$. On a donc
\m{C_n=C(\ov{g})}. On note \m{C_i} la courbe de multiplicit\'e $i$
sous-jacente, pour \m{2\leq i<n} (ce sont bien entendu des courbes scind\'ees,
munies de la restriction de $\pi$), $L$ le fibr\'e en droites sur $C$
associ\'e, et \m{\ov{g}_2} l'image de \m{\ov{g}} dans \m{ H^1(C,\kg_2)}. Alors
on a
\[E(\ov{g}_2) \ \simeq \ L^*\oplus T_C \]
(une courbe double scind\'ee est triviale). Notons que la projection
\m{\pi:C_n\to C} induit une d\'ecomposition canonique en somme directe : le
fibr\'e \m{E(\ov{g}_2)} s'identifie au {\em fibr\'e tangent restreint} \m{\T_1}
(cf. \ref{fa_der} et la proposition \ref{der_rat3}), et le morphisme tangent
associ\'e \`a $\pi$ induit un morphisme \m{E(\ov{g}_2)\to T_C}. L'inclusion
\m{T_C\to E(\ov{g}_2)} est d\'efinie elle canoniquement pour toute courbe
multiple. On montre comme dans les cas des courbes multiples qu'on a un
isomorphisme de faisceaux de groupes \m{(\ko_C)^g\simeq L^{n-2}}. Comme dans le
cas des courbes multiples on obtient si \m{n\geq 3} une application surjective
\[l_g:H^1(L^{n-2})\lra H^1(\rho'_n)^{-1}(\rho'_n(g)) .\]

La relation avec les faisceaux de groupes \m{\kg_n} est donn\'ee par le
diagramme commutatif avec lignes exactes
\xmat{0\ar[r] & \ko_C\ar[r]\flinc[d] & \ks_n\ar[r]\flinc[d] &
\ks_{n-1}\ar[r]\flinc[d] & 0\\
0\ar[r] & \ko_C\oplus T_C\ar[r] & \kg_n\ar[r] & \kg_{n-1}\ar[r] & 0}

Le faisceau de groupes \m{(\ks_n)^g} s'identifie au faisceau \m{\AAut^S_C(C_n)}
des automorphismes de \m{C_n} laissant la projection sur $C$ invariante. Le
diagramme commutatif pr\'ec\'edent induit le suivant
\xmat{0\ar[r] & L^{n-2}\ar[r]\flinc[d] & \AAut^S_C(C_n)\ar[r]\flinc[d] &
\AAut^S_C(C_{n-1})\ar[r]\flinc[d] & 0\\
0\ar[r] & L^{n-2}\oplus(T_C\ot L^{n-1})\ar[r] & \AAut_C(C_n)\ar[r] &
\AAut_C(C_{n-1})\ar[r] & 0}
On a un diagramme commutatif
\xmat{
H^1(L^{n-2})\ar[r]^-{l_g}\flinc[d] &
H^1(\rho'_n)^{-1}(\rho'_n(g))\ar[d]^{H^1(s_n)} \\
H^1(L^{n-2})\oplus H^1(T_C\ot L^{n-1})\ar[r]^-{\lambda_g} &
H^1(\rho_n)^{-1}(\rho_n(\ov{g}))
}
Je ne sais pas si \m{H^1(s_n)} est injective. Ce probl\`eme est li\'e \`a celui
de l'existence de paires \m{(C_n,\pi)} non isomorphes mais donc les courbes
multiples associ\'ees le sont.

On donne au chapitre \ref{mult_3} la classification des courbes multiples
scind\'ees de multiplicit\'e 3.
\end{sub}

\sepsub

\Ssect{Classification des prolongements de courbes multiples dans les cas
simples}{cl_simp}

Soient \m{n\geq 3} un entier et \m{C_{n-1}} une courbe primitive de
multiplicit\'e \m{n-1} de courbe r\'eduite associ\'ee $C$ et de fibr\'e en
droites associ\'e $L$. Soient \m{\gamma\in H^1(C,\kg_{n-1})} tel que
\m{C_{n-1}=C(\gamma)} et \m{\gamma_2\in H^1(C,\kg_2)} l'image de $\gamma$. On
s'int\'eresse aux prolongements de \m{C_{n-1}} en courbes de multiplicit\'e
$n$. Il existe deux cas simples :
\begin{itemize}
\item[--] Le cas o\`u \ \m{h^1(E(\gamma_2)\ot L^{n-1})=0} : il y a dans ce cas
un prolongement unique de \m{C_{n-1}} en courbe de multiplicit\'e $n$.
\item[--] Le cas o\`u \m{\Aut_C(C_{n-1})} est trivial. Les prolongements sont
dans ce cas param\'etr\'es par \m{H^1(C,E(\gamma_2)\ot L^{n-1})}.
\end{itemize}
Soit $g$ le genre de $C$. Le premier cas se produit par exemple quand
\m{E(\gamma_2)} est stable et
\[\deg(L) \ \geq \ \frac{3g-3}{n-\frac{3}{2}} .\]
Des cas o\`u le second cas se produit sont donn\'es dans le corollaire
\ref{fa2_3_7}.

\end{sub}

\newpage

\section{Faisceaux des d\'erivations et fibr\'es tangents restreints}
\label{fa_der}

\Ssect{D\'efinitions et propri\'et\'es \'el\'ementaires}{def_el}

Soit \m{n\geq 2} un entier. Soit \m{C_n} une courbe multiple primitive de
multiplicit\'e $n$ et de courbe r\'eduite associ\'ee $C$. On a une
filtration canonique \m{C=C_1\subset\cdots\subset C_{n-1}\subset C_n}, o\`u
\m{C_i} est multiple primitive de multiplicit\'e $i$. Comme d'habitude on
note \m{\ko_i} le faisceau structural de \m{C_i}, et $L$ le fibr\'e en droites
\m{\ki_C/\ki_{C_2}} sur $C$.

On peut supposer que \m{C_n=C(g)}, avec \m{g\in H^1(C,\kg_n)}. Si \m{2\leq
i\leq n} on notera \m{g_i} l'image de \m{g_n} dans \m{H^1(C,\kg_i)}. On a donc
\m{C(g_i)=C_i}.

\sepsubsub

\begin{subsub}\label{f_rat}Fonctions rationnelles -- \rm On peut voir $C$ comme
un point du sch\'ema \m{C_n}. On appelle {\em fonction rationnelle} sur \m{C_n}
un \'el\'ement de l'anneau \m{\ko_{C_n,C}}. Soient \m{P\in C} et
\m{z\in\ko_{nP}} une \'equation locale de $C$. Il existe d'apr\`es le
th\'eor\`eme \ref{pr0} un isomorphisme
\[\theta:\ko_{nP}\simeq\ko_{CP}\DB{t}/(t^n)\]
compatible avec la projection \m{\ko_{nP}\to\ko_{CP}}. L'anneau des fonctions
rationnelles sur \m{C_n} s'identifie alors au localis\'e
\m{\big(\ko_{CP}\DB{t}/(t^n)\big)_{(t)}}. On peut donc repr\'esenter les
fonctions rationnelles sur \m{C_n} par des sommes
\m{\sigg_{i=0}^{n-1}\alpha_it^i}, o\`u pour \m{1\leq i\leq n}, \m{\alpha_i} est
une fonction rationnelle sur $C$. Mais cette repr\'esentation d\'epend bien
s\^ur de l'isomorphisme $\theta$. On notera \m{\C(C_n)} l'anneau des fonctions
rationnelles sur \m{C_n}.
\end{subsub}

\sepsubsub

\begin{subsub}\label{der_rat}D\'erivations -- \rm Pour tout ouvert $U$ de
\m{C_n} on note \m{\kd_n^0(U)} le \m{\ko_n(U)}-module des d\'erivations de
\m{\ko_n(U)} dans lui-m\^eme. Si $U$ est non vide et distinct de \m{C_n}, tout
\'el\'ement de \m{\kd_n^0(U)} induit de mani\`ere \'evidente une d\'erivation
$D$ de \m{\C(C_n)} telle que pour tout \m{P\in U} on ait
\m{D(\ko_{nP})\subset\ko_{nP}}. Il en d\'ecoule que pour tout ouvert
\m{V\subset U} on a une restriction bien d\'efinie \m{\kd_n^0(U)\to\kd_n^0(V)}.
On obtient ainsi un pr\'efaisceau de \m{\ko_n}-modules \m{\kd^0_n}. Le faisceau
associ\'e \m{\kd_n} est coh\'erent et s'appelle le {\em faisceau des
d\'erivations} sur \m{C_n}. Notons tout ouvert $U$ de \m{C_n} distinct de
\m{C_n} est affine et qu'on a donc \m{\kd_n(U)=\kd_n^0(U)}. On note de m\^eme
\m{\kd_i} le faisceau des d\'erivations sur \m{C_i} si \m{2\leq i<n}.

Le faisceau \m{\kd_n} est le dual du faisceau des diff\'erentielles
\m{\Omega_{C_n}}.
\end{subsub}

\sepprop

\begin{subsub}\label{der_rat1}{\bf Proposition : }
Le faisceau \m{\kd_n} est quasi localement libre de type \m{(0,\ldots,0,1,1)}.
\end{subsub}

Autrement dit, \m{\kd_n} est localement isomorphe \`a \m{\ko_n\oplus\ko_{n-1}}.

\begin{proof}
Cela du fait que, comme indiqu\'e dans \ref{autr_0} les d\'erivations de
\m{\C\DB{x,t}} sont de la forme \m{D=a\frac{\partial}{\partial
x}+tb\frac{\partial}{\partial t}}, avec \m{a,b\in\C\DB{x,t}}.
\end{proof}

\sepprop

\begin{subsub}\label{der_cocyc} Repr\'esentation de \m{\kd_n} par des
1-cocycles -- \rm Soit \m{(g_{ij})} un 1-cocycle d'un recouvrement ouvert
\m{(U_i)} de $C$ repr\'esentant $g$. Donc \m{g_{ij}} est un automorphisme de
\m{\ko_C(U_{ij})[t]/(t^n)} tel que \m{\rho g_{ij}=\rho}, $\rho$ d\'esignant le
morphisme canonique \m{\ko_C(U_{ij})[t]/(t^n)\to\ko_C(U_{ij})}. On peut
supposer que pour tout indice $i$, \m{\omega_{C\mid U_i}} est trivial. On
notera abusivement \m{\frac{\partial}{\partial x}} une section de
\m{\omega_{C\mid U_i}} engendrant ce fibr\'e. Alors le faisceau \m{\kd_{n,i}} 
des d\'erivations de \m{\ko_C(U_{ij})[t]/(t^n)} est isomorphe \`a
\m{\ko_{n\mid U_i}\oplus\ko_{n-1\mid U_i}}, engendr\'e par
\m{\frac{\partial}{\partial x}} et \m{t\frac{\partial}{\partial t}}.
Le faisceau \m{\kd_n} s'obtient en recollant les \m{\kd_{n,i}} au moyen des
automorphismes de \m{\kd_{n,i\mid U_{ij}}}
\begin{equation}\label{def_delt}
\xymatrix{\Delta_{ij} : D\fmaps[r] & g_{ij}\circ D\circ g_{ij}^{-1} \ ,}
\end{equation}
par le proc\'ed\'e indiqu\'e dans \ref{cpa_def3}. Supposons que
\m{g_{ij}=\phi_{\mu\nu}} , avec \m{\mu,\nu\in\ko_{n-1}(U_{ij})}. Soient
\m{\mu',\nu'} les \'el\'ements de \m{\ko_{n-1}(U_{ij})} tels que
\m{g_{ij}^{-1}=\phi_{\mu'\nu'}} . Alors un calcul simple montre qu'on a
\begin{equation}\label{def_delt2}
\Delta_{ij}(\frac{\partial}{\partial x}) \ = \
\big(1+g_{ij}(\DIV{\mu'}{x})\nu t\big).\frac{\partial}{\partial x} +
\nu g_{ij}(\DIV{\nu'}{x}).t\frac{\partial}{\partial t},
\end{equation}
\begin{equation}\label{def_delt3}
\Delta_{ij}(t\frac{\partial}{\partial t}) \ = \
\big(g_{ij}(\DIV{\mu'}{t})\nu^2t-\mu\big).t\frac{\partial}{\partial x} +
\big(g_{ij}(\DIV{\nu'}{t})\nu^2t+1\big).t\frac{\partial}{\partial t} \quad .
\end{equation}
On consid\`ere maintenant l'isomorphisme de recollement
\[\Theta_{ij}:\kd_{n,j\mid U_{ij}}\lra g_{ij}^*(\kd_{n,i\mid U_{ij}}) ,\]
(cf. \ref{cpa_def3}). Compte tenu de la structure de module de
\m{g_{ij}^*(\kd_{n,i\mid U_{ij}})}, la matrice de \m{(\Theta_{ij})}
relativement \`a la base \m{(\frac{\partial}{\partial x},
t\frac{\partial}{\partial t})} est obtenue en appliquant \m{g_{ij}^{-1}} aux
coefficients de \m{\frac{\partial}{\partial x}} et
\m{t\frac{\partial}{\partial t}} dans les formules pr\'ec\'edentes.
\end{subsub}

\sepsubsub

\begin{subsub} Restriction \`a $C$ -- \rm On note \m{E_C} le fibr\'e vectoriel
de rang 2 sur $C$ d\'efini par l'unique extension non triviale \
\m{0\to\ko_C\to E_C\to T_C\to 0} .
\end{subsub}

\sepprop

\begin{subsub}\label{der_rat2}{\bf Proposition : }
On a \ \m{\kd_{n\mid C}\simeq E_C} \ si \m{\deg(L)\not=0} et \ \m{\kd_{n\mid
C}\simeq \ko_C\oplus T_C} \ si \m{\deg(L)=0} .
\end{subsub}

\begin{proof}
On utilise la description pr\'ec\'edente de \m{\kd_n} en termes de 1-cocycles.
On a \m{g_{ij}(t)=\tau_{ij}t}, avec \m{\tau_{ij}} inversible, donc
\m{\nu_{ij}=\rho(\tau_{ij})\in\ko_C^*(U_{ij})}, et d'apr\`es la proposition
\ref{pr0_b}, le 1-cocycle \m{(\nu_{ij})} repr\'esente le fibr\'e en droites $L$
sur $C$. Un calcul simple (utilisant \ref{K_neq2}, (\ref{def_delt2}) et
(\ref{def_delt3})) montre que relativement \`a la base
\m{(\frac{\partial}{\partial x},t\frac{\partial}{\partial t})} de
\m{\kd_{n,i\mid U_{ij}}\ot\ko_{C\mid U_{ij}}}, l'automorphisme induit par
\m{\Delta_{ij}} est d\'efini par la matrice
\[\begin{pmatrix}
1 & 0 \\ -\nu_{ij}\frac{d}{dx}(\frac{1}{\nu_{ij}}) & 1
  \end{pmatrix} \ . \]
Le 1-cocycle \m{(\nu_{ij}\frac{d}{dx}(\frac{1}{\nu_{ij}})dx)} repr\'esente un
\'el\'ement $\lambda$ de \m{H^1(C,\omega_C)\simeq\C}. D'apr\`es \cite{ha},
Chapter III, ex. 7.4 et \cite{mat} on a \m{\lambda=\deg(L)}. Le r\'esultat
d\'ecoule alors ais\'ement de \ref{const_ext}.
\end{proof}

\sepprop

La {\em premi\`ere filtration canonique} de \m{\kd_n} (cf. \cite{dr2}) est 
\[T_C\ot \ki_C^{n-1} \ \subset\kd_2\ot \ki_C^{n-2}\subset\cdots\subset \
\kd_{n-1}\ot\ki_C\subset \ \kd_n \ \]
(c'est-\`a-dire que \m{\ki_C^i\kd_n\simeq\kd_{n-i}\ot\ki_C^i}). Les gradu\'es
se calculent ais\'ement avec la proposition \ref{der_rat2}.

\sepsubsub

\begin{subsub}\label{tan_rest} Fibr\'e tangent restreint -- \rm
Soit $E$ le sous-faisceau de \m{\kd_n} annulateur de \m{\ki_C^{n-1}}. C'est un
fibr\'e vectoriel de rang 2 sur \m{C_{n-1}}. On pose
\[\T_{n-1} \ = \ E\ot(\ki_C)^* ,\]
o\`u le dual de \m{\ki_C} est pris sur \m{C_{n-1}} (le dual de \m{\ki_C} sur
\m{C_n} \'etant isomorphe lui \`a \m{\ko_{n-1}}). On appelle \m{\T_{n-1}} le
{\em fibr\'e tangent restreint de \m{C_n}}. Bien que ce soit un fibr\'e sur
\m{C_{n-1}} il d\'epend effectivement de \m{C_n}, et plus pr\'ecis\'ement de
l'inclusion \m{C_{n-1}\subset C_n}. La classe d'isomorphisme de \m{C_n} seule
d\'etermine \m{\T_{n-1}} \`a l'action de \m{\Aut_C(C_{n-1})} pr\`es.

Consid\'erons les 1-cocycles \m{(g_{ij})} repr\'esentant $g$ et
\m{(\Delta_{ij})} repr\'esentant \m{\kd_n} (cf. \ref{der_cocyc}). Alors le
\m{\ko_{n-1\mid U_{ij}}}-module libre \m{t\ko_{n\mid U_{ij}}\oplus
\ko_{n-1\mid U_{ij}}} est invariant par \m{\Delta_{ij}}, et \m{\T_{n-1}} est
obtenu en recollant ces fibr\'es au moyen des \m{\Delta_{ij}}. Sur l'ouvert
\m{U_i}, \m{\T_{n-1}} est engendr\'e par \m{t\frac{\partial}{\partial x}} et
\m{t\frac{\partial}{\partial t}}.
\end{subsub}

\sepprop

\begin{subsub}\label{der_rat3}{\bf Proposition : }
On a \ \m{\T_i=\T_{n-1\mid C_i}} \ pour \m{2\leq i<n-1}, et \
\m{\T_1\simeq E(g_2)} .
\end{subsub}

\begin{proof}
La premi\`ere assertion est imm\'ediate. La seconde se d\'emontre en utilisant
les formules (\ref{def_delt2}) et (\ref{def_delt3}).
Soit \m{(\psi_{ij})=(\phi_{\mu_{ij}\nu_{ij}})} le 1-cocycle \`a valeurs dans
\m{\kg_2} image de \m{(g_{ij})}. Ce 1-cocycle repr\'esente \m{g_2}. Le
fibr\'e \m{\T_1} sur $C$ est obtenu en recollant les
\[\ko_{C\mid U_i}.t\frac{\partial}{\partial x}\oplus\ko_{C\mid
U_i}.t\frac{\partial}{\partial t}\]
au moyen des matrices \m{\begin{pmatrix}\nu_{ij} & -\mu_{ij}\\
0 & 1\end{pmatrix}} . Donc par d\'efinition de \m{E(g_2)} et d'apr\`es
\ref{desc_ext} on a bien \m{\T_1\simeq E(g_2)} (cf. \ref{g_aut2_l0}) .
\end{proof}

\sepprop

La {\em seconde filtration canonique} de \m{\kd_n} (cf. \cite{dr2}) est
\[\T_1\ot\ki_C^{n-1} \ \subset\cdots\subset \ \T_{n-2}\ot\ki_C^2 \ \subset
\ \T_{n-1}\ot\ki_C \ \subset \ \kd_n\]
(c'est-\`a-dire que l'annulateur de \m{\ki_C^i} dans \m{\kd_n} est \m{\T_i\ot
\ki^{n-i}}).

On a une filtration mixte
\[\T_C\ot \ki_C^{n-1}\subset\T_1\ot\ki_C^{n-1}\subset\cdots\subset
\kd_i\ot\ki_C^{n-i}\subset\T_i\ot\ki_C^{n-i}\subset\cdots\subset
\T_{n-1}\ot\ki_C\subset\kd_n \ .\]
dont les gradu\'es sont
\[(\T_i\ot\ki_C^{n-i})/(\kd_i\ot\ki_C^{n-i})\simeq L^{n-i-1} , \ \ \ \
(\kd_i\ot\ki_C^{n-i})/(\T_{i-1}\ot\ki_C^{n-i+1})\simeq T_C\ot L^{n-i} .\]

\sepsubsub

\begin{subsub} Cas d'une courbe plong\'ee dans une surface lisse -- \rm
On suppose que \m{C_n} est plong\'ee dans une surface lisse $S$. Alors on a
\[\T_{n-1} \ = \ T_{S\mid C_{n-1}} .\]
\end{subsub}

\end{sub}

\sepsub

\Ssect{Prolongements du faisceau des d\'erivations en fibr\'e vectoriel de rang
2 et extensions de la courbe multiple}{prol_dn}

On utilise les notations de \ref{def_el}. On suppose que \m{n\geq 3}.

On montre dans ce qui suit que sous certaines hypoth\`eses, prolonger
\m{C_{n-1}} en une courbe de multiplicit\'e $n$ revient essentiellement \`a
\'etendre \m{\kd_{n-1}} en un fibr\'e vectoriel de rang 2 sur \m{C_{n-1}}
(correspondant au fibr\'e \m{\T_{n-1}} sur la courbe de multiplicit\'e $n$
extension de \m{C_{n-1}}). Dans tout ce qui suit on suppose que le groupe
\m{\Aut_C(C_{n-1})} des automorphismes de \m{C_{n-1}} laissant $C$ invariante
(cf. \ref{fa2}) est trivial.

Le fibr\'e vectoriel \m{\T_{n-1}} de rang 2 sur \m{C_{n-1}} n'est d\'efini
qu'une fois choisie l'extension \m{C(g)=C_n} de \m{C_{n-1}}. C'est un
prolongement du faisceau des d\'erivations de \m{C_{n-1}}. On est donc dans la
situation de \ref{pl_ll} : on a un diagramme commutatif avec lignes et colonnes
exactes
\xmat{ & & 0\ar[d] & 0\ar[d]\\
0\ar[r] & \T_{n-2}\ot\ki_C\ar[r]\fleq[d] & \kd_{n-1}\ar[r]\ar[d] &
T_C\ar[r]\ar[d] & 0\\
0\ar[r] & \T_{n-2}\ot\ki_C\ar[r] & \T_{n-1}\ar[r]\ar[d] &
\T_1=E(g_2)\ar[r]\ar[d] & 0\\
& & L^*\fleq[r]\ar[d] & L^*\ar[d] \\ & & 0 & 0}
Les prolongements possibles de \m{\kd_{n-1}} sont d\'ecrits par le diagramme
commutatif avec lignes et colonnes exactes
\xmat{ & \Hom(T_C,\T_1\ot L^{n-2})\ar[d] \\
 & H^1(\T_1\ot L^{n-1})\ar[d]^\theta \\
0\ar[r] & \Ext^1_{\ko_C}(\T_1,\T_1\ot L^{n-2})\ar[r]^-i\ar[d] &
\Ext^1_{\ko_{n-1}}(\T_1,\T_{n-2}\ot\ki_C)\ar[r]\ar[d] &
\End(\T_1)\ar[r]\ar[d] & 0 \\
0 \ar[r] & \Ext^1_{\ko_C}(T_C,\T_1\ot L^{n-2})\ar[r]\ar[d] &
\Ext^1_{\ko_{n-1}}(T_C,\T_{n-2}\ot\ki_C)\ar[r] & \Hom(T_C,\T_1)\ar[r] & 0 \\
& 0
}
Soit \m{\sigma\in\Ext^1_{\ko_{n-1}}(T_C,\T_{n-2}\ot\ki_C)} correspondant \`a
\m{\T_{n-1}}. Alors d'apr\`es \ref{pl_cl2b} les \'el\'ements de
\m{\Ext^1_{\ko_{n-1}}(T_C,\T_{n-2}\ot\ki_C)} correspondant aux prolongements de
\m{\kd_{n-1}} sont de la forme \m{\sigma+i\circ\theta(u)}, avec \m{u\in
H^1(\T_1\ot L^{n-1})}. On note \m{\T(u)} le fibr\'e vectoriel prolongement de
\m{\kd_{n-1}} associ\'e \`a \m{\sigma+i\circ\theta(u)}.

D'apr\`es \ref{fib_gn} les courbes de multiplicit\'e $n$ extensions de
\m{C_{n-1}} sont associ\'ees aux \'el\'ements de \m{H^1(C,\kg_{n-1})} de la
forme \m{\lambda_g(u)}, avec \m{u\in H^1(C,\T_1\ot L^{n-1})}. On note
\m{\T_{n-1}(u)} le fibr\'e \m{\T_{n-1}} correspondant \`a la courbe
\m{C(\lambda_g(u))}.

\sepprop

\begin{subsub}\label{prol_dn1}{\bf Th\'eor\`eme : } On a \ \m{\T(-(n-1)u)\simeq
\T_{n-1}(u)} \ pour tout \m{u\in H^1(C,\T_1\ot L^{n-1})}.
\end{subsub}

\begin{proof}
On pose \m{g_{ij}=\phi_{\mu_{ij},\nu_{ij}}}, avec \m{\mu_{ij},
\nu_{ij}\in\ko_C(U_{ij})/(t^{n-1})}. Soit \m{\ov{g}=\lambda_g(u)}. Alors
$\ov{g}$ est repr\'esent\'e par un cocycle de la forme \m{(\ov{g}_{ij})}, avec
\[\ov{g}_{ij} \ = \ \lambda_{\theta_{ij},\beta_{ij}}\circ g_{ij} ,\]
avec \m{\theta_{ij},\beta_{ij}\in\ko_C(U_{ij})}. La famille \m{(\theta_{ij},
\beta_{ij})} repr\'esente $u$ au sens de \ref{desc_ext} et \ref{g_aut2_p1}.
Mais au sens de \ref{pl_cl2c}, $u$ est repr\'esent\'e par \m{(
\frac{\theta_{ij}}{(\nu_{ij})_0^{n-1}},\frac{\beta_{ij}}{(\nu_{ij})_0^{n-1}})}.

Un calcul simple utilisant \ref{prox0} montre que
\m{\ov{g}_{ij}=\phi_{\ov{\mu}_{ij},\ov{\nu}_{ij}}}, avec
\[\ov{\mu}_{ij}=\mu_{ij}+\delta_{ij}t^{n-2} , \ \ \ \
\ov{\nu}_{ij}=\nu_{ij}+\epsilon_{ij}t^{n-2},\]
o\`u
\[\delta_{ij}=(\mu_{ij})_0\beta_{ij}+\theta_{ij}, \ \ \ \
\epsilon_{ij}=(\nu_{ij})_0\beta_{ij} .\]
On pose
\[g_{ij}^{-1}=\phi_{\mu'_{ij},\nu'_{ij}}, \ \ \ \
\ov{g}_{ij}^{-1}=\phi_{\ov{\mu}'_{ij},\ov{\nu}'_{ij}} .\]
Alors, en utilisant (\ref{aut_xt5}) on voit que
\[\ov{\mu}'_{ij}=\mu'_{ij}+\delta'_{ij}t^{n-2}, \ \ \ \
\ov{\nu}'_{ij}=\nu'_{ij}+\epsilon'_{ij}t^{n-2},\]
avec
\[\delta'_{ij}=\frac{1}{(\nu_{ij})_0^n}((\mu_{ij})_0\epsilon_{ij}-(\nu_{ij})_0
\delta_{ij})=-\frac{1}{(\nu_{ij})_0^{n-1}}\theta_{ij} ,\]
\[\epsilon'_{ij}=-\frac{1}{(\nu_{ij})_0^n}\epsilon_{ij}=
-\frac{1}{(\nu_{ij})_0^{n-1}}\beta_{ij} .\]
Compte tenu de \ref{cpa_def3}, les formules (\ref{def_delt2}), (\ref{def_delt3})
montrent que \m{\T_{n-1}\ot\ki_C} est d\'efini par les matrices
\[\begin{pmatrix}\frac{1}{\nu'_{ij}}(1+\DIV{\mu'_{ij}}{x}t) &
\frac{1}{\nu'_{ij}}(\mu'_{ij}+\DIV{\mu'_{ij}}{t}t) \\
\frac{1}{\nu'_{ij}}\DIV{\nu'_{ij}}{x}t &
1+\frac{1}{\nu'_{ij}}\DIV{\nu'_{ij}}{t}t
\end{pmatrix} \ \ .
\]
Il en d\'ecoule que \m{\T_{n-1}=\T_{n-1}(0)} est d\'efini par les matrices
\[\begin{pmatrix}A_{ij} & B_{ij} \\ C_{ij} & D_{ij}\end{pmatrix} \ = \
\begin{pmatrix}1+\DIV{\mu'_{ij}}{x}t & \mu'_{ij}+\DIV{\mu'_{ij}}{t}t\\
\DIV{\nu'_{ij}}{x}t & \nu'_{ij}+\DIV{\nu'_{ij}}{t}t\end{pmatrix} \ \ .\]
De m\^eme, \m{\T_{n-1}(u)} est d\'efini par les matrices
\[\begin{pmatrix}\ov{A}_{ij} & \ov{B}_{ij} \\ \ov{C}_{ij} & \ov{D}_{ij}
\end{pmatrix} \ = \ \begin{pmatrix}1+\DIV{\ov{\mu}'_{ij}}{x}t &
\ov{\mu}'_{ij}+\DIV{\ov{\mu}'_{ij}}{t}t\\
\DIV{\ov{\nu}'_{ij}}{x}t & \ov{\nu}'_{ij}+\DIV{\ov{\nu}'_{ij}}{t}t\end{pmatrix}
\ \ .\]
On d\'eduit des calculs pr\'ec\'edents que \m{\ov{A}_{ij}=A_{ij}},
\m{\ov{C}_{ij}=C_{ij}}, et
\[\ov{B}_{ij}=B_{ij}-(n-1)\frac{1}{(\nu_{ij})_0^{n-1}}\theta_{ij}t^{n-2}, \ \ \
\ \ov{D}_{ij}=D_{ij}-(n-1)\frac{1}{(\nu_{ij})_0^{n-1}}\beta_{ij}t^{n-2} ,\]
d'o\`u le r\'esultat d'apr\`es \ref{pl_cl2c}.
\end{proof}

\sepprop

On en d\'eduit le

\sepprop

\begin{subsub}\label{prol_dn2}{\bf Corollaire : } On suppose que
\m{E(g_2)=\T_1} est stable, que \m{\deg(L)\leq 0} et \m{L^{k}\not\simeq\ko_C}
pour \m{1\leq k<2n}. Alors deux prolongements de
\m{C_{n-1}} en courbes de multiplicit\'e $n$ sont isomorphes si et seulement si
les prolongements correspondants \m{\T_{n-1}} de \m{\kd_{n-1}} en faisceau
localement libre de rang 2 sur \m{C_{n-1}} le sont.
\end{subsub}

\begin{proof}
D'apr\`es le corollaire \ref{fa2_3_7} les conditions pr\'ec\'edentes
entrainent la trivialit\'e de \m{\Aut_C(C_{n-1})}. D'autre part
la stabilit\'e de \m{E(g_2)} implique que $C$ est de genre positif et que
\m{h^0(\omega_C\ot E(g_2)\ot L^{n-2})=0}. On a donc \m{\coker(\theta)=
H^1(C,\T_1\ot L^{n-2})} (cf. \ref{pl_cl2b}). Le r\'esultat d\'ecoule donc du
th\'eor\`eme \ref{prol_dn1}.
\end{proof}
\end{sub}

\newpage

\section{Automorphismes des courbes multiples primitives}\label{fa2}

Soit \m{n\geq 2} un entier. Soit \m{C_n} une courbe multiple primitive de
multiplicit\'e $n$ et de courbe r\'eduite associ\'ee $C$ projective. On
supposera que \m{C_n=C(\gamma)}, avec \m{\gamma\in H^1(C,\kg_n)}. Si \m{2\leq
i\leq n} on notera \m{\gamma_i} l'image de $\gamma$ dans \m{H^1(C,\kg_i)}.
On a donc \m{C(\gamma_i)=C_i}. On supposera que $\gamma$ est repr\'esent\'e par
un 1-cocycle \m{(g_{ij})} d'un recouvrement ouvert \m{(U_i)} tel que pour tout
$i$, \m{\omega_{C\mid U_i}} soit trivial.

Soient \m{\AAut_C(C_n)} le faisceau de groupes des automorphismes de \m{C_n}
laissant $C$ invariante, et \m{\Aut_C(C_n)} le groupe de ses sections globales.
On a un isomorphisme canonique
\[\AAut_C(C_n) \simeq \ \kg_n^\gamma .\]
On a donc d'apr\`es la proposition \ref{g_aut2_p1}, si \m{n\geq 3}, une suite
exacte de faisceaux de groupes
\begin{equation}\label{fa2_1_0}0\lra E(\gamma_2)\ot
L^{n-1}\lra\AAut_C(C_n)\lra\AAut_C(C_{n-1})\lra 0 .\end{equation}

\sepsub

\Ssect{Automorphismes des courbes doubles}{aut_db}

\begin{subsub}\label{fa2_1}{\bf Proposition : } On suppose que \m{C_2} est non
triviale. Alors

1 - Les automorphismes de \m{C_2} induisant l'identit\'e sur $C$ sont de la
forme \m{\chi_D=I_{C_2}+D}, o\`u $D$ est une section globale de \m{T_C\ot L}.

2 - On a un isomorphisme canonique de groupes \ \m{\Aut_C(C_2)\simeq
H^0(C,T_C\ot L)} .
\end{subsub}

\begin{proof}
Les automorphismes de \m{C_2} sont de la forme \m{\chi_D=I_{C_2}+D}, o\`u $D$
est une d\'erivation \`a valeurs dans \m{\ki_C}, c'est-\`a-dire un \'el\'ement
de \m{H^0(\T_1\ot\ki_C)}. D'apr\`es la suite exacte (\ref{fib_x2b}) et le fait
que cette suite est non scind\'ee, \m{C_2} \'etant non triviale, on a
\[H^0(\T_1\ot\ki_C) \ \simeq \ H^0(T_C\ot L) .\]
On obtient un isomorphisme de groupes dans 2- \`a cause de l'\'egalit\'e
\m{\chi_D\circ\chi_{D'}=\chi_{D+D'}} .
\end{proof}

\sepprop

Si \m{C_2} est triviale, on a un isomorphisme canonique
\[\Aut(C_2) \ \simeq \ H^0(T_C\ot L)\times\C^* ,\]
la loi de groupe \'etant
\[(\mu,\lambda).(\mu',\lambda') \ = \ (\lambda\mu'+\mu,\lambda\lambda') .\]
Le sous-groupe \m{H^0(T_C\ot L)} est le m\^eme que pour une courbe non
triviale. Pour interpr\`eter le groupe \m{\C^*} on utilise la construction
de \m{C_2} (cf. \ref{pr0_c}). L'automorphisme de \m{C_2} correspondant \`a
\m{\lambda\in\C^*} provient de l'homoth\'etie du fibr\'e \m{L^*} de rapport
$\lambda$.
\end{sub}

\sepsub

\Ssect{Automorphismes des courbes de multiplicit\'e $n$}{aut_n}

\begin{subsub}\label{fa2_2}Le morphisme \m{\Aut_C(C_n)\to\C^*} -- \rm
(Cf. \ref{mor_1}). Comme dans le cas \m{n=2} (cf. \ref{g_aut2_l0}) on a
\m{(\ko_C^*)^\gamma=\ko_C^*}, donc le morphisme canonique \
\m{\xi_n:\kg_n\to\ko_C^*} \ induit un morphisme surjectif
\m{\xi_n^\gamma:(\kg_n)^\gamma\to\ko_C^*}, d'o\`u
\[H^0(\xi_n^\gamma):\Aut_C(C_n)\lra\C^* ,\]
qu'on peut aussi d\'ecrire de la fa\c con suivante : soit
\m{\Phi\in\Aut_C(C_n)}, d\'efini par une famille \m{(\phi_i)}, o\`u \m{\phi_i}
est un automorphisme de \m{\ko_C(U_i)[t]/(t^n)}, telle que \m{\phi_jg_{ij}=
g_{ij}\phi_i} sur \m{U_{ij}} pour tous $i$, $j$, et supposons que
\m{\phi_i=\phi_{\mu_i,\nu_i}}. Alors on a pour tout i
\[H^0(\xi_n^\gamma)(\Phi) \ = \ (\nu_i)_0 .\]
Le rapport avec \ref{chgt}, \ref{K_n} est le suivant : pour tout point
ferm\'e $P$ de \m{C_n} on a 
\[\xi^\gamma_{n,P} \ = \ \xi_{\ko_{n,P}}^n .\]
Si \m{C_n} est triviale, on a une inclusion naturelle
\m{\C^*\subset\Aut_C(C_n)} (cela d\'ecoule de la description des courbes
triviales, cf. \ref{pr0_c}), et on a \'evidemment
\m{H^0(\xi_n^\gamma))(\lambda)=\lambda} pour tout \m{\lambda\in\C^*}.

\end{subsub}

\sepprop

\begin{subsub}\label{fa2_3}{\bf Th\'eor\`eme : } Si \m{\imm(H^0(\xi_n^\gamma))}
contient un \m{\lambda} tel que \m{\lambda^i\not=1} pour \m{1\leq i\leq n-1},
alors \m{C_n} est triviale.
\end{subsub}

\begin{proof} On proc\`ede par r\'ecurrence sur $n$. Le cas \m{n=2} est
r\`egl\'e par la proposition \ref{fa2_1}. Supposons que \m{n\geq 3} et que le
r\'esultat soit vrai pour \m{n-1}. Soient \m{\lambda\in\imm(H^0(\xi_n^\gamma))}
tel que \m{\lambda^i\not=1} pour \m{1\leq i\leq n-1}, et \m{\Phi\in\Aut_C(C_n)}
tel que \m{H^0(\xi_n^\gamma)(\Phi)=\lambda}. En consid\'erant l'automorphisme
de \m{C_{n-1}} induit par $\Phi$ on voit que \m{C_{n-1}} est triviale. On peut
donc mettre les \m{g_{ij}} sous la forme
\[g_{ij} \ = \ \lambda_{\theta_{ij}\beta_{ij}}\circ\phi_{0,\nu_{ij}} \ = \
\phi_{\theta_{ij}t^{n-2},\nu_{ij}(1+\beta_{ij}t^{n-2})} \ ,\]
o\`u \m{(\nu_{ij})} repr\'esente $L$ et \m{(\theta_{ij},\beta_{ij})} un
\'el\'ement $\omega$ de \m{H^1(E(\gamma_2)\ot L^{n-1})} (puisque \m{C_2} est
triviale on a \m{E(\gamma_2)=T_C\oplus L^*}). Il faut montrer que \m{\omega=0}.
On peut supposer que les fonctions \m{\nu_{ij}} sont non constantes. Compte
tenu de \ref{g_aut2_p1} on a \m{\theta_{ij}=-\nu_{ij}^{n-1}\theta_{ji}} et
\m{\beta_{ij}=-\nu_{ij}^{n-2}\beta_{ji}}.

L'automorphisme $\Phi$ est d\'efini par une famille \m{(\phi_i)}, o\`u
\m{\phi_i} est un automorphisme de\Nligne \m{\ko_C(U_i)[t]/(t^n)}, telle que
\begin{equation}\label{fa2_3_1}\phi_ig_{ij} \ = \ g_{ij}\phi_j\end{equation}
pour tous $i$, $j$. Posons \m{\phi=\phi_{\mu_i,\lambda}}. Alors
d'apr\`es la proposition \ref{prox0}, (\ref{fa2_3_1}) se traduit par les deux
relations suivantes dans \m{\ko_C(U_{ij})[t]/(t^{n-1})} :
\begin{equation}\label{fa2_3_2}\mu_i+\lambda^{n-1}\theta_{ij}t^{n-2} \ = \
\nu_{ij}\big(\sigg_{k=0}^{n-2}(\mu_j)_k\nu_{ij}^kt^k\big)+
(\theta_{ij}+\nu_{ij}\beta_{ij}(\mu_i)_0)t^{n-2} ,
\end{equation}
\begin{equation}\label{fa2_3_3}\lambda\phi_{\mu_i,\lambda}(\nu_{ij})+
\lambda^{n-1}\nu_{ij}\beta_{ij}t^{n-2} \ = \ \lambda\nu_{ij}+
\lambda\nu_{ij}\beta_{ij}t^{n-2} .
\end{equation}
De (\ref{fa2_3_3}) on d\'eduit
\begin{equation}\label{fa2_3_4}\phi_{\mu_i,\lambda}(\nu_{ij}) \ = \
\nu_{ij}+\beta_{ij}\nu_{ij}(1-\lambda^{n-2})t^{n-2} .\end{equation}
En utilisant l'\'egalit\'e
\[\phi_{\mu_i,\lambda}(\nu_{ij}) \ = \ \sigg_{i=0}^{n-1}\frac{1}{k!}
(\mu_it)^k\frac{d^k\nu_{ij}}{dx}\]
et le fait que \m{\nu_{ij}} n'est pas constant on voit que \m{\mu_i} est de la
forme
\[\mu_i \ = \ a_it^{n-3}+b_it^{n-2} ,\]
avec \m{a_i,b_i\in\ko_C(U_i)}, et (\ref{fa2_3_4}) devient
\begin{equation}\label{fa2_3_5}\beta_{ij} \ = \ \frac{1}{(1- \lambda^{n-2})
\nu_{ij}}\frac{d\nu_{ij}}{dx}a_i \ .\end{equation}
En \'ecrivant cette relation avec \m{(j,i)} au lieu de \m{(i,j)} et en
utilisant les \'egalit\'es \m{\beta_{ij}=-\nu_{ij}^{n-2}\beta_{ji}} et
\m{\nu_{ji}=1/\nu_{ij}}, on obtient \m{\nu_{ij}^{n-2}=a_i/a_j}, d'o\`u en
d\'erivant
\[(n-2)\nu_{ij}^{n-3}\frac{d\nu_{ij}}{dx} \ = \ \frac{1}{a_j^2}
(a_i\frac{da_j}{dx}-a_j\frac{da_i}{dx}) \ .\]
Cette relation, utilis\'ee dans (\ref{fa2_3_5}) donne
\[\beta_{ij} \ = \ \frac{1}{(n-2)(1- \lambda^{n-2})}
(\nu_{ij}^{n-2}\frac{1}{a_j}\frac{da_j}{dx}-\frac{1}{a_i}\frac{da_i}{dx}) \ ,\]
donc \m{(\beta_{ij})} repr\'esente $0$ dans \m{H^1(L^{n-2})}, et on peut
supposer que \m{\beta_{ij}=0} pour tous $i$, $j$. Il faut bien entendu
v\'erifier qu'on peut obtenir une expression ind\'ependante du choix de la base
\m{\DIV{}{x}} de \m{T_{C\mid U_i\cup U_j}}. Cela se voit ais\'ement en
utilisant la proposition \ref{lemx1}.

La relation (\ref{fa2_3_2}) s'\'ecrit, en utilisant ce qui pr\'ec\`ede,
\[\theta_{ij} \ = \ \frac{1}{1- \lambda^{n-1}}(b_i-\nu_{ij}^{n-1}b_j) \ ,\]
donc \m{(\theta_{ij})} repr\'esente $0$ dans \m{H^1(T_C\ot L^{n-1})}, et on a
bien \m{\omega=0}.
\end{proof}

\sepprop

On va en d\'eduire que \m{\Aut_C(C_n)} est particuli\`erement simple, voire
trivial, dans un certain nombre de cas. Rappelons qu'on a une suite exacte
\begin{equation}\label{fa2_3_6}0\lra T_C\lra E(\gamma_2)\lra L^*\lra
0\end{equation}
(cf. \ref{fib_x2}).

\sepprop

\begin{subsub}\label{fa2_3_7}{\bf Corollaire : } Soient $g$ le genre de $C$ et
\m{d=\deg(L)}. On suppose que les conditions suivantes sont r\'ealis\'ees :
\begin{itemize}
\item[--] Si $g\geq 1$ : $d\leq 0$, et en cas d'\'egalit\'e \
$L^k\not\simeq\ko_C$ pour \m{1\leq k\leq n-1}.
\item[--] Si $g=0$ : $d<-2$.
\end{itemize}
Alors :

1 - Si la suite exacte (\ref{fa2_3_6}) n'est pas scind\'ee, alors
\m{\Aut_C(C_n)} est trivial.

2 - Si la suite exacte (\ref{fa2_3_6}) est scind\'ee et \m{C_n} non triviale,
alors \m{\Aut_C(C_n)} est fini, et si \m{C_n} est triviale, alors on a
\m{\Aut_C(C_n)\simeq\C^*}.
\end{subsub}

\begin{proof}
D\'ecoule de la proposition \ref{fa2_1} et des suites exactes (\ref{fa2_1_0}).
\end{proof}

\sepprop

\begin{subsub}Le noyau de \m{H^0(\xi_n^\gamma)} -- \rm On notera \
\m{\Aut^0_C(C_n)} \ le noyau de \m{H^0(\xi_n^\gamma)}. Soit\Nligne \m{D\in
H^0(C_{n-1},\kd_{n-1}\ot\ki_C)}. Pour tout ouvert $U$ de \m{C_n}, on en d\'eduit
une d\'erivation $D$ de \m{\ko_n(U)} telle que \ \m{\imm(D)\subset\ki_C} \ et
\ \m{D(\ki_{C\mid U}))\subset\ki_{C\mid U}^2}, et un \'el\'ement \m{\chi_D} de
\m{\Aut_C^0(C_n)} d\'efini par : pour tout ouvert $U$ de \m{C_n} et tout
\m{\lambda\in\ko_n(U)},
\[\chi_D(\lambda) \ = \ \sigg_{k\geq 0}\frac{1}{k!}D^k(\lambda) .\]
Cette notation est coh\'erente avec elle de \ref{fa2_1}. Du th\'eor\`eme
\ref{autr_2} on d\'eduit ais\'ement le
\end{subsub}

\sepprop

\begin{subsub}\label{fa2_3_8}{\bf Th\'eor\`eme : } Pour tout
\m{\phi\in\Aut^0_C(C_n)} il existe un unique \m{D\in
H^0(C_{n-1},\kd_{n-1}\ot\ki_C)} tel que \m{\phi=\chi_D}.
\end{subsub}

\sepprop

De la suite exacte (\ref{fa2_1_0}) on d\'eduit la suivante :
\[0\lra H^0(E(g_2)\ot L^{n-1})\lra\Aut_C(C_n)\lra\Aut_C(C_{n-1}) .\]
Le sous-groupe \m{H^0(E(g_2)\ot L^{n-1})} est contenu dans \m{\Aut^0_C(C_n)},
on a donc la suite exacte
\[0\lra H^0(E(g_2)\ot L^{n-1})\lra\Aut_C^0(C_n)\lra\Aut_C^0(C_{n-1}) .\]
Cette suite exacte est la m\^eme que la suivante
\[0\lra H^0(E(g_2)\ot L^{n-1})\lra H^0(C_{n-1},\kd_{n-1}\ot\ki_C)\lra
H^0(C_{n-2},\kd_{n-2}\ot\ki_C)\]
d\'eduite de 
\[0\lra E(g_2)\ot L^{n-1}\lra\kd_{n-1}\ot\ki_C\lra(\kd_{n-1}\ot\ki_C)_{\mid
C_{n-1}}\lra 0 .\]
Mais bien entendu la structure de groupe sur \m{H^0(C_{n-1},\kd_{n-1}\ot\ki_C)}
d\'eduite de celle de \m{\Aut_C(C_{n-1})} n'est pas l'addition.

\end{sub}

\sepsub

\Ssect{L'action de \m{\Aut_C(C_{n-1})}}{act_cn1}

On suppose que \m{n\geq 3}.
On d\'ecrit ici l'action de \m{\Aut_C(C_{n-1})} sur \m{H^1(E(\gamma_2)\ot
L^{n-1})}. Rappelons que si \m{\kc_n(C_{n-1})} d\'esigne l'ensemble des courbes
primitives de multiplicit\'e $n$ qui sont des prolongements de \m{C_{n-1}},
alors on a une application surjective
\[\lambda_\gamma : H^1(E(\gamma_2)\ot L^{n-1})\lra\kc_n(C_{n-1})\]
dont les fibres sont pr\'ecis\'ement les orbites de cette action (cf.
\ref{suite_c}, \ref{pr0_f}).

Soit \m{\psi\in\Aut_C(C_{n-1})}. On suppose que pour tout $i$, \m{\psi_{U_i}}
se rel\`eve en un automorphisme \m{\ov{\psi}_i} de \m{C_{n\mid U_i}}. Soit
\m{\psi_i} l'automorphisme de \m{U_i\times{\bf Z}_n} d\'efini par le carr\'e
commutatif
\xmat{C_{n\mid U_i}\ar[r]^-\simeq\ar[d]^{\ov{\psi}_i} &
U_i\times{\bf Z}_n\ar[d]^{\psi_i} \\
C_{n\mid U_i}\ar[r]^-\simeq & U_i\times{\bf Z}_n}
Il en d\'ecoule que \m{\psi_i\circ g_{ij}\circ\psi_j^{-1}\circ g_{ij}^{-1}} est
trivial sur \m{C_{n-1\mid U_i}} pour tous $i$, $j$, et s'identifie donc \`a une
section de \m{\big(E(\gamma_2)\ot L^{n-1}\big)_{\mid U_i}}. On posera donc
\[\psi_i\circ g_{ij}\circ\psi_j^{-1}\circ g_{ij}^{-1} \ = \
\lambda_{v_{ij}w_{ij}} \ .\]
Soit \ \m{\delta=H^0(\xi_{n-1}^\gamma)(\psi)\in\C^*} . Alors \m{\psi_i} se met
sous la forme \ \m{\psi_i=\phi_{\mu_i\delta}} .

\sepprop

\begin{subsub}\label{act_cn1_1}{\bf Proposition : }
Soit \m{u\in H^1(E(\gamma_2)\ot L^{n-1})}, repr\'esent\'e par un
cocycle \m{(\lambda_{\theta_{ij}\beta_{ij}})} (cf. \ref{ker_rho}, \ref{cas_n}).
Alors \m{\psi.u} est repr\'esent\'e par le
cocycle \m{(\lambda_{\theta'_{ij}\beta'_{ij}})}, avec
\[\theta'_{ij} \ = \ \delta^{n-1}\theta_{ij}-(\mu_i)_0.\delta^{n-2}\beta_{ij}
+ v_{ij} , \ \ \ \
\beta'_{ij} \ = \ \delta^{n-2}\beta_{ij}+w_{ij} \ .\]
\end{subsub}

\begin{proof}
D'apr\`es \ref{suite_c}, \m{\psi.u} est repr\'esent\'e par le
cocycle
\[\big(\psi_i\circ\lambda_{\theta_{ij}\beta_{ij}}\circ g_{ij}\circ\psi_j^{-1}
\circ g_{ij}^{-1}\big)\ = \ \big(\psi_i\circ\lambda_{\theta_{ij}\beta_{ij}}
\circ\psi_i^{-1}\circ\lambda_{v_{ij}w_{ij}}\big)  \ ,\]
et le r\'esultat d\'ecoule ais\'ement du lemme \ref{g_aut2_l1}.
\end{proof}

\sepprop

\begin{subsub}\label{act_cn1_2}Cas o\`u \m{C_{n-1}} est triviale -- \rm 
Soient $g$ le genre de $C$ et \m{d=\deg(L)}. On suppose que \m{C_{n-1}} est
triviale et que les conditions suivantes sont r\'ealis\'ees :
\begin{itemize}
\item[--] Si $g\geq 1$ : $d<0$, ou $d=0$ et $L^k\not=\ko_C$ pour $1\leq k\leq
n-1$.
\item[--] Si $g=0$ : $d<-2$.
\end{itemize}
On a dans ce cas \ \m{E(\gamma_2)=T_C\oplus L^*}, donc
\[H^1(E(\gamma_2)\ot L^{n-1}) \ \simeq \ H^1(T_C\ot L^{n-1})\oplus H^1(L^{n-2})
\ .\]
D'apr\`es le corollaire \ref{fa2_3_7} on a \m{\Aut_C(C_{n-1})=\C^*}, et son
action sur \m{H^1(E(\gamma_2)\ot L^{n-1})} est donn\'ee par :
\[{\xymatrix@R=2pt{
\C^*\times \big(H^1(T_C\ot L^{n-1})\oplus H^1(L^{n-2})\big)\ar[r] &
H^1(T_C\ot L^{n-1})\oplus H^1(L^{n-2}) \\
(\delta,(\theta,\beta))\fmaps[r] & (\delta^{n-1}\theta,\delta^{n-2}\beta)
}}\]
Il en d\'ecoule qu'en g\'en\'eral \m{\kc_n(C_{n-1})} est constitu\'e de la
courbe triviale de multiplicit\'e $n$ extension de \m{C_{n-1}} et des courbes
non triviales, qui s'identifient aux points ferm\'es d'un {\em espace projectif
tordu} (appel\'e aussi {\em espace projectif anisotrope}, cf. \cite{de},
\cite{be_ro}) de dimension \ \m{4g-5-(2n-3)d}. Le lieu singulier de cet
espace projectif tordu est constitu\'e de deux espaces projectifs de dimensions
respectives \ \m{g-2-(n-2)d} \ et \ \m{3g-4-(n-1)d}. Le premier espace
projectif correspond aux prolongements de \m{C_{n-1}} qui sont des courbes {\em
scind\'ees} non triviales (cf. \ref{c_scind}). Le second correspond aux courbes
non triviales \m{C_n} telles que \m{\ki_C} soit isomorphe \`a \m{\pi^*(L)},
$\pi$ d\'esignant la projection \m{C_{n-1}\to C}.
\end{subsub}

\end{sub}

\sepsec

\section{Courbes primitives de multiplicit\'e 3}\label{mult_3}

\Ssect{Les cas simples}{mult_3_1}

Les courbes de multiplicit\'e 2 sont bien connues (cf. \cite{ba_ei} ou
\ref{pr0_d}). Soit $C$ une courbe projective lisse irr\'eductible de genre $g$.
Rappelons que les courbes doubles non scind\'ees de courbe r\'eduite associ\'e
$C$ et telles que le faisceau d'id\'eaux de $C$ soit \m{L\in\Pic(C)} sont
param\'etr\'ees par \m{\P(H^1(T_C\ot L))}. Il existe une seule courbe scind\'ee
de courbe r\'eduite associ\'e $C$ et telle que le faisceau d'id\'eaux de $C$
soit $L$, c'est la courbe triviale d\'ecrite en \ref{pr0_c}.

Soit \m{C_2} une courbe double de courbe r\'eduite associ\'ee $C$ et de fibr\'e
en droites sur $C$ associ\'e $L$ de degr\'e $d$. On s'int\'eresse aux courbes
triples (c'est-\`a-dire de multiplicit\'e 3) qui sont des prolongements de
\m{C_2}. Soit \m{\kc_3(C_2)} l'ensemble de ces courbes. Soit \m{\gamma\in
H^1(C,\kg_2)} l'\'el\'ement correspondant \`a \m{C_2}. D'apr\`es \ref{pr0_f},
il existe une surjection
\begin{equation}\label{m3_0}H^1(E(\gamma)\ot L^2)\lra\kc_3(C_2)\end{equation}
qui d\'epend du choix d'un \'el\'ement de \m{\kc_3(C_2)}, et dont les fibres
sont les orbites de l'action de \m{\Aut_C(C_2)}.

\sepprop

\begin{subsub}\label{m3_1}{\bf Proposition : } S'il existe plusieurs
prolongements possibles de \m{C_2} en courbes de multiplicit\'e 3, alors on a
\m{d<2g-2}, ou \m{d=2g-2} et \m{L\simeq\omega_C}.
\end{subsub}

\begin{proof} S'il existe plusieurs prolongements possibles de \m{C_2} on doit
avoir \Nligne \m{h^1(E(\gamma)\ot L^2)>0} d'apr\`es la description
pr\'ec\'edente de \m{\kc_3(C_2)}. De la suite exacte \Nligne \m{0\to T_C\to
E(\gamma)\to L^*\to 0} \ on d\'eduit la suivante
\[H^1(T_C\ot L^2)\lra H^1(E(\gamma)\ot L^2)\lra H^1(L)\lra 0 .\]
Donc si \m{h^1(E(\gamma)\ot L^2)>0}, on a \m{h^1(T_C\ot L^2)>0} ou
\m{h^1(L)>0}. Le r\'esultat d\'ecoule du fait que si \m{D\in\Pic(C)} est tel
que \m{h^1(D)>0}, alors on a \m{\deg(L)<2g-2} ou \m{\deg(L)=2g-2} et
\m{L=\omega_C}.
\end{proof}

\sepprop

On suppose dans toute la suite que \m{d<2g-2}, ou \m{d=2g-2} et
\m{L\simeq\omega_C}. Les cas les plus simples sont ceux pour lesquels
\m{\Aut_C(C_2)} est trivial. Dans ce cas l'application (\ref{m3_0}) est une
bijection et \m{\kc_3(C_2)} s'identifie \`a \m{H^1(E(\gamma)\ot L^2)} de
mani\`ere non canonique (\`a une translation pr\`es d'apr\`es \ref{cas_com}).
Le groupe \m{\Aut_C(C_2)} a \'et\'e calcul\'e en \ref{fa2}. On en d\'eduit
ais\'ement que les seuls cas qui restent o\`u il est non trivial sont les
suivants :
\begin{enumerate}
\item[--] $d < 2g-2$ et $C_2$ est triviale.
\item[--] $L=\omega_C$.
\end{enumerate}

\end{sub}

\sepsub

\Ssect{Les cas o\`u \m{d<2g-2} et o\`u \m{C_2} est triviale}{mult_3_2}

On a alors \m{\Aut_C(C_2)=\C^*} et \ \m{H^1(E(\gamma)\ot L^2)=H^1(T_C\ot L^2)
\oplus H^1(L)}. Soit \m{C_3^0} la courbe triviale de multiplicit\'e 3
prolongement de \m{C_2}.

Soit \m{(\nu_{ij})} un cocycle repr\'esentant $L$, relativement \`a un
recouvrement \m{(U_i)} de $C$ tel que chaque \m{\omega_{C\mid U_i}} soit
trivial. Alors \m{C^0_3} est repr\'esent\'ee par le 1-cocycle
\m{(\phi_{0,\nu_{ij}})}. Si \m{\tau\in\C^*}, l'automorphisme correspondant
\m{\psi_\tau} de \m{C_2} s'\'etend en un automorphisme de \m{C_3^0}
repr\'esent\'e par la famille \m{(\phi_{0,\tau})} d'automorphismes des
\m{U_i\times{\bf Z}_n}.

Soit \m{\sigma\in H^1(E(\gamma)\ot L^2)}, repr\'esent\'e par une famille
\m{(\lambda_{\theta_{ij}\beta_{ij}})} (cf. \ref{ker_rho}, \ref{g_aut2_p1}).
Alors d'apr\`es la proposition \ref{act_cn1_1}, \m{\psi_\tau.\sigma} est
repr\'esent\'e par la famille
\[(\phi_{0,\tau}\circ\lambda_{\theta_{ij}\beta_{ij}}\circ\phi_{0,\nu_{ij}}
\circ\phi_{0,1/\tau}\circ\phi_{0,1/\nu_{ij}}) \ = \
(\lambda_{\tau^2\theta_{ij},\tau\beta_{ij}}) \ .\]
On en d\'eduit l'action de \m{\Aut_C(C_2)} sur
\m{H^1(E(\gamma)\ot L^2)} :
\[{\xymatrix@R=2pt{
\C^*\times(H^1(T_C\ot L^2)\oplus H^1(L))\ar[r] &  H^1(T_C\ot L^2)\oplus H^1(L)\\
(\tau,(\theta,\beta))\fmaps[r] & (\tau^2\theta,\tau\beta) }}\]
Comme attendu, \m{\C^*} laisse $0$ invariant, ce qui correspond au fait que les
automorphismes de \m{C_2} se prolongent en automorphismes de \m{C_3^0}.

Les prolongements de \m{C_2} sont donc d'une part la courbe triviale \m{C_3^0}
et d'autre part les prolongements non triviaux, param\'etr\'es par le quotient
de \m{(H^1(T_C\ot L^2)\oplus H^1(L))\backslash\lbrace 0\rbrace} par l'action
pr\'ec\'edente de \m{\C^*}. Si \m{h^1(E(\gamma)\ot L^2)\not=0}, ce quotient est
un {\em espace projectif tordu}. Les courbes scind\'ees forment un
sous-espace projectif isomorphe \`a \m{P(H^1(L))}, si \m{h^1(L)\not=0}.
\end{sub}

\sepsub

\Ssect{Le cas o\`u \m{L=\omega_C}}{mult_3_3}

\begin{subsub} Le cas o\`u \m{C_2} est triviale -- \rm
On a alors \ \m{H^1(E(\gamma)\ot L^2)=H^1(\omega_C)\oplus H^1(\omega_C)\simeq
\C^2}~, et \ \m{H^0(T_C\ot L)=\C}.
Comme pr\'ec\'edemment \m{C^0_3} est repr\'esent\'ee par le 1-cocycle
\m{(\phi_{0,\nu_{ij}})}. Soit \m{D_i=\frac{d}{dx_i}} une base de \m{T_{C\mid
U_i}}. Sur \m{U_{ij}} on a \m{D_j=\frac{1}{\nu_{ij}}D_i} (car \m{L=\omega_C}).
Soit \m{\mu\in H^0(T_C\ot L)} et \m{\psi_\mu} l'automorphisme associ\'e de
\m{C_2} (cf. \ref{aut_db}). Soit \m{\psi_{\mu,1}^{D_i}} l'automorphisme de
\m{U_i\times{\bf Z}_2} tel que le diagramme suivant soit commutatif :
\xmat{C_{2\mid U_i}\ar[r]^-\simeq\ar[d]^{\psi_{\mu,1}^{D_i}} &
U_i\times{\bf Z}_2\ar[d]^{\psi_\mu} \\
C_{2\mid U_i}\ar[r]^-\simeq & U_i\times{\bf Z}_2}
(cf. \ref{chgt}). On consid\`ere maintenant les automorphismes aussi not\'es
\m{\psi_{\mu,1}^{D_i}} de \m{U_i\times{\bf Z}_3} extensions des pr\'ec\'edents.
Soit
\[\lambda_{v_{ij}w_{ij}} \ = \ \psi_{\mu,1}^{D_i}\circ\phi_{0,\nu_{ij}}
\circ(\psi_{\mu,1}^{D_j})^{-1}\circ\phi_{0,1/\nu_{ij}} \ . \]

\sepprop

\begin{subsub} {\bf Lemme : } On a
\[v_{ij} \ = \ \frac{D_i(\nu_{ij})}{\nu_{ij}}\mu \ , \ \ \ \
w_{ij} \ = \ -\frac{1}{2}\frac{D_i(\nu_{ij})}{\nu_{ij}}\mu^2 \ .\]
\end{subsub}

\begin{proof}
Cela d\'ecoule du calcul simple
\[\lambda_{v_{ij}w_{ij}}(x_j)=x_j-\frac{1}{2}\frac{D_i(\nu_{ij})}{\nu_{ij}}
\mu^2 t^2 \ , \ \ \ \
\lambda_{v_{ij}w_{ij}}(t)=t+\frac{D_i(\nu_{ij})}{\nu_{ij}}\mu t^2 \ .\]
\end{proof}

\sepprop

On en d\'eduit avec la proposition \ref{act_cn1_1} l'action de \m{\Aut_C(C_2)}
sur \m{H^1(E(\gamma)\ot L^2)} :
\[{\xymatrix@R=2pt{
\C\times\C^2\ar[r] &  \C^2\\
(\delta,(\theta,\beta))\fmaps[r] &
(\theta-\beta\delta-(g-1)\delta^2,\beta+(2g-2)\delta) }}\]
Les facteurs \m{g-1} et \m{2g-2} sont caus\'es par le cocycle 
\m{(\frac{D_i(\nu_{ij})}{\nu_{ij}}dx_i)} (cf. \ref{der_rat2}). Il faut
aussi consid\'erer l'action du sous-groupe \m{\C^*} de \m{\Aut_C(C_2)}, qui a
\'et\'e calcul\'ee en \ref{mult_3_2}. Il en d\'ecoule ais\'ement qu'il n'y a
que deux \m{\Aut_C(C_2)}-orbites : celle de $0$ correspond \`a la courbe
triviale prolongement de \m{C_2} et l'autre est un prolongement non trivial qui
est une courbe scind\'ee.
\end{subsub}

\sepsubsub

\begin{subsub} Le cas o\`u \m{C_2} n'est pas triviale -- \rm
Ici on a \m{H^1(E(\gamma)\ot L^2)\simeq\C}, et l'action de \Nligne
\m{H^0(T_C\ot L)=\C} est donn\'ee par
\[{\xymatrix@R=2pt{
\C\times\C\ar[r] &  \C\\
(\delta,\beta)\fmaps[r] & \beta+(2g-2)\delta
}}\]
On en d\'eduit que si \m{g\not=1} il existe un unique prolongement de
multiplicit\'e 3 de \m{C_2}, alors que si \m{g=1} il en existe une famille
param\'etr\'ee par $\C$.
\end{subsub}
\end{sub}

\sepsub

\Ssect{Courbes scind\'ees de multiplicit\'e 3}{mult_3_4}

Elles sont enti\`erement classifi\'ees d'apr\`es ce qui pr\'ec\`ede, mais on
peut en donner une autre description. Soit \m{L\in\Pic(C)}. On consid\`ere une
extension
\[0\lra\ko_C\lra E\lra L^*\lra 0\]
sur $C$ correspondant \`a \m{\sigma\in\Ext^1_{\ko_C}(L^*,\ko_C)=H^1(L)}. De
l'inclusion \m{\ko_C\subset E} on d\'eduit un plongement \m{C\subset\P(E)}. On
note \m{C_3(\sigma)} la courbe de multiplicit\'e 3 associ\'ee, qui est
scind\'ee et dont le fibr\'e en droites associ\'e est $L$.

Soit \m{C_2^0} la courbe triviale de multiplicit\'e 2, de courbe r\'eduite $C$
et de fibr\'e en droites associ\'e $L$. Soit \m{\kc^S(C_2^0)} l'ensemble des
classes d'isomorphisme de prolongements de \m{C_2^0} en courbe scind\'ee de
multiplicit\'e $3$. D'apr\`es \ref{c_scind} on a une application surjective
canonique
\[l_g : H^1(L)\lra\kc^S(C_2^0) .\]

\sepprop

\begin{subsub}\label{mult_3_4_1}{\bf Proposition : } 1 - Les fibres de \m{l_g}
sont les orbites de l'action de \m{\C^*} par multiplication.

2 - Pour tout \m{\sigma\in H^1(L)}, on a \ \m{l_g(\sigma)=C_3(\sigma)}.
\end{subsub}

\begin{proof}
La premi\`ere assertion d\'ecoule imm\'ediatement de la classification
compl\`ete des courbes de multiplicit\'e $3$. La seconde se d\'emontre
ais\'ement en partant d'un cocycle d\'efinissant $\sigma$, d'o\`u on d\'eduit
des trivialisations locales de $E$ permettant de d\'efinir \m{C_3(\sigma)} par
un cocycle ad\'equat.
\end{proof}

\end{sub}

\end{document}